\documentclass[11pt,a4paper]{amsart}
\usepackage{amsthm,amsmath} 
\usepackage{amssymb} 
\usepackage{mathrsfs} 
\usepackage{color}

\RequirePackage[numbers]{natbib}
\RequirePackage[colorlinks,citecolor=blue,urlcolor=blue]{hyperref}
\RequirePackage{graphicx}

\numberwithin{equation}{section}
\theoremstyle{plain}

\newtheorem{thm}{Theorem}[section]
\newtheorem{lem}[thm]{Lemma}
\newtheorem{prop}[thm]{Proposition}

\theoremstyle{remark}

\newtheorem{rem}[thm]{Remark}

\newtheorem{defn}[thm]{Definition}
\newtheorem{assum}[thm]{Assumption}

\newcommand{\la}{\langle}
\newcommand{\ra}{\rangle}

\newcommand{\ba}{\begin{array}}
\newcommand{\ea}{\end{array}}
\newcommand{\be}{\begin{equation}}
\newcommand{\ee}{\end{equation}}
\newcommand{\bea}{\begin{eqnarray}}
\newcommand{\eea}{\end{eqnarray}}
\newcommand{\beaa}{\begin{eqnarray*}}
\newcommand{\eeaa}{\end{eqnarray*}}

\def\dbE{\mathbb{E}}
\def\dbF{\mathbb{F}}

\def\dbL{\mathbb{L}}

\def\dbP{\mathbb{P}}
\def\dbR{\mathbb{R}}

\def\sL{\mathscr{L}}

\def\a{\alpha}
\def\b{\beta}

\def\d{\delta}
\def\e{\varepsilon}

\def\l{\lambda}

\def\si{\sigma}

\def\f{\varphi}
\def\th{\theta}
\def\o{\omega}
\def\h{\widehat}
\def\D{\Delta}
\def\Th{\Theta}
\def\L{\Lambda}

\def\O{\Omega}
\def\cA{{\mathcal A}}

\def\cC{{\mathcal C}}

\def\cF{{\mathcal F}}
\def\cG{{\mathcal G}}

\def\cL{{\mathcal L}}

\def\cN{{\mathcal N}}

\def\cP{{\mathcal P}}

\def\cU{{\mathcal U}}

\def\no{\noindent}

\def\ss{\smallskip}
\def\ms{\medskip}
\def\bs{\bigskip}
\def\q{\quad}
\def\qq{\qquad}

\def\pa{\partial}
\def\cd{\cdot}
\def\cds{\cdots}

\def\td{\nabla}

\def\tr{\hbox{\rm tr}}

\newcommand{\basa}{\begin{assumption}}
\newcommand{\easa}{\end{assumption}}

\newcommand{\bas}{\begin{assum}}
\newcommand{\eas}{\end{assum}}

\def\limsup{\mathop{\overline{\rm lim}}}

\def\pa{\partial}
\def\h{\widehat}

 \def\cd{\cdot}
\def\cds{\cdots}

\def\tr{\hbox{\rm tr$\,$}}

\def\dis{\displaystyle}

\def\wh{\widehat}

\def\1{{\bf 1}}

\def\:{\!:\!}
\def\reff{\eqref}

at 9pt

\definecolor{alp}{rgb}{0.0, 0.5, 0.0}

\def\R{\mathbb{R}}

\renewcommand {\theequation}{\arabic{section}.\arabic{equation}}

\def\thesection{\arabic{section}}

\begin{document}


\title[Master Equations with Non-separable Hamiltonians]{Mean Field Games Master Equations with Non-separable Hamiltonians and Displacement Monotonicity} 

\author{Wilfrid Gangbo}
\address{Department of Mathematics, UCLA, California, USA}
\email{wgangbo@math.ucla.edu}

\author{Alp\'ar R. M\'esz\'aros}
\address{Department of Mathematical Sciences, University of Durham, Durham, UK}
\email{alpar.r.meszaros@durham.ac.uk}

\author{Chenchen Mou}
\address{Department of Mathematics, City University of Hong Kong, Hong Kong SAR, China}
\email{chencmou@cityu.edu.hk}

\author{Jianfeng Zhang}
\address{Department of Mathematics, USC, California, USA}
\email{jianfenz@usc.edu}


\makeatletter
\@namedef{subjclassname@2020}{%
  \textup{2020} Mathematics Subject Classification}
\makeatother

\begin{abstract} 
In this manuscript, we propose a structural condition on  non-separable Hamiltonians, which we term {\it displacement monotonicity} condition, to study second order mean field games master equations.  A rate of dissipation of a bilinear form is brought to bear a global (in time) well-posedness theory, based on a priori uniform Lipschitz estimates on the solution in the measure variable. Displacement monotonicity being sometimes in dichotomy with the widely used Lasry-Lions monotonicity condition,  the novelties of this work persist even when restricted to separable Hamiltonians.
\end{abstract}

\date{}
\maketitle

%



\section{Introduction}\label{sect-Introduction}
\setcounter{equation}{0} 

In this manuscript, $T>0$ is a given arbitrary time horizon and $\beta\geq0.$ We consider evolutive equations which represent  games where the players are in motion in the space $\mathbb R^d$ and their distributions at each time are represented by elements of  $\cP_2(\R^d)$, the set of Borel probability measures on $\R^d$, with finite second moments.
The data governing the game are a  Hamiltonian $H$ and a terminal cost function $G$ such that 
$$
H:\R^d\times\cP_2(\R^d)\times\R^d\to\R\quad \text{and} \quad G:\R^d\times\cP_2(\R^d)\to\R.
$$ 
For our description, we assume to be given a rich enough underlying probability measure space $(\Omega, \mathcal F, \mathbb P)$.  The problem at hand is to find a real valued function $V$ which depends on the time variable $t$, the space variable $x$ and the probability measure variable $\mu,$ such that 
\begin{equation}\label{master}
\left\{
\begin{array}{ll}
-\pa_t V -\frac{\h\b^2}{2} \tr(\pa_{xx} V) + H(x,\mu,\partial_x V)  - \cN V =0, & {\rm{in}}\ (0,T)\times\R^d\times\cP_2(\R^d),\\
\dis V(T,x,\mu) = G(x,\mu), & {\rm{in}}\ \R^d\times\cP_2(\R^d).
\end{array}
\right.
\end{equation} 
This second order equation is called the master equation in mean field games, in presence of both idiosyncratic and common noise (if $\beta>0$), where  $\cN$ is the non--local operator defined by 
\begin{align}\label{eq:M-operator}
\left.\ba{c}
\dis \cN V(t,x,\mu):= \tr\bigg(\bar{\tilde  \dbE}\Big[\frac{\h\b^2}{2} \pa_{\tilde x} \pa_\mu V(t,x, \mu, \tilde \xi) +\frac{\b^2}{2}\pa_{\mu\mu}V(t,x,\mu,\bar\xi,\tilde\xi) 
 \\
 +\b^2 \pa_x\pa_\mu V(t,x,\mu,\tilde \xi) - \pa_\mu V(t, x, \mu, \tilde \xi)(\pa_pH)^\top(\tilde \xi,\mu, \pa_x V(t, \tilde \xi, \mu))\Big]\bigg).
 \ea\right.
\end{align}
Above, $\b$ stands for the intensity of the common noise, the idiosyncratic noise is supposed to be non-degenerate (for simplicity, its intensity is set to be 1) and we use the notation $\h \b^2 := 1+\b^2$.  We always assume $H(x, \mu, \cdot)$ to be convex, however, we emphasize already at this point the fact that in general it can have a `non-separable structure', i.e. we {\it do not} assume to have a decomposition of the form  
\bea
\label{separable}
H(x, \mu, p) = H_0(x, p) - F(x, \mu).
\eea
In \eqref{master}-\eqref{eq:M-operator}, $\partial_t$ stands for the time derivative while $\partial_x$ stands for the gradient operator on $\mathbb R^d.$ We postpone to Section \ref{sect-Preliminary}, comments on the $W_2$--Wasserstein gradient $\partial_{\mu}$ and the $W_2$--Wasserstein second gradient $\partial_{\mu\mu}$. 
Given $\mu\in\cP_2(\R^d)$, $\tilde \xi$ and $\bar\xi$ are independent random variables with the same law $\mu$, and  $\bar{\tilde \dbE}$ is the expectation with respect to their joint law.  

First introduced by Lions in lectures \cite{Lions}, the master equation appeared in the context of the theory of mean field games, a theory initiated independently by Lasry-Lions \cite{LL06a, LL06b, LL07a} and Caines-Huang-Malham\'e \cite{HCM06}. It is a time dependent equation which serves to describe the interaction between an individual agent and a continuum of other agents. The master equation characterizes the equilibrium cost of a representative agent within a continuum of players, provided there is a unique mean field equilibrium. Roughly speaking, it plays the role of the Hamilton-Jacobi-Bellman equation in the stochastic control theory. We refer the reader to  
\cite{Cardaliaguet,CD1,CD2}  for a comprehensive exposition on the  subject.

The master equation \reff{master} is known to admit a local (in time) classical solution when the data $H$ and $G$ are sufficiently smooth, even when the noises are absent (cf. \cite{BY,GangboS2015,Mayorga}). Local solutions are known 
to be unique 
(cf.  \cite{CCP,CD2}),  including cases where the Hamiltonians are local functions of the measure variable (cf. \cite{AmbMes}). Nevertheless, it is much more challenging to obtain global classical solutions, as they are expected to exist only under additional structural assumptions on the data. Such a sufficient condition is typically a sort of {\it monotonicity condition},  that provides uniqueness of solutions to the underlying mean field game system (a phenomenon that heuristically corresponds to the non-crossing of generalized characteristics of the master equation). For a non-exhaustible  
list of results on the global in time well-posedness theory of mean field games master equations in various settings, we refer the reader to 
\cite{CDLL,CD1,CD2, CCD}, and in the realm of potential mean field games, to  \cite{BGY1,BGY2,GM}. We also refer to  \cite{Bertucci2, CarSou, MZ2} for global existence and uniqueness of weak solutions 
and to  \cite{BC,BCCD1,BCCD2,Bertucci,BerLasLio,BLL} for finite state mean field games master equations. All the above global well-posedness results require the Hamiltonian $H$ to be separable in $\mu$ and $p$, i.e.  it is of the form \eqref{separable}, for some $H_0$ and $F$.
Moreover, as highlighted above, $F$ and $G$ need to satisfy a certain monotonicity condition, which in particular ensures the uniqueness of mean field equilibria of the corresponding mean field games. We remark that non-separable Hamiltonians appear naturally in applications (such as economical models, problems involving congestions effects, etc., see e.g. \cite{AHLLM,AP,Ambrose,GV}). We shall also note that \cite{BLPR} establishes the global in time well-posedness result for a linear master equation, without requiring separability or monotonicity conditions. However, since the Hamiltonian $H$ is linear in $p$, there is no underlying game involved in \cite{BLPR}.

A typical condition, extensively used in the literature \cite{BC, BLL, CDLL, CD1, CD2, CCD, MZ2}, is the so-called {\it Lasry-Lions monotonicity} condition. For a function $G:\R^d\times\cP_2(\R^d)\to\R$, this can be formulated as
\bea
\label{LLmon}
\mathbb E\Big[G(\xi_1,\mathcal{L}_{\xi_1})+G(\xi_2,\mathcal{L}_{\xi_2})-G(\xi_1,\cL_{\xi_2})-G(\xi_2,\mathcal{L}_{\xi_1})\Big]\ge 0,
\eea
for any random variables $\xi_1, \xi_2$ with appropriate integrability assumptions. Here, $\cL_\xi:=\xi_\#\mathbb{P}$ stands for the law of the random variable $\xi$.

In this manuscript, we turn to a different condition. The main condition we impose here on $G,$ is what we term the {\it displacement monotonicity} condition, which can be formulated as
\bea
\label{dismon}
\mathbb E\Big[\big[\pa_xG(\xi_1,\mathcal{L}_{\xi_1})-\pa_xG(\xi_2,\mathcal{L}_{\xi_2})\big][\xi_1-\xi_2]\Big]\ge 0.
\eea
When $G$ is sufficiently smooth, displacement monotonicity means that  the bilinear form 
\begin{equation}\label{eq:dec25.2020.1}
(\eta_1, \eta_2) \mapsto (d_x d)_\xi G(\eta_1, \eta_2):=  \tilde \dbE \big[\big \langle \partial_{x \mu} G(\xi, \mu, \tilde \xi) \tilde \eta_1, \eta_2\big \rangle\big]+
\mathbb E\big[ \langle \partial_{xx}G(\xi, \mu)\eta_1, \eta_2 \rangle\big]
\end{equation}
is non--negative definite for all square integrable random variables $\xi$. Here $(\tilde \xi, \tilde \eta_1)$ is an independent copy of $(\xi,\eta_1)$ and $\tilde{\mathbb E}$ is the expectation with respect to the joint law of $(\xi,\eta_1,\eta_2,\tilde\xi,\tilde{\eta_1})$. 

Our terminology is inspired by  the so--called {\it displacement convexity} condition, a popular notion in the theory of optimal transport theory (cf. \cite{McCann}). Indeed, when $G$ is derived from a potential, i.e. there exists $g:\cP_2(\R^d)\to\R$ such that $\partial_x G=\partial_\mu g$, then \eqref{dismon} is equivalent to the displacement convexity of $g$. Let us underline that in the current study, we never  need to require that $G$ is derived from a potential.  

Displacement convexity and monotonicity have some sparse history in the framework of mean field games and control problems of McKean-Vlasov type. In the context of mean field game systems, the first work using this seems to be the one of Ahuja \cite{Ahuja} (see also \cite{ARY}), whose weak monotonicity condition is essentially equivalent to the displacement monotonicity. In the context of control problems of McKean-Vlasov type, displacement convexity assumptions appeared first in  \cite{CD:15} and \cite{CCD}. It seems  that  \cite{CCD} is the first work that relied on displacement convexity in the study of well-posedness of a master equation arising in a McKean-Vlasov control problem. However, let us emphasize that this master equation is different from the master equation appearing in the theory of mean field games,
and the techniques developed in \cite{CCD} are not applicable in our setting. In the framework of potential master equations and in particular in more classical infinite dimensional control problems  on Hilbert spaces,  the displacement convexity condition has been used in \cite{BGY1,BGY2}\footnote{These references essentially used a notion of  $\lambda$-convexity in displacement sense, and obtained local in time classical solutions for the master equation. However, it is clear from their results that the solution is global when the data are actually displacement convex.} and \cite{GM}.

Our main contribution  in this manuscript is the discovery of a condition on $H$, which allows a global well--posedness theory of classical solution for the master equation \reff{master}. This condition, which we continue to term {\it displacement monotonicity condition for Hamiltonians}, amounts to impose that the bilinear form
\begin{align*}
 (\eta_1, \eta_2)&\mapsto ({ {\rm{displ}}}^{\f}_\xi H)(\eta_1, \eta_2):=(d_x d)_\xi H\big(\cdot,{\f(\xi)}\big) (\eta_1,  \eta_2) \\
& +{1\over 4} \tilde \dbE\Big[\Big\langle \big(\pa_{pp} H(\xi, \mu, \f(\xi))\big)^{-1} \pa_{p\mu} H(\xi, \mu, \tilde \xi, \f(\xi)) \tilde \eta_1,  ~\pa_{p\mu} H(\xi, \mu, \tilde \xi, \f(\xi)) \tilde \eta_2\Big\rangle\Big]
\end{align*}
is non--positive definite for all $\mu \in \cP_2$, $\xi\in \dbL^2(\cF, \mu)$ and   
all appropriate  $\f\in C^1(\dbR^d; \dbR^d)$, see  Definition \ref{defn:displace-mono-H} below for the precise condition.  In the previous formula, we clearly assume strict convexity on $H$ the $p$ variable. This condition is instrumental for our global well--posedness theory of classical solution to the master equation \reff{master}. To the best of our knowledge, this is the first global well-posedness result in the literature of mean field games master equations, for non-separable Hamiltonians. When $H$ is separable (i.e. of the form \eqref{separable}) and $H_0 = H_0(p)$, the non--positive definiteness assumption on ${ {\rm{displ}}}^{\f}_\xi H$, is equivalent to \reff{dismon} for $F$. For certain Hamiltonians, displacement monotonicity is in dichotomy with the Lasry-Lions monotonicity. Thus, not only our well-posedness  results are new  for a wide class of data functions, but we shall soon see that the novelty in our results, extends to a class of separable Hamiltonians. For discussions on displacement monotone functions that fail to be Lasry-Lions monotone we refer to \cite{Ahuja, GM} and to Subsection \ref{subsec:LL} below.

We show, at the heart of our analysis, that under the displacement monotonicity condition on $H$ and \eqref{dismon} on $G$, $V(t,\cdot,\cdot),$ the solution to the master equation \reff{master}, which has sufficient a priori regularity, also satisfies \reff{dismon} for all $t\in[0,T]$. Let us recall that, when $H$ is separable and both $G$ and $F$ satisfy the Lasry-Lions monotonicity condition \reff{LLmon}, then $V(t,\cdot,\cdot)$ inherits \reff{LLmon} as well for all $t\in[0,T]$. However, when $H$ is non-separable, it remains a challenge to find an appropriate condition on $H$ which ensures that if $G$ satisfies the Lasry-Lions monotonicity, so  does  $V(t,\cdot,\cdot)$ for all $t\in[0,T]$ (see Remark \ref{rem-Vconvex}(iii) below). 

For separable $H$, the Lasry-Lions monotonicity of $V(t,\cdot,\cdot)$ is typically proven through the mean field game system, the corresponding coupled system of forward backward (stochastic) PDEs or SDEs for which $V$ serves as the decoupling field (see Remark \ref{rem-MFS} below). We instead follow a different route and derive the displacement monotonicity of $V(t,\cdot,\cdot)$ by using the master equation itself. We show that if $V$ is a smooth solution to the master equation then for any $\mu \in \cP_2$ and $\eta \in \mathbb L^2(\mathcal F_0)$ there  exists a  path $t \mapsto (X_t, \delta X_t)$ of random variables starting at $(\xi, \eta)$, with $\mu=\cL_\xi$, such that 
%
\begin{equation}\label{eq:dec24.2020.2}
 (d_x d)_{X_T}V_T\big(\delta X_T, {\delta X}_T\big)-  \int_0^T  \big({\rm{displ}}_{X_t}^{\f_t} H\big) \big(\delta X_t, {\delta X}_t\big) dt \leq (d_x d)_{X_0} V_0(\eta,  \eta).
\end{equation}
Here, (see Remark \ref{rem-Vconvex} for a more accurate formulation), 
\[
\f_t=\partial_x V(t, \cdot, \mu_t), \quad \mu_t=X_{t\, \#} \mathbb P.
\]
Note that \eqref{eq:dec24.2020.2} provides us an explicit ``rate of dissipation of displacement monotonicity'' of the bilinear form $(d_x d)V(t, \cdot, \cdot)$, from smaller to larger times. This favors our terminal value problem, as we are provided with a ``rate at which the displacement monotonicity is built in'' from larger to smaller times. 

Our approach seems new, even when restricted to separable $H.$ We are also able to obtain a variant of \eqref{eq:dec24.2020.2} that is applicable to the Lasry-Lions monotonicity case, but only for separable $H.$ One trade--off in our approach is that, since we apply It\^{o}'s formula on the derivatives of $V$, we need higher order a priori regularity estimates on $V$ and consequently require regularity of the data slightly higher than what is needed for the existence of local classical solutions (cf. \cite{CD2}). We 
believe that, thanks to the smooth mollification technique developed in \cite{MZ2}, one could relax these regularity requirements. In fact, we expect a well-posedness theory of weak solutions in the sense of \cite{MZ2}. In this work, our main goal is to overcome the challenge of dealing with non-separable Hamiltonians and so, this manuscript postpones the optimal regularity issue to future studies. 

The displacement monotonicity of $V(t,\cdot,\cdot)$  has a noticeable implication: it yields an a priori uniform $W_2$--Lipschitz continuity estimate  for $V$ in the $\mu$ variable.  Here is the main principle to  emphasize: any possible alternative condition to the non--positive definiteness assumption on ${ {\rm{displ}}}^{\f}_\xi H$, which ensures the monotonity of $V$ (either in Lasry-Lions sense or in displacement sense), will also provide the uniform Lipschitz continuity of $V$ in $\mu$ (with respect to either  $W_1$ or $W_2$). As a consequence, this yields the global well-posedness of the master equation. We shall next elaborate on this observation which seems to be new in the literature and interesting on its own right.

Uniform $W_1$--Lipschitz continuity of $V$ is known to be the key ingredient for constructing  even local in time classical solutions of the master equation (cf. \cite{CD2,MZ2}) in mean field games with common noise. The uniform $W_2$--Lipschitz property we obtain, is not final. We complement this in light of a crucial observation: when the data $H$ and $G$ are uniformly $W_1$--Lipschitz continuous in $\mu$, we can show that the uniform  $W_2$--Lipschitz continuity of $V$ actually implies its uniform $W_1$--Lipschitz continuity in the $\mu$ variable. We achieve this by a delicate analysis on the pointwise representation formula for $\pa_\mu V$, developed in \cite{MZ2}, tailored to our setting.


In our final step to establish the global well-posedness of the master equation, we follow the by now standard approach in  \cite{CD2,CCD,MZ2}. That is, based on the a priori uniform Lipschitz continuity property of $V$ in the $\mu$ variable (with respect to $W_1$), we construct the local classical solution and then  extend it backwardly in time.  Another important point in our argument is that the length of the time intervals used for the local solutions, depends only on the $W_2$--Lipschitz constants of the data.

The rest of the paper is organized as follows. Section \ref{sect-Preliminary} contains the setting of our problem and some preliminary results. In Section \ref{sect-assum} we present our technical assumptions and introduce the new notion of displacement monotonicity for non-separable $H$. 
In Section \ref{sect-Vconvex} we show that any solution  of the master equation which is regular enough, preserves the displacement monotonicity property. Section \ref{sect-Lipschitz} is devoted to uniform a priori $W_2$--Lipschitz estimates on $V$. In Section \ref{sect-global} we derive the uniform $W_1$--Lipschitz estimates and establish the global well-posedness of the master equation \reff{master}. 

\medskip

\section{Preliminaries} 
\label{sect-Preliminary}
\setcounter{equation}{0}

\subsection{The product probability space}
\label{sect-space}
In this paper we shall use the probabilistic approach. In order to reach out to the largest community of people working on mean field games, in this subsection we present our probabilistic setting in details, which we think will facilitate the reading of those who are not experts in stochastic analysis.

Throughout the paper, we fix $T>0$ to be a given arbitrary time horizon. Let $(\O_0, \dbF^0, \dbP_0)$ and $(\O_1, \dbF^1, \dbP_1)$ be two filtered probability spaces, on which there are defined $d$-dimensional Brownian motions $B^0$ and $B$, respectively. For $\dbF^i =\{\cF^i_t\}_{0\le t\le T}$, $i=0,1$, we assume $\cF^0_t=\cF^{B^0}_t$,  $\cF^1_t =\cF^1_0 \vee \cF^{B}_t$, and $\dbP_1$ has no atom in $\cF^1_0$ so it can support any measure on $\dbR^d$ with finite second moment. Consider the product spaces 
\begin{equation}\label{product}
\O := \O_0 \times \O_1,\q \dbF = \{\cF_t\}_{0\le t\le T} := \{\cF^0_t \otimes \cF^1_t\}_{0\le t\le T},\q \dbP := \dbP_0\otimes \dbP_1,\q \dbE:= \dbE^\dbP.
\end{equation}
In particular, $\cF_t := \si(A_0\times A_1: A_0\in \cF^0_t, A_1\in \cF^1_t\}$ and $\dbP(A_0\times A_1) = \dbP_0(A_0) \dbP_1(A_1)$. 
 We shall automatically extend $B^0, B, \dbF^0, \dbF^1$ to the product space in the obvious sense, but using the same notation. For example,  $B^0(\o) = B^0(\o^0)$ for $\o=(\o^0, \o^1)\in \O$, and $\cF^0_t = \{A_0 \times \O_1: A_0\in \cF^0_t\}$. In particular, this implies that $B^0$ and $B^1$ are independent $\dbP$-Brownian motions and are independent of $\cF_0$. 

It is convenient to introduce another filtered probability space $(\tilde \O_1, \tilde \dbF^1, \tilde B, \tilde \dbP_1)$  in the same manner as $(\O_1,  \dbF^1, B, \dbP_1)$, and consider the larger filtered probability space given by
\begin{equation}\label{product2}
\tilde \O := \O\times \tilde\O_1 ,\q \tilde\dbF = \{\tilde \cF_t\}_{0\le t\le T} := \{\cF_t \otimes \tilde \cF^1_t\}_{0\le t\le T},\q \tilde \dbP := \dbP\otimes \tilde\dbP_1,\q \tilde \dbE:= \dbE^{\tilde \dbP}.
\end{equation}
Given an $\cF_t$-measurable random variable $\xi = \xi(\o^0, \o^1)$, we say $\tilde \xi = \tilde \xi(\o^0, \tilde \o^1)$ is a conditionally independent copy of $\xi$ if, for each $\o^0$, the $\dbP_1$-distribution of $\xi(\o^0, \cd)$ is equal to the $\tilde\dbP_1$-distribution of $\tilde\xi(\o^0, \cd)$.  That is, conditional on $\cF^0_t$, by extending to $\tilde \O$ the random variables $\xi$ and $\tilde \xi$ are conditionally independent and have the same conditional distribution under $\tilde \dbP$. Note that, for any appropriate deterministic function $\f$,
\bea
\label{conditional expectation}
\left.\ba{c}
\dis \tilde \dbE_{\cF^0_t} \big[ \f(\xi, \tilde \xi)\big] (\o^0) = \dbE^{\dbP_1\otimes \tilde\dbP_1}\Big[\f\big(\xi(\o^0, \cdot), \tilde\xi(\o^0, \tilde\cdot)\big)\Big],\q \dbP_0-\mbox{a.e.}~\o^0;\\
\dis \tilde \dbE_{\cF_t} \big[ \f(\xi, \tilde \xi)\big] (\o^0,\o^1) = \dbE^{\tilde\dbP_1}\Big[\f\big(\xi(\o^0, \o^1), \tilde\xi(\o^0, \tilde \cdot)\big)\Big],\q \dbP-\mbox{a.e.}~(\o^0, \o^1).
\ea\right.
\eea
Here $ \dbE^{\tilde\dbP_1}$ is the expectation on $\tilde \o^1$, and $\dbE^{\dbP_1\times \tilde\dbP_1}$ is on $(\o^1, \tilde \o^1)$. Throughout the paper, we will use the probability space $(\O, \dbF, \dbP)$. However, when conditionally independent copies of random variables or processes are needed, we will tacitly use the extension to the larger space $(\tilde \O, \tilde \dbF, \tilde \dbP)$ without mentioning.

When we need two conditionally independent copies, we introduce further $(\bar \O_1, \bar \dbF^1, \bar B, \bar \dbP_1)$ and the product space $(\bar\O, \bar \dbF, \bar\dbP, \bar\dbE)$  as in \reff{product2}, and set the joint product space
\begin{align}\label{product3}
&\bar{\tilde \O}  := \O\times \tilde\O_1\times \bar\O_1 ,\q \bar{\tilde \dbF} = \{\bar{\tilde \cF}_t\}_{0\le t\le T} := \{\cF_t\otimes\tilde\cF^1_t \otimes \bar \cF^1_t\}_{0\le t\le T},\\ 
\nonumber&\bar {\tilde \dbP} := \dbP\otimes \tilde\dbP_1\otimes \bar\dbP_1,\q \bar{\tilde \dbE}:= \dbE^{\bar{\tilde \dbP}}.
\end{align}
Then, given $\cF_t$-measurable $\xi=\xi(\o^0, \o^1)$,  we may have two conditionally independent copies  under $\bar {\tilde \dbP}$: $\tilde \xi(\bar{\tilde \o}) = \tilde\xi(\o^0, \tilde \o^1)$ and $\bar \xi(\bar{\tilde \o}) =\bar \xi(\o^0, \bar \o^1)$,  $\bar{\tilde \o}=(\o^0, \o^1, \tilde \o^1, \bar \o^1)\in \bar{\tilde \O}$.  

To avoid possible notation confusion, we emphasize that
\begin{align}
&\mbox{when $\xi = \xi(\o^1)$ is $\cF^1_t$-measurable, then $\tilde\xi, \bar\xi$ are independent copies of $\xi$ under $\bar {\tilde \dbP}$};\nonumber\\
\label{product4}
&\mbox{the expectation $\tilde \dbE$ is on $\tilde \o = (\o^0, \o^1, \tilde \o^1)$, not just on $\tilde \o^1$; similarly $\bar \dbE$  is acting on}\\ 
\nonumber&\mbox{$\bar \o = (\o^0, \o^1, \bar \o^1)$, not just on $\bar \o^1,$}\ \mbox{and $\bar{\tilde \dbE}$ is an expectation on $\bar{\tilde \o} = (\o^0, \o^1, \tilde \o^1, \bar \o^1)$}.\nonumber
\end{align}


\subsection{Preliminary analysis on the Wasserstein space} 
\label{sect:Wasserstein}
Let $\cP:=\cP(\dbR^d)$ be the set of all probability measures on $\mathbb R^d$ and  $\d_x\in \cP$ denotes the Dirac mass at $x\in \mathbb R^d$. For any $ q\geq 1$ and any measure $\mu\in \cP$, we set 
\begin{equation}\label{cP}
M_q(\mu):=\left( \int_{\mathbb R^d} |x|^q \mu(dx)\right)^{1\over q}\q\text{and}\q \mathcal{P}_q:=\mathcal{P}_q(\mathbb R^d):=\{\mu\in\mathcal{P}:~ M_q(\mu) < \infty\}.
\end{equation}
%
%
%
For any sub-$\si$-field $\cG\subset \cF_T$ and $\mu\in \cP_q$, denote by $\dbL^q(\cG)$ the set of $\dbR^d$-valued, $\cG$-measurable, and $q$-integrable random variables $\xi$;  and $\dbL^q(\cG;\mu)$ the set of $\xi\in \dbL^q(\cG)$ such that $\cL_\xi=\mu$. Here $\cL_\xi=\xi_\# \dbP$ is the law of $\xi$,   obtained as the push--forward of $\dbP$ by $\xi$. Also, for $\mu\in \cP_q$, let $\dbL^q_\mu(\dbR^d; \dbR^d)$ denote the set of Borel measurable functions $v:\dbR^d\to \dbR^d$ such that $\|v\|_{\dbL_{\mu}^q}^q:=\int_{\dbR^d} |v(x)|^q \mu(dx) < \infty$.
%
%
Moreover, for any $\mu,\nu\in \cP_q$,  their $W_q$--Wasserstein distance  is defined as follows: 
\bea
\label{Wp}
W_q(\mu, \nu) := \inf\Big\{\big(\dbE[|\xi-\eta|^q]\big)^{1\over q}:~\mbox{for all $\xi\in \dbL^q(\cF_T; \mu)$, $\eta\in \dbL^q(\cF_T; \nu)$}\Big\}.
\eea

According to the terminology in \cite{AGS}, the Wasserstein gradient of a function $U: \cP_2 \to \dbR$ at $\mu$, is an element $\pa_\mu U(\mu,\cdot)$ of $\overline{\nabla C_c^\infty(\mathbb R^d)}^{\dbL^2_\mu}$ (the closure of gradients of $C_c^\infty$ functions in $\dbL^2_\mu(\R^d;\R^d)$) and so, it is a priori defined $\mu$--almost everywhere. The theory developed in \cite{Cardaliaguet, GT, Lions, WZ} 
shows that $\pa_\mu U(\mu,\cdot)$ can be characterized by the property 
\bea
\label{pamu}
U(\cL_{\xi +  \eta}) - U(\mu) = \dbE\big[\langle \pa_\mu U(\mu, \xi), \eta \rangle \big] + o(\|\eta\|_2), \ \forall\ \xi,\eta,\ {\rm{with}}\ \cL_\xi=\mu.
\eea
Let $\cC^0(\cP_2)$ denote the set of $W_2$--continuous functions $U: \cP_2 \to \dbR$. For $k \in \{1,2\}$ we next define a subset of $\cC^k(\cP_2)$, referred to  as functions of {\it full $\cC^k$} regularity in \cite[Chapter 5]{CD1}),  as follows. By $\cC^1(\cP_2),$ we mean the space of functions $U\in C^0(\cP_2)$ such that $\pa_\mu U$ exists for all $\mu\in\cP_2$ and it has a unique jointly continuous extension to $\cP_2\times \dbR^d$, which we continue to denote by 
\vskip-0.15in
$$ \mathbb R^d\times\cP_2\ni (\tilde x, \mu) \mapsto \pa_\mu U (\mu, \tilde x)\in \mathbb R^d.$$ 
\vskip-0.03in
\no We sometimes refer to the extension as the global version, and we note that our requirement of pointwise continuity property of this global version is stronger than the $\dbL^2$--continuity requirement made in some of the mean field game literature (cf. e.g. \cite{CCD}). Similarly, $\cC^2(\cP_2)$ stands for the set of functions $U\in \cC^1(\mathcal{P}_2)$ such that the global version of $\partial_\mu U$ is differentiable in the sense that the following maps exist and have unique jointly continuous extensions: 
\begin{align*}
&\mathbb R^d\times\cP_2\ni (\tilde x, \mu) \mapsto  \pa_{\tilde x\mu}U(\mu, \tilde x) \in \mathbb R^d\ \ {\rm{and}}\\ 
&\mathbb R^{2d}\times\cP_2\ni (\tilde x, \bar x,  \mu) \mapsto  \pa_{\mu\mu}U(\mu, \tilde x, \bar x) \in \mathbb R^{d \times d}.
\end{align*}
\noindent $\cC^{2}(\mathbb R^d\times\mathcal{P}_2)$ is the set of continuous functions $U:\mathbb R^d\times\cP_2\to\dbR$ satisfying the following: 

(i) $\pa_xU, \pa_{xx}U$ exist and are jointly continuous on $\R^d\times\cP_2$;

(ii) The following maps exist and have unique jointly continuous extensions
\begin{align*}
&\mathbb R^{2d}\times\cP_2\ni (x, \tilde x, \mu) \mapsto \pa_\mu U(x, \mu, \tilde x)\in \mathbb R^d\ \ {\rm{and}}\\
&\mathbb R^{2d}\times\cP_2\ni (x, \tilde x, \mu) \mapsto\pa_{x\mu}U(x, \mu, \tilde x) \in \mathbb R^{d \times d};
\end{align*}

(iii) Finally, the following maps exist and have unique jointly continuous extensions 
\begin{align*}
&\mathbb R^{2d}\times\cP_2\ni (x, \tilde x, \mu) \mapsto \pa_{\tilde x\mu}U(x, \mu, \tilde x) \in \mathbb R^{d \times d} \ \ \text{and}\\
&\mathbb R^{3d}\times\cP_2\ni (x, \tilde x, \bar x, \mu) \mapsto \pa_{\mu\mu}U(x, \mu, \tilde x, \bar x) \in \mathbb R^{d \times d}.
\end{align*}
Lastly, we fix the state space for our master equation as
\[
\Th:= [0, T]\times \dbR^d \times \cP_2, 
\]
and let $\cC^{1,2,2}(\Th)$ denote the set of  $U\in C^0( \Th; \dbR)$ such that the following maps exist and have a unique jointly continuous extensions as previously described: 
$\pa_t U$, $\pa_x U$, $\pa_{xx} U$, $\pa_\mu U$, $\pa_{x}  \pa_\mu U$, $\pa_{\tilde x}  \pa_\mu U,$ $ \pa_{\mu\mu} U.$

We underline that for notational conventions, we always denote the `new spacial variables' appearing in Wasserstein derivatives with tilde symbols (for first order Wasserstein derivatives), with ``bar'' symbols (for second order Wasserstein derivatives) and so on, and we place them right after the corresponding measures variables. For example, when $U:\R^d\times\cP_2\times\R^d\to\R$ is typically evaluated as $U(x,\mu,p)$, we use the notations $\partial_\mu U(x,\mu,\tilde x,p)$, $\partial_{\tilde x}\partial_\mu U(x,\mu,\tilde x, p)$, $\partial_\mu\partial_\mu U(x,\mu,\tilde x,\bar x, p)$, and so on. This convention will be carried through to compositions with random variables too, for example $\partial_\mu U(x,\mu,\tilde\xi,p)$, when $\tilde\xi$ is an $\R^d$-valued random variable.

Throughout the paper, we shall also use the following notations: for any $R>0$,
\begin{align}\label{BR}
B_R^o :=  \big\{p\in \dbR^d: |p|< R\big\},\q B_R := \big\{p\in \dbR^d: |p|\le R\big\},\q D_R:= \mathbb R^d\times \mathcal{P}_2\times B_R.
\end{align}

The following simple technical lemma (not to confuse with  \cite[Remark 4.16]{CD2}) is useful.
\begin{lem}\label{lem-order} 
For any $U\in \mathcal C^2(\mathcal{P}_2)$ and $(\mu,\tilde x)\in\mathcal{P}_2\times\mathbb R^d$, $\partial_{\tilde x \mu} U(\mu,\tilde x)$ is a symmetric matrix.
\end{lem}

\begin{proof} 
Since $U$ is of class $\cC^2(\cP_2)$, we may assume without loss of generality that $\mu$ is supported by a closed ball $B_R$, it is absolutely continuous and has a smooth density $\rho$ with $c:= \inf_{x\in B_R} \rho(x)>0$.
By the fact that $\partial_\mu U(\mu, \cdot) \in \overline{\nabla C_c^\infty(\mathbb R^d)}^{\dbL^2_\mu}$, there exists a sequence $(\varphi_n)_n \subset \nabla C_c^\infty(\mathbb R^d)$ such that 
\begin{equation}\label{1}
0=\lim_{n}\|\partial_{\tilde x} \varphi_n-\partial_\mu U(\mu, \cdot)\|_{\dbL^2_\mu} \geq c \lim_{n}\|\partial_{\tilde x} \varphi_n-\partial_\mu U(\mu, \cdot)\|_{L^2({B_R})},
\end{equation}
where 
$L^2({B_R})$ stands for the standard Lebesgue space. Set 
\[
\bar \varphi_n:= \varphi_n-{1 \over \mathcal L^d(B_R)}\int_{B_R}  \varphi_n(x) dx.
\]
By the Poincar\'e--Wirtinger inequality, there exists a universal constant $c_d$ such that 
\[
\|\bar \varphi_n\|_{L^2(B_R)} \leq c_d \|\partial_{\tilde x} \bar \varphi_n\|_{L^2({B_R})}.
\]
Thanks to the Sobolev Embedding Theorem and the strong convergence of $(\partial_{\tilde x} \bar \varphi_n)_n$ in $L^2(B_R)$, we conclude that there exists $\varphi$ in the Sobolev space $H^1(B_R)$ such that $( \bar \varphi_n)_n$ converges to $\varphi$ in $H^1(B_R)$. By \eqref{1} we have 
\begin{equation}\label{2}
\partial_\mu U(\mu, \cdot)=\partial_{\tilde x} \varphi, \qquad \int_{B_R}  \varphi(x) dx=0.
\end{equation}
Since $\partial_\mu U(\mu, \cdot)$ is continuously differentiable, the representation formula 
\[
\varphi({\tilde x})= \varphi(0)+ \int_0^1 \partial_\mu U(\mu, t {\tilde x}) \cdot {\tilde x} dt
\]
implies that $\varphi$ is continuously differentiable. Since $\partial_{\tilde x} \varphi= \partial_\mu U(\mu, \cdot)$ is continuously differentiable, we conclude that $\varphi$ is twice continuously differentiable. Thus, $\partial_{\tilde x}\partial_\mu U(\mu, \cdot)=\partial_{\tilde x \tilde x} \varphi$ is symmetric. 
\end{proof}

An interesting property of  $\cC^{1,2,2}(\Th)$ functions is their use in a general It\^{o}  formula. Let $U\in C^{1,2,2}(\Th)$ be such that for any compact subset $K\subset \mathbb R^d\times\mathcal{P}_2$
\beaa
&&\sup_{(t,x,\mu)\in [0,T]\times K}\Big[\int_{\mathbb R^d}\Big(|\pa_\mu U(t,x,\mu,\tilde x)|^2+|\pa_{\tilde x}\pa_{\mu}U(t,x,\mu,\tilde x)|^2
+|\pa_{x}\pa_{\mu}U(t,x,\mu,\tilde x)|^2\Big)\mu(d\tilde x)\\
&&\qquad \qquad \qquad  +\int_{\mathbb R^{2d}}|\pa_{\mu\mu}U(t,x,\mu,\tilde x,\bar x)|^2\mu(\tilde x)\mu(\bar x)\Big]<+\infty.
\eeaa
For $i=1,2$,  consider $\dbF$-progressively measurable and bounded processes
$$b^i: [0, T]\times \O\to \dbR^d \quad \text{and} \quad \si^i, \si^{i,0}: [0, T]\times \O\to \dbR^{d\times d}.$$
Set
$$d X^i_t := b^i_t dt + \si^i_t dB_t + \si^{i,0}_t d B^0_t,\q\mbox{and introduce the conditional law}~\rho_t:=\mathcal{L}_{X^2_t|\cF^{0}_t}.
$$ 
Then (cf., e.g.,  \cite[Theorem 4.17]{CD2}, \cite{BLPR,CCD}), recalling the notations for conditionally independent copies and \reff{product4},
\small
\begin{align}\label{Ito}
d U(t, X^1_t, &\rho_t) =  \Big[\pa_t U + \pa_x U\cdot b^1_t + \frac{1}{2} \tr\big(\pa_{xx} U [\si_t^1 (\si_t^1)^\top + \si_t^{1,0}(\si_t^{1,0})^\top]\big)\Big](t, X^1_t, \rho_t) dt\nonumber\\
&+\pa_xU(t,X^1_t,\rho_t)\cd\si_t^1dB_t +\Big[(\si^{1,0}_t)^\top\pa_xU+\tilde{\mathbb E}_{\cF_t}\big[(\tilde \si^{2,0}_t)^\top\pa_\mu U(\cd,\tilde X^2_t) \big]\Big](t,X^1_t,\rho_t)\cd dB_t^0 
\nonumber\\
&+ \tr \bigg(\bar{\tilde \dbE}_{\cF_t}\Big[\pa_\mu U(\cd,\tilde X^2_t) (\tilde b^{2}_t)^\top+\frac{1}{2} \pa_{\tilde x}\pa_\mu U(\cd, \tilde X^2_t)\big[\tilde \si_t^2 (\tilde \si_t^2)^\top + \tilde \si_t^{2,0}(\tilde \si_t^{2,0})^\top\big]\\
&\q +\pa_x\pa_\mu U(\cd,\tilde X^2_t)\si^{1,0}_t (\tilde \si_t^{2,0})^\top+\frac{1}{2}\pa_{\mu\mu}U(\cd,\tilde X^2_t,\bar X^2_t) \tilde\si_t^{2,0}(\bar \si_t^{2,0})^\top\Big](t,X^1_t,\rho_t)\bigg) dt.\nonumber
\end{align}\normalsize
Throughout this paper, the elements of $\dbR^d$ are viewed as column vectors; $\pa_x U, \pa_{\mu} U\in\R^d$ are also column vectors; $\pa_{x\mu}U:= \pa_x \pa_\mu U := \pa_x \big[(\pa_\mu U)^\top\big]\in \dbR^{d\times d}$, where $^\top$ denotes the transpose, and similarly for the other second order derivatives; both the notations ``$\cd$'' and $\langle\cdot,\cdot\rangle$ denote the inner product of column vectors. Moreover, the term $\pa_x U \cd \si^1_t dB_t$ means $\pa_x U \cd (\si^1_t dB_t)$, but we omit the parentheses for notational simplicity. 

\subsection{The Lasry-Lions monotonicity and the displacement monotonicity}\label{subsec:LL}
In this subsection, we discuss two types of monotonicity conditions and provide more convenient alternative formulations.
\begin{defn} Let $U :\dbR^d\times \cP_2\to\mathbb R$.
\begin{enumerate}
\item[(i)] $U$ is called \emph{Lasry-Lions monotone}, if for any $\xi_1, \xi_2\in \dbL^2(\mathcal{F}^1_T)$,
\bea
\label{mon1}
\mathbb E\Big[U(\xi_1,\mathcal{L}_{\xi_1})+U(\xi_2,\mathcal{L}_{\xi_2})-U(\xi_1,\mathcal{L}_{\xi_2})-U(\xi_2,\mathcal{L}_{\xi_1})\Big] \ge 0.
\eea
\item[(ii)] $U$ is called \emph{displacement monotone} if  $U(\cdot,\mu)\in C^1(\R^d)$ for all $\mu\in\cP_2$ and for any $\xi_1, \xi_2\in \dbL^2(\mathcal{F}^1_T)$,
\bea
\label{displacement1}
\dbE\Big[ \big\langle\pa_x U(\xi_1, \mathcal{L}_{\xi_1}) - \pa_x U(\xi_2, \mathcal{L}_{\xi_2}), \xi_1-\xi_2\big\rangle\Big] \ge 0.
\eea
\end{enumerate}
\end{defn}
\begin{rem}\label{rem:implicationDisp} Assume $U \in \cC^1(\dbR^d\times \cP_2)$
\begin{enumerate}
\item[(i)] If $\partial_\mu U(\cdot,\mu,\tilde x)\in C^1(\R^d)$, for all $(\mu,\tilde x)\in\cP_2\times\R^d$, then the inequality \eqref{mon1} implies,
\beaa
0&\le& \dbE\Big[U(\xi, \cL_\xi) + U(\xi+ \e \eta, \cL_{\xi+ \e\eta}) - U(\xi, \cL_{\xi+\e\eta}) - U(\xi + \e\eta, \cL_\xi)\Big]\\
&=&\e^2 \int_0^1 \int_0^1 \tilde\dbE\Big[\big\langle \pa_{x\mu} U\big(\xi+\th_1 \e \eta, \cL_{\xi+\th_2 \e\eta}, \tilde \xi+\th_2\e \tilde \eta\big)\tilde \eta, \eta\big\rangle \Big] d\th_1 d\th_2,
\eeaa
for any $\xi, \eta\in \dbL^2(\mathcal{F}^1_T)$ and any $\e>0$, where $(\tilde \xi, \tilde \eta)$ is an independent copy of $(\xi, \eta)$. Thus
\bea
\label{mon2} 
\tilde\dbE\Big[\big\langle \pa_{x\mu} U\big(\xi, \cL_\xi, \tilde \xi\big)\tilde \eta, \eta\big\rangle\Big] \ge 0,\q\forall \xi, \eta\in \dbL^2(\cF^1_T).
\eea
\item[(ii)] If $\partial_x U\in \mathcal C^1(\dbR^d\times \cP_2)$, then the inequality \eqref{displacement1} implies 
\beaa
&&\dis 0\le \dbE\Big[ \langle\pa_x U(\xi+\e \eta, \cL_{\xi+\e\eta}) - \pa_x U(\xi, \cL_\xi), \e\eta\rangle \Big]\\
&&\dis= \e^2 \int_0^1 \tilde\dbE\Big[ \big\langle\pa_{xx} U(\xi+\th \e \eta, \cL_{\xi+\th\e\eta}) \eta, \eta\big\rangle +\big \langle \pa_{x\mu} U(\xi+\th \e \eta, \cL_{\xi+\th\e\eta},\tilde \xi+\th \e\tilde \eta)\tilde \eta, \eta \big\rangle \Big] d\th,
\eeaa
for any $\xi, \eta\in \dbL^2(\mathcal{F}^1_T)$ and  $\e>0$, where $(\tilde \xi, \tilde \eta)$ is an independent copy of $(\xi, \eta)$, and thus, 
\begin{align}\label{displacement2} 
(d_x d)_\xi U(\eta, \eta):=
\tilde\dbE\Big[\big\langle \pa_{x\mu} U\big(\xi, \cL_\xi, \tilde \xi\big)\tilde \eta, \eta\big\rangle  +\big\langle \pa_{xx} U\big(\xi, \cL_\xi\big)\eta, \eta\big\rangle  \Big] \ge 0. 
\end{align} 
\item[(iii)] Assume $\cU\in\cC^2(\cP_2)$ and $U\in \cC^1(\R^d\times\cP_2)$ are such that $\partial_\mu\cU \equiv \partial_x U(x,\mu)$ on $\R^d\times\cP_2.$ 
Then $U$ is displacement monotone if and only if $\cU$ is \emph{displacement convex},  cf. \cite{McCann}. 
\end{enumerate}
\end{rem}

\begin{rem}\label{rem-monequiv}  Throughout this manuscript, given $U\in \cC^2(\dbR^d\times \cP_2)$, we call \reff{mon2} the Lasry-Lions monotonicity condition and call \reff{displacement2}  the displacement monotonicity condition.  Indeed, it is obvious that \reff{displacement1} and \reff{displacement2} are equivalent. We prove in the appendix that \reff{mon1} and \reff{mon2} are also equivalent. 
\end{rem}

\begin{rem}\label{rem-displacement}
(i) \reff{displacement2} implies that $U$ is convex in $x$, namely $\pa_{xx} U$ is nonnegative definite. We provide a simple proof in Lemma \ref{lem:Uxxconvex} below, and we refer to \cite[Proposition B.6]{GM} for a more general result.
Note that in particular, \reff{mon2} does not imply \reff{displacement2}. Indeed, let $U(x,\mu) = U_0(x) + U_1(\mu)$ such that $\partial_{xx}U_0$ is not nonnegative definite. Then $\pa_{x\mu} U (x, \mu, \tilde x) \equiv 0$ and so, \reff{mon2} holds while  $\pa_{xx} U = \pa_{xx} U_0$ is not nonnegative definite.  Thus \reff{displacement2} fails.

(ii) For any function $U\in \cC^2(\dbR^d\times \cP_2)$ with $|\pa_{xx} U|$ and $|\pa_{x\mu} U|$ bounded above by $C>0$, the function $\bar U(x,\mu):= U(x,\mu)+C|x|^2$ will always satisfy \reff{displacement2}:
\beaa
&\dis 
(d_x d)_\xi \bar U(\eta, \eta)=\tilde \dbE\Big[\big\langle \pa_{x\mu} U(\xi, \mu, \tilde \xi)\tilde \eta, \eta\big\rangle  +\big\langle \pa_{xx} U(\xi,\mu) \eta,\eta\big\rangle + 2C|\eta|^2\Big]\\
&\dis\ge \tilde\dbE\Big[-C |\eta| |\tilde \eta| - C |\eta|^2 + 2C|\eta|^2\Big] = C\Big[\dbE[|\eta|^2] - |\dbE[\eta]|^2\Big]\ge 0.
\eeaa
This means that \reff{displacement2}  does not imply \reff{mon2} either. Indeed,  if $U$ is a function violating \reff{mon2} but having bounded derivatives, then the above $\bar U$ satisfies \reff{displacement2}. But, since $\pa_{x\mu} \bar U = \pa_{x\mu} U$, $\bar U$ violates  \reff{mon2}.

(iii) We note that, for the function $\bar U$ above, $\pa_x \bar U$ is unbounded. For our main results later, we need displacement monotone functions with bounded derivatives. One can construct such an example as follows.
Let $\phi: \dbR^d \to \dbR$ be convex, even and smooth with bounded derivatives. 
Set $U(x, \mu) := \int_{\dbR^d} \phi(x-y) \mu(dy)$. Then $U$ satisfies \reff{displacement2} and its derivatives are bounded.
\end{rem}

\begin{lem}\label{lem:Uxxconvex}
Assume $\partial_x U\in \mathcal C^1(\dbR^d\times \cP_2)$ and $U$ satisfies \reff{displacement2}. Then $\pa_{xx} U$ is non-negative definite.
\end{lem}

\begin{proof} 
Without loss of generality we assume that $\mu$ has a positive and smooth density $\rho$. For $\xi \in \dbL^2(\cF^1_T, \mu)$, $x_0\in \dbR^d$, and $\eta_\e = v_\e(\xi)$, where 
$v \in C_c^\infty(\mathbb R^d; \dbR^d)$ and for $\e>0$, denote $v_\e(x) := \e^{-d}v({x - x_0 \over \e})$.
We see that $\tilde \eta_\e = v_\e(\tilde \xi)$. Then, straightforward calculation reveals 
\begin{align*}
 (d_xd)_\xi U(\eta_\e, \eta_\e)
&=\int_{\dbR^{2d}} \big\langle \partial_{x\mu}U(x_0+\e z, \mu,  x_0+\e \tilde z) v(\tilde z), v(z) \big\rangle \rho(x_0+\e z)\rho(x_0+\e\tilde z) dzd\tilde z\\
&+ 
\e^{-d}\int_{\dbR^d}\big\langle \partial_{xx}U(x_0+\e z,\mu) v(z), v(z) \big\rangle \rho(x_0+\e z) dz.
\end{align*}
Thus, by \reff{displacement2} we have
\beaa
0 \le \lim_{\e\to 0} \Big[\e^d (d_xd)_\xi U(\eta_\e, \eta_\e)\Big] = \rho(x_0)\int_{\dbR^d}\big\langle \partial_{xx}U(x_0,\mu) v(z), v(z) \big\rangle  dz.
\eeaa
Since $\rho(x_0)>0$ and $v$ is arbitrary, this implies immediately that $\partial_{xx}U(x_0,\mu)$ is non-negative definite.
\end{proof}


In Section \ref{sect-Lipschitz} below, we will also use the following notion of displacement semi-monotonicity, inspired by the displacement semi-convexity in potential games (cf. \cite{AGS, BGY2}). 
\begin{defn} 
Assume  $U,\pa_x U \in \cC^1(\dbR^d\times \cP_2)$. We say $U$ is \emph{displacement semi-monotone} if there exists a constant $\l\ge 0$ such that, for any $\xi, \eta\in \dbL^2(\mathcal{F}^1_T)$,
\bea
\label{semidisplacement1}
(d_x d)_\xi U(\eta, \eta) \ge -\lambda \mathbb E\big[|\eta|^2\big].
\eea
\end{defn}

\begin{rem}
\label{rem-semi}
It is obvious that displacement semi-monotonicity is weaker than the displacement monotonicity. Moreover,  when $\pa_{xx} U$ is bounded, the Lasry-Lions monotonicity \reff{mon1} also implies the displacement semi-monotonicity.
\end{rem}

\subsection{The master equation and mean field games}
In this subsection we summarize in an informal and elementary way, the well-known connection between the solutions of the master equation \reff{master} and the value functions arising in mean field games (cf. e.g. \cite{CD1,CD2}). We recall $\b\geq 0$ represents the intensity of the common noise and  $L, G$ are two given functions:
$$
 L:\dbR^d\times\cP_2\times\dbR^d\to \dbR, \quad \text{and} \quad G: \dbR^d\times \cP_2\to \dbR$$ that are continuous in all variables. As usual, the Legendre-Fenchel transform of  the Lagrangian $L$ with respect to the last variable is  the Hamiltonian $H$  defined  as
\bea
\label{H}
H(x,\mu, p) := \sup_{a\in \dbR^d} [  -\langle a, p\rangle  - L(x,\mu, a) ],\q (x, p, \mu) \in \dbR^{2d} \times \cP_2.
\eea

Given $t\in [0, T]$, we set 
$$
B^t_s:= B_s-B_t, \quad B^{0,t}_s:=B_s^{0}-B_t^0,\qquad \forall s\in [t, T],
$$ and denote by $\cA_t$ the set of admissible  controls $\a: [t, T]\times\mathbb R^d\times C([t, T]; \dbR^{d}) \to \dbR^d$ that are uniformly Lipschitz continuous in the second variable, progressively measurable, and adapted.  For any $\xi\in \dbL^2(\cF_t)$ and $\a\in  \cA_t$, by the Lipschitz continuity property of $\alpha$, the SDE 
\bea
\label{Xth}
\dis X^{t,\xi, \a}_s = \xi + \int_t^s  \a_r( X^{t,\xi, \a}_r, B^{0, t}_\cd) dr + B^t_s+\b B^{0,t}_s,\q s\in [t, T],
\eea
has a unique strong solution. We note that, by the adaptedness, the control $\a$ actually takes the form  $\a_r( X^{t,\xi, \a}_{r}, B^{0, t}_{[t, r]})$, where, $B^{0, t}_{[t, r]}$ stands for the restriction of $B^{0, t}$ to the interval $[t, r]$.  Consider the conditionally expected cost functional for the mean field game:
\bea
\label{Jxi}
\left.\ba{lll}
\dis J(t,x,\xi; \a,  \a'):= \dbE^{\dbP}_{\cF^0_t}\Big[G(X^{t, x, \a'}_T,  \cL_{X^{t,\xi,\a}_T|\cF^0_T}) \\
\dis\qq\qq + \int_t^T   L(X_s^{t,x,\alpha'},\cL_{X^{t,\xi,\a}_s|\cF^0_s},\a'_s (X^{t,x,\a'}_\cd, B^{0, t}_\cd))  ds\Big].
\ea\right.
\eea
 Here, $\xi$ represents the initial state of the ``other" players, $\a$ is the common control of the other players, and $(x, \a')$ is the initial state and control of the individual player. When $\xi\in \dbL^2(\cF^1_t)$ is independent of $\cF^0_t$,  it is clear that $J(t,x,\xi; \a,  \a')$ is deterministic. One shows that 
\[
\xi'\in \dbL^2(\cF^1_t),\;\; \cL_{\xi'}=\cL_\xi \quad \implies \quad J(t,x,\xi'; \a,  \a')=J(t,x,\xi; \a,  \a') \qquad \forall x, \a,\a'.
\]
Therefore,  we may define  
\bea
\label{Jmu}
J(t,x,\mu; \a,  \a') := J(t,x,\xi; \a,  \a'),\quad \q \xi\in \dbL^2(\cF^1_t, \mu).
\eea
Now for any $ (t,x,\mu)\in \Th$ and $\a\in \cA_t$,  we consider the infimum 
\bea
\label{Va0}
\dis V( t,x,\mu;\a) := \inf_{\a'\in \cA_t}  J(t,x,\mu; \a, \a').
\eea

\begin{defn}
\label{defn-NE}
We say $\a^*\in \cA_t$ is a mean field Nash equilibrium of \reff{Va0} at $(t,\mu)$ if 
\beaa
V( t,x,\mu; \a^*) = J(t,x,\mu; \a^*, \a^*) \quad \mbox{for $\mu$-a.e. $x\in \dbR^d$}.
\eeaa
\end{defn}
When there is a unique mean field equilibrium  for each $(t,\mu)$, denoted as $\a^*(t,\mu)$, it makes sense to define 
\bea
\label{V}
V(t,x,\mu) := V( t,x,\mu; \a^*(t,\mu)).
\eea 
Using the It\^{o} formula \reff{Ito}, one shows that if $V$ if sufficiently regular, then it is a classical solution to the master equation \reff{master}. However, we would like to point out that the theory of the global well-posedness of \eqref{master} that we develop will not rely explicitly on this connection.


\medskip

The master equation \eqref{master} is also associated to the following forward backward McKean-Vlasov SDEs on $[t_0,T]$: given $t_0$ and $\xi\in \dbL^2(\cF_{t_0})$,
\begin{equation}\label{FBSDE1} 
\dis\left\{\ba{ll}
\dis X^{\xi}_t & \dis =\xi - \int_{t_0}^t\pa_pH(X_s^\xi,\rho_s,Z_s^{\xi})ds+B^{t_0}_t+\b B_t^{0,t_0}; \\
\dis Y_t^{\xi}&\dis =G(X_T^{\xi},\rho_T)+\int_t^T\wh L(X_s^{\xi},\rho_s,Z_s^{\xi})ds- \int_t^TZ_s^{\xi}\cd dB_s-\int_t^TZ_s^{0,\xi}\cd dB_s^{0},
\ea\right.
\end{equation}
where
$$\q \wh L(x,\mu, p) := L(x, \mu,\pa_p H(x,\mu,p))=p\cdot \pa_pH(x,\mu,p)-H(x,\mu,p),\q  \rho_t := \rho^\xi_t:= \cL_{X_t^\xi|\cF^0_t}.$$
Given $\rho$ as above  and $x\in\R^d$, we consider on $[t_0,T]$,  the standard decoupled FBSDE
\small\begin{equation}\label{FBSDE2}
\left\{\ba{ll}
\dis X_t^{x}&\dis=x+ B^{t_0}_t+\b B_t^{0,t_0}; \\
\dis Y_t^{x,\xi}&\dis=G(X_T^{x},\rho_T)- \int_t^T H(X_s^x,\rho_s,Z_s^{x,\xi})ds-\int_t^T Z_s^{x,\xi}\cd dB_s-\int_t^T Z_s^{0,x,\xi}\cd dB_s^{0}.
\ea\right.
\end{equation}
\normalsize
Alternatively, we may consider the coupled FBSDE instead of the decoupled one \reff{FBSDE2}:
\small\begin{equation}\label{FBSDE3}
\left\{
\begin{array}{l}
 \dis X_t^{\xi,x}=x-\int_{t_0}^t\pa_pH(X_s^{\xi,x},\rho_s,Z_s^{\xi,x})ds+ B^{t_0}_t+\b B_t^{0,t_0}; \\
\dis Y_t^{\xi,x}=G(X_T^{\xi,x},\rho_T)+\int_t^T \wh L(X_s^{\xi,x},\rho_s,Z_s^{\xi,x})ds-\int_t^T Z_s^{\xi,x}\cd dB_s-\int_t^T Z_s^{0,\xi,x}\cd dB_s^{0}.
\end{array}
\right.
\end{equation}
\normalsize
These  FBSDEs connect to the master equation \reff{master} as follows: if $V$ is a classical solution to \eqref{master} and if the above FBSDEs have strong solution, then 
\bea\label{YXV}
\left.\ba{c}
Y^\xi_t = V(t, X^\xi_t, \rho_t),\q Y^{x,\xi}_t = V(t, X^{x}_t, \rho_t),\q Y^{\xi,x}_t = V(t, X^{\xi,x}_t, \rho_t),\\
Z^\xi_t = \pa_x V(t, X^\xi_t, \rho_t),\q Z^{ x,\xi}_t = \pa_xV(t, X^{x}_t, \rho_t),\q Z^{\xi,x}_t = \pa_xV(t, X^{\xi,x}_t, \rho_t).
\ea\right.
\eea

\begin{rem}
\label{rem-MFS}
 (i) The forward-backward SDE system \reff{FBSDE1}-\reff{FBSDE2} or \reff{FBSDE1}-\reff{FBSDE3} is called the mean field game system. Equivalently, one may also consider the following forward-backward stochastic PDE system as the mean field system on $[t_0, T]$:
\small
\begin{equation}
\label{SPDE}
\left\{
\begin{array}{ll}
\dis d \rho(t,x)&\dis =  \Big[\frac{\h\beta^2}{2}\tr\big( \pa_{xx} \rho(t,x)\big) + div(\rho(t,x) \pa_p H(x, \rho(t,\cd), \pa_x u(t,x)))\Big]dt-\b\partial_x\rho(t,x)\cd d B_t^0\\
\dis d u(t, x)&\dis =  v(t,x)\cd dB_t^0 - \Big[\tr\big(\frac{\h\beta^2}{2} \pa_{xx} u(t,x) + \b\partial_x v^\top(t,x)\big) -H(x,\rho(t,\cdot),\pa_x u(t,x))\Big]dt\\
\dis&\dis \rho(t_0,\cd) = \cL_\xi,\q  u(T,x) = G(x, \rho(T,\cdot)).
\end{array}
\right.
\end{equation}\normalsize
Here the solution triple $(\rho, u, v)$ is $\dbF^0$-progressively measurable and $\rho(t,\cd,\o)$ is a (random) probability measure. The solution $V$ to the master equation also serves as the decoupling field for this forward-backward system, i.e.
\bea
\label{uvV}
u(t,x,\o) = V(t,x, \rho(t,\cd,\o)).
\eea

(ii) In this paper we focus on the well-posedness of the master equation \eqref{master}. It is now a folklore in the literature that once we obtain a classical solution $V$ (with suitably bounded derivatives), we immediately get existence and uniqueness of a mean field equilibrium $\a^*$ in \eqref{Va0} in the sense of Definition \ref{defn-NE}. Indeed, given $V$, in light of \reff{YXV} we may decouple the forward backward system \reff{FBSDE1} (or similarly decouple \reff{SPDE}) as
\begin{equation}
\label{FBSDE4}
X^{\xi}_t=\xi - \int_{t_0}^t\pa_pH\big(X_s^\xi,\rho_s, \pa_x V(s, X^\xi_s, \rho_s)\big)ds+B^{t_0}_t+\b B_t^{0,t_0},\q  \rho_t := \cL_{X_t^\xi|\cF^0_t}.
\end{equation}
If $V$ is sufficiently regular, this SDE has a unique solution $(X^\xi, \rho)$, and then we can easily see that 
\beaa
\a^*(t, x, \o):= -\pa_pH\big(x,\rho_t(\o), \pa_x V(t, x, \rho_t(\o)\big)
\eeaa
is the unique mean field equilibrium of the game.

(iii) Given a classical solution $V$ with bounded derivatives, in particular with bounded $\pa_{\mu\mu} V$, we can show the convergence of the corresponding $N$-player game. The arguments are more or less standard, see \cite{CDLL, CD2}, and we leave the details to interested readers. 
\end{rem}

\section{The displacement monotonicity of non-separable $H$} 
\label{sect-assum}
\setcounter{equation}{0}
In this section we collect all our standing assumptions on the data that are used in this manuscript to prove our main theorems. In particular, we shall introduce our new notion of displacement monotonicity for non-separable $H$. Under appropriate condition on $H$ and recalling \reff{eq:M-operator} for the non-local operator $\cN$, it is convenient in the sequel to define the operator 
\[ 
 \sL V(t,x,\mu):=-\pa_t V -\frac{\h\b^2}{2} \tr(\pa_{xx} V) + H(x,\mu,\partial_x V)  - \cN V,
\]
which acts on the set of smooth functions on $[0,T] \times \mathbb R^d\times \mathcal{P}_2.$ 

We first specify the technical conditions on $G$ and $H$. Recall the $B_R$ and $D_R$ in \reff{BR}. 

\begin{assum}\label{assum-regG} We make the following  assumptions on $G$.

\ms
(i) $G\in \cC^2(\mathbb R^d\times\mathcal{P}_2)$ with $|\pa_x G|,|\pa_{xx}G|\le L^G_0$ and $|\pa_\mu G|,|\pa_{x\mu}G| \le L^G_1$;


\ms
(ii) 
$
G,\pa_xG, \pa_{xx}G\in \cC^2(\mathbb R^d\times\mathcal{P}_2)$, and $\pa_\mu G, \pa_{x\mu}G \in C^2(\mathbb R^d\times\mathcal{P}_2\times\mathbb R^d),
$
and the supremum norms of all their derivatives are uniformly bounded.
\end{assum}

\begin{assum} \label{assum-regH} We make the following  assumptions on $H$.

(i) $H\in \cC^2(\mathbb R^d\times\mathcal{P}_2\times\mathbb R^d)$ and, for any $R>0$,  there exists $L^H(R)$ such that 
\beaa
&|\pa_xH|,|\pa_pH|, |\pa_{xx}H|,|\pa_{xp}H|,|\pa_{pp}H|\le L^H(R),\ \ {\rm{on}}\ D_R;\ms\\
&|\pa_\mu H|,|\pa_{x\mu}H|,|\pa_{p\mu}H|\le L^H(R),\ \ {\rm{on}}\ \mathbb R^d\times\mathcal{P}_2\times\mathbb R^d\times B_R;
\eeaa

(ii) $H\in \cC^3(\mathbb R^d\times\mathcal{P}_2\times\mathbb R^d)$, and
$$ 
H,\; \pa_xH,\;  \pa_pH,\; \pa_{xx}H,\;  \pa_{xp}H,\;  \pa_{pp}H,\;  \pa_{xxp}H,\;  \pa_{xpp}H,\;  \pa_{ppp}H\; \in \cC^2(\mathbb R^d\times\mathcal{P}_2\times\mathbb R^d), 
$$ 
the supremum norms of all their derivatives are uniformly bounded  on $D_R$ and  
$$
\pa_\mu H,\pa_{x\mu}H, \pa_{p\mu}H,\pa_{ xp\mu}H,\pa_{ pp\mu}H\in \cC^2(\mathbb R^d\times\mathcal{P}_2\times\mathbb R^{2d})
$$ and the supremum norms of all their derivatives are bounded  on $\mathbb R^d\times\mathcal{P}_2\times\mathbb R^d\times B_R;$

(iii) There exists $C_0>0$, such that 
$$|\pa_x H(x,\mu,p)|\leq C_0(1+|p|),\ \  {\rm{for\ any\ }} (x,\mu,p)\in\mathbb R^d\times\mathcal{P}_2\times\mathbb R^d;$$

(iv) $H$ is strictly convex in $p$, and, for any $R>0$, there exists $L^H(R)$ such that 
\[
\Big|\Big(\pa_{pp} H(x, \mu, p)\Big)^{-{1\over 2}}\partial_{p\mu}H(x,\mu,\tilde{x},p)\Big|\leq L^H(R),\ \  {\rm{for\ any\ }} (x,\mu,\tilde x,p)\in\mathbb R^d\times\mathcal{P}_2\times\mathbb R^d\times B_R.
\] 
\end{assum}
\begin{rem}
\label{rem:W2}
 (i) Given a function $U\in \cC^1(\cP_2)$, one can easily see that $U$ is uniformly $W_1$--Lipschitz continuous   if and only if $\pa_\mu U$ is bounded. 
 

(ii) Under Assumption \ref{assum-regG} and  by the above remark, we see that $G$ and $\pa_x G$ are uniformly Lipschitz continuous in $\mu$ under $W_1$ on $\mathbb R^d\times\mathcal{P}_2$ with  Lipschitz constant $L^G_1$. This implies further the Lipchitz continuity of $G, \pa_x G$ in $\mu$ under $W_2$  on $\mathbb R^d\times \mathcal{P}_2$, and we denote the Lipschitz constant by $L^G_2 \le L^G_1$:
\bea
\label{L2G}
\tilde \dbE\Big[ |\pa_\mu G(x, \mu, \tilde \xi) \tilde \eta|\Big] \le L^G_2 \Big(\dbE[|\eta|^2]\Big)^{1\over 2},\q \tilde \dbE\Big[ |\pa_{x\mu} G(x, \mu, \tilde \xi) \tilde \eta|\Big] \le L^G_2 \Big(\dbE[|\eta|^2]\Big)^{1\over 2},
\eea
for all $\xi\in \dbL^2(\cF^1_T, \mu)$, $\eta\in \dbL^2(\cF^1_T)$. Similarly, under Assumption \ref{assum-regH}, 
$H, \pa_x H, \pa_pH$ are uniformly Lipschitz continuous in  $\mu$ under $W_1$ (or $W_2$) on  $\mathbb R^d\times \mathcal{P}_2\times B_R$ with Lipschitz constant $L^H(R)$. 
\end{rem}

We now introduce the crucial notion of displacement monotonicity for non-separable $H$.  
\begin{defn} \label{defn:displace-mono-H} Let $H$ be a Hamiltonian satisfying \ref{assum-regH}(i) and (iv). We say that $H$ is \emph{displacement monotone} if  for any $\mu \in \cP_2$, $\xi\in \dbL^2(\cF^1_T, \mu)$, and any bounded Lipschitz continuous function $\f\in C^1(\dbR^d; \dbR^d)$, the following bilinear form is non-positive definite on $\eta\in \dbL^2(\cF^1_T)$:

\bea
\label{Hconvex} 
\left.\ba{lll}
\dis ({ {\rm{displ}}}^{\f}_\xi H)(\eta, \eta):=  (d_x d)^\f_\xi H(\eta, \eta) + Q^{\f}_\xi H(\eta,  \eta) \le 0,\q  {\rm{where}}\\
\dis (d_x d)^\f_\xi H(\eta, \eta):= \tilde \dbE\Big[\big\langle \pa_{x\mu} H(\xi, \mu, \tilde \xi, \f(\xi))\tilde \eta + \pa_{xx} H(\xi, \mu, \f(\xi)) \eta,\; \eta\big\rangle  \Big];\\
\dis Q^{\f}_\xi H(\eta, \eta):={1\over 4}\dbE\Big[\Big| \big(\pa_{pp} H(\xi, \mu, \f(\xi))\big)^{-{1\over 2}} \tilde \dbE_{\cF^1_T}\big[\pa_{p\mu} H(\xi, \mu, \tilde \xi, \f(\xi)) \tilde \eta\big] \Big|^2\Big].
\ea\right.
\eea

\end{defn}
The following assumptions are central in our work.

\begin{assum}\label{assum-convex} 
(i) $G$ satisfies Assumption \ref{assum-regG} (i) and it is displacement monotone, namely it satisfies \reff{displacement2}.

(ii) $H$ satisfies Assumptions \ref{assum-regH}(i), (iv) and is displacement monotone, namely \eqref{Hconvex} holds.
\end{assum}
\begin{rem}
\label{rem-Hconvex}
 (i) When $H(x,\mu, p) = H_0(p)-F(x,\mu) $, \reff{Hconvex} reads off 
\beaa
({ {\rm{displ}}}^{\f}_\xi H)(\eta, \eta) = - (d_x d)_\xi F(\eta, \eta) \le 0.
\eeaa
This is precisely the displacement monotonicity condition \reff{displacement2} on $F$ and so,  \reff{Hconvex} is an extension of the displacement monotonicity to the functions on $\mathbb R^d \times \cP_2\times \dbR^d.$ 

(ii) Under Assumptions \ref{assum-regG} and \ref{assum-regH}, we may weaken the requirement in Definition \ref{defn:displace-mono-H} such that \reff{Hconvex} holds true only for those $\f$ satisfying $|\f|\le C^x_1$, $|\pa_x \f|\le C^x_2$, for the constants $C^x_1, C^x_2$ determined in \reff{pauest1} below. All the results in this paper will remain true under this weaker condition.

(iii) As in Remark \ref{rem-displacement} (i), one can easily see that  \reff{Hconvex} implies $\pa_{xx} H$ is non-positive definite. This will be useful in the proof of Proposition \ref{prop-Hconvex}.
\end{rem}

\begin{prop}
\label{prop-Hconvex}
Under Assumptions \ref{assum-regH}(i) and (iv),  $H$ is displacement monotone if and only if \reff{Hconvex} holds true for $\si(\xi)$-measurable $\eta$, namely $\eta = v(\xi)$ for some deterministic function $v$. That is, by writing in integral form,  $H$ is displacement monotone if and only if, for any $\mu \in \cP_2$,  $v\in \dbL^2_\mu(\dbR^d; \dbR^d)$ (defined in Section \ref{sect:Wasserstein}), and any bounded Lipschitz continuous function $\f\in C^1(\dbR^d; \dbR^d)$, it holds
\bea
\label{Hconvex2}
\left.\ba{c}
\dis\int_{\dbR^{2d}} \Big\langle \pa_{x\mu} H(x, \mu, \tilde x, \f(x)) v(\tilde x) + \pa_{xx} H(x, \mu, \f(x)) v(x),\; v(x)\Big\rangle \mu(dx)\mu(d\tilde x)\\
\dis+ {1\over 4}\int_{\dbR^d}\bigg[\Big| \Big(\pa_{pp} H(x, \mu, \f(x))\Big)^{-{1\over 2}} \int_{\dbR^d}\big[\pa_{p\mu} H(x, \mu, \tilde x, \f(x)) v(\tilde x)\big] \mu(d\tilde x)\Big|^2\bigg]\mu(dx)\le 0.
\ea\right.
\eea
In particular, when $H$ is separable, namely $\pa_{p\mu} H=0$ and hence $Q^{\f}_\xi H(\eta, \eta)=0$, then \reff{Hconvex2} reduces to 
\bea
\label{Hconvex3}
\dis\int_{\dbR^{2d}} \Big\langle \pa_{x\mu} H(x, \mu, \tilde x, \f(x)) v(\tilde x) + \pa_{xx} H(x, \mu, \f(x)) v(x),\; v(x)\Big\rangle \mu(dx)\mu(d\tilde x)\le 0.
\eea
\end{prop}
\proof First assume \reff{Hconvex} holds. For any desired $\mu, v, \f$,  let $\xi\in \dbL^2(\cF^1_T, \mu)$ and $\eta := v(\xi)$.  Note that $\tilde \eta = v(\tilde \xi)$ for the same function $v$. Then  \reff{Hconvex2} is exactly the integral form of \reff{Hconvex}.

We now prove the opposite direction. Assume  \reff{Hconvex2} holds true. Following the same line of arguments as in the proof of  Remark \ref{rem-displacement} (i), one first shows that $\pa_{xx} H$ is non-positive definite.  Now for any $\xi\in \dbL^2(\cF^1_T,\mu)$ and $\eta \in \dbL^2(\cF^1_T)$. Denote $\eta' := \dbE[\eta |\xi]$. Then there exists $v\in \dbL^2_\mu(\dbR^d; \dbR^d)$ such that $\eta' = v(\xi)$. Note that $$\tilde \eta':= \tilde\dbE[\tilde \eta |\cF^1_T, \tilde \xi] = \tilde\dbE[\tilde \eta |\tilde \xi]= v(\tilde \xi)$$ for the same function $v$. Then \reff{Hconvex2} implies that \reff{Hconvex} holds for $(\eta', \tilde \eta')$. Note that, by the independence of $(\tilde \xi, \tilde \eta)$ and $(\xi, \eta)$, we have 
\beaa
&&\dis\tilde \dbE\Big[\big\langle \pa_{x\mu} H(\xi, \mu, \tilde \xi, \f(\xi))\tilde \eta,\; \eta\big\rangle  \Big]=\tilde \dbE\Big[\big\langle \pa_{x\mu} H(\xi, \mu, \tilde \xi, \f(\xi))\tilde \eta',\; \eta'\big\rangle  \Big]\\
&&\dis \tilde \dbE_{\cF^1_T}\Big[\pa_{p\mu} H(\xi, \mu, \tilde \xi, \f(\xi)) \tilde \eta\Big] = \tilde \dbE_{\cF^1_T}\Big[\pa_{p\mu} H(\xi, \mu, \tilde \xi, \f(\xi)) \tilde \eta'\Big]. \\
\eeaa
Since $\pa_{xx}H$ is non-positive definite, we have 
\beaa
&&\dis \dbE\Big[\big\langle  \pa_{xx} H(\xi, \mu, \f(\xi)) \eta,\; \eta\big\rangle  \Big] \le  \dbE\Big[\big\langle  \pa_{xx} H(\xi, \mu, \f(\xi)) \eta',\; \eta'\big\rangle  \Big]. 
\eeaa
We combine all these to obtain 
\beaa
({\rm{displ}}^{\f}_\xi H)(\eta, \eta)\le ({ displ}^{\f}_\xi H)(\eta', \eta')\le 0,
\eeaa
which completes the proof. 
\qed

We next provide an example of non-separable $H$ which satisfies all our assumptions. We first note that,  similar to Remark \ref{rem-displacement} (ii), for any $H\in \cC^2(\dbR^d\times \cP_2\times \dbR^d)$ with bounded second order derivates, the function $H(x,\mu, p) - C |x|^2$ always satisfies \reff{Hconvex} for $C>0$ large enough. However, this function $H(x,\mu, p) - C |x|^2$ fails to be Lipschitz in $x$. We thus modify it as follows.

Let $H_0(x, \mu, p)$ be any smooth function with bounded derivatives up to the appropriate order  so that $H_0$ satisfies Assumption \ref{assum-regH}(i)-(ii). Suppose  for some constant $R_0>0$,
\bea
\label{H0supp}
H_0(x, \mu, p) = 0 ~\mbox{when}~|x|> R_0,\q \mbox{and}\q \pa_\mu H_0(x, \mu, \tilde x, p) = 0 ~\mbox{when}~|\tilde x|> R_0. \eea
A particular example of $H_0$ satisfying both conditions in \eqref{H0supp} is
$$H_0(x, \mu, p) = h\left(x, p, \int_{\dbR^d} f(x, \tilde x, p) \mu(d\tilde x)\right)$$ 
where $f$ and $h$ are smooth, $h(x,p,r)=0$ for $|x|>R_0$ and $\partial_{\tilde x}f(x, \tilde x, p)=0$ for $|\tilde x|\ge R_0$.  Let $\psi_C: \dbR^d\to \dbR$ be a smooth and convex function such that $\psi_C(x) = C|x|^2$ when $|x|\le R_0$ and $\psi_C(x)$ growth linearly when $|x|\ge R_0+1$.  Then we have the following result. 

\begin{lem}\label{lem:exampleH} If $C_0$ is sufficiently large then the Hamiltonian 
\bea
\label{Hconstruct}
H(x, \mu, p) := H_0(x,\mu, p) + C_0|p|^2 - \psi_{C_0}(x).
\eea
 satisfies Assumption \ref{assum-regH} and is displacement monotone.
\end{lem}
\proof It is straightforward to verify Assumption \ref{assum-regH} (i), (ii), (iii), and $H$ also satisfies Assumption \ref{assum-regH}  (iv) when $C_0$ is large enough.  Then it remains to prove \reff{Hconvex}. Let $C>0$ be a bound of $\pa_{x\mu} H_0$, $\pa_{xx} H_0$, $\pa_{p\mu} H_0$, and 
choose $C_0$ such that 
$$
2C_0>3C, \quad \pa_{pp} H_0 + 2C_0I_d \ge  I_d.
$$
We first note that 
\beaa
&\pa_{x\mu} H = \pa_{x\mu} H_0,\q \pa_{p\mu} H = \pa_{p\mu} H_0, \q \pa_{pp} H = \pa_{pp} H_0 + 2C_0I_d \ge  I_d,\\
&\pa_{xx} H(x,\mu, p) = \pa_{xx} H_0(x,\mu, p) - 2C_0 I_d \1_{\{|x|\le R_0\}} - \pa_{xx} \psi_{C_0}(x) \1_{\{|x|>R_0\}}.
\eeaa
By \reff{H0supp} we have 
\beaa
&&\dis ({ {\rm{displ}}}^{\f}_\xi H)(\eta, \eta)=\dbE\Big[ \big\langle\tilde\dbE_{\cF^1_T}\big[\pa_{x\mu} H_0(\xi, \mu, \tilde \xi, \f(\xi))\tilde \eta\big], \eta\big\rangle \\
&&\dis\qq + \1_{\{|\xi|\le R_0\}} \big\langle[\pa_{xx} H_0(\xi, \mu, \f(\xi)) - 2C_0 I_d]\eta,\eta\big\rangle -\1_{\{|\xi|> R_0\}} \big\langle[\pa_{xx} \psi_{C_0}(\xi)\eta,\eta\big\rangle \\
&&\dis\qq +\frac{1}{4}\Big|[2C_0I_d+\pa_{pp} H_0(\xi, \mu, \f(\xi))]^{-{1\over 2}} \tilde \dbE_{\cF^1_T}[\pa_{p\mu} H_0(\xi, \mu, \tilde \xi, \f(\xi)) \tilde \eta\big]\Big|^2\Big].
\eeaa
We use Jensen's inequality, the assumption on $C_0$ and by the convexity of $\psi_{C_0}$ to obtain 
\begin{align*} 
\dis({ {\rm{displ}}}^{\f}_\xi H)(\eta,  \eta) &\le \dbE\Big[ C\1_{\{|\xi|\le R_0\}}|\eta| \tilde \dbE\big[\1_{\{|\tilde \xi|\le R_0\}}|\tilde \eta|\big]+[C-2C_0] \1_{\{|\xi|\le R_0\}} |\eta|^2\\
&\dis + C\big[\tilde E[ \1_{\{|\tilde \xi|\le R_0\}}|\tilde \eta|]\big]^2 - \1_{\{|\xi|> R_0\}} \langle[\pa_{xx} \psi_{C_0}(\xi)\eta,\eta\rangle \Big]\\
&\dis \le [C-2C_0] \dbE\big[\1_{\{|\xi|\le R_0\}} |\eta|^2\big] + 2C \Big(\dbE\big[\1_{\{|\xi|\le R_0\}} |\eta|\big]\Big)^2\\
&\dis- \dbE\Big[ \1_{\{|\xi|> R_0\}} \langle[\pa_{xx} \psi_{C_0}(\xi)\eta,\eta\rangle \Big] \leq 0.
\end{align*}
Thus, $H$ satisfies \reff{Hconvex}. 
\qed

We next  express the displacement monotonicity of $H$ in terms of $L,$ defined through \reff{H}. 
\begin{prop} 
\label{prop-Lconvex}
Let $H$ be such that Assumptions \ref{assum-regH} (i) and (iv) hold. Let $\mu \in \mathcal P_2.$ 

(i) $H$ satisfies \eqref{Hconvex} if and only if  $L$ satisfies the following:
\bea
\label{Lconvex}
\left.\ba{c}
\dis \tilde\dbE\Big[\big\langle \pa_{x\mu} L(\xi, \mu, \tilde \xi, \psi(\xi)) \tilde \eta,\eta\big\rangle +\big\langle\pa_{xx} L(\xi, \mu, \psi(\xi))\eta,\eta\big\rangle \Big]  \\
\dis \ge \dbE\Big[\Big|[\pa_{aa}L(\xi,\mu,\psi(\xi))]^{-{1\over 2}}\big[{1\over 2} \tilde \dbE_{\cF^1_T}\big[\pa_{a\mu} L(\xi, \mu, \tilde \xi, \psi(\xi)) \tilde \eta\big] +\pa_{ax}L(\xi, \mu, \psi(\xi)) \eta \big]\Big|^2\Big],
\ea\right.
\eea
for all $\mu, \xi, \eta, \f$ as in Definition \ref{defn:displace-mono-H} and 
$\psi(x) := -\pa_p H(x, \mu, \f(x))$.  

(ii) A sufficient condition for $L$ to satisfy \reff{Lconvex} and hence for $H$  to satisfy \eqref{Hconvex} is:
\bea
\label{Lconvex3}
\L:= {d^2\over d\e d\d}\dbE\Big[ L \big(\xi + (\e+\d)\eta, \cL_{\xi+\e\eta}, \xi'+(\e+\d)\eta'\big)\Big]\Big|_{(\e,\d)=(0,0)} \ge 0,
\eea
for all $\xi, \xi', \eta, \eta'\in \dbL^2(\cF^1_T)$.
\end{prop} 
\proof (i)  First, standard convex analysis theory ensures regularity properties of $L.$ The optimal argument $a^* = a^*(x,\mu, p)$ satisfies:
\beaa
H(x, \mu, p) =-L(x, \mu, a^*) -  \langle a^*, p\rangle, \q \pa_a L(x, \mu, a^*) +p=0,\q  a^* = -\pa_p H(x, \mu, p).
\eeaa
One can easily derive further the following identities (some of them are well-known in convex analysis)
\begin{align}\label{HL}
\left.\ba{lll}
\dis \pa_x H(x,\mu, p) =- \pa_x L(x, \mu, a^*);\q \pa_\mu H(x,\mu, \tilde x, p) = -\pa_\mu L(x, \mu, \tilde x, a^*);\\
\dis\pa_{aa}L \ge {1\over L^H(R)} I_d\; \text{on} \; D_R\q\mbox{and}\q \pa_{pp} H(x, \mu, p) =  [\pa_{aa}L(x, \mu, a^*)]^{-1}; \\
\dis \pa_{xp} H(x, \mu, p) =   \pa_{xa}L(x, \mu, a^*)[\pa_{aa}L(x, \mu, a^*)]^{-1}; \\
\dis\pa_{xx} H(x,\mu, p) =- \pa_{xx} L(x, \mu, a^*) + \pa_{xp} H(x,\mu, p)\pa_{ax} L(x, \mu, a^*)\\
\dis \qq\qq\qq =\Big[-\pa_{xx} L +  \pa_{xa} L [\pa_{aa}L]^{-1}\pa_{ax} L\Big](x, \mu, a^*);\\
\dis \pa_{x\mu} H(x,\mu, \tilde x, p) =-  \pa_{x\mu} L(x, \mu, \tilde x, a^*) +\pa_{xp} H(x,\mu, p) \pa_{a\mu} L(x, \mu, \tilde x, a^*)\\
\dis \qq\qq\qq =\Big[-\pa_{x\mu} L +  \pa_{xa} L [\pa_{aa}L]^{-1}\pa_{a\mu} L\Big](x, \mu, \tilde x, a^*);\\
\dis \pa_{p\mu} H(x,\mu, \tilde x, p) = \pa_{pp} H(x,\mu, p) \pa_{a\mu} L(x, \mu, \tilde x, a^*)\\
\dis\qq\qq\qq\q = [\pa_{aa}L(x, \mu, a^*)]^{-1}\pa_{a\mu} L(x, \mu, \tilde x, a^*).
\ea\right.
\end{align}

Now let $\f$ be chosen as in Definition \ref{defn:displace-mono-H}  and $\psi(x) := -\pa_p H(x, \mu, \f(x))$,  then we have
\beaa
&&-({\rm{displ}}^{\f}_\xi H)(\eta, \eta)=   \tilde\dbE\Big[ \big\langle\big[\pa_{x\mu} L -  \pa_{xa} L [\pa_{aa}L]^{-1}\pa_{a\mu} L\big](\xi, \mu, \tilde \xi, \psi(\xi)) ~\tilde\eta, \eta\big\rangle \\
&& +\big\langle\big[\pa_{xx} L -  \pa_{xa} L [\pa_{aa}L]^{-1}\pa_{ax} L\big](\xi, \mu, \psi(\xi)) ~\eta,\eta\big\rangle \\ 
&& -\frac{1}{4}\Big|[\pa_{aa}L(\xi,\mu,\psi(\xi))]^{-{1\over 2}}\tilde {\mathbb E}_{\cF^1_T}\big[ \pa_{a\mu } L(\xi, \mu, \tilde \xi, \psi(\xi)) \tilde \eta\big]\Big|^2 \Big]\\
&&=\tilde\dbE\Big[\big\langle \pa_{x\mu} L(\xi, \mu, \tilde \xi, \psi(\xi))\tilde \eta,  \eta\big\rangle +\big\langle\pa_{xx} L(\xi, \mu, \psi(\xi))\eta,\eta\big\rangle   \\
&&\qq -\Big|[\pa_{aa}L(\xi,\mu,\psi(\xi))]^{-{1\over 2}}\big[{1\over 2} \tilde \dbE_{\cF^1_T}\big[\pa_{a\mu} L(\xi, \mu, \tilde \xi, \psi(\xi)) \tilde \eta\big] +\pa_{ax}L(\xi, \mu, \psi(\xi)) \eta \big]\Big|^2\Big].
\eeaa
Then clearly  \eqref{Hconvex} is equivalent to \reff{Lconvex}.

(ii) Assume \eqref{Lconvex3} holds and $\xi, \xi', \eta, \eta'\in \dbL^2(\cF^1_T)$. By straightforward calculations we have
\beaa
\dis \L &=& {d\over d\e}\dbE\Big[ \big\la \pa_x L \big(\xi + \e\eta, \cL_{\xi+\e\eta}, \xi'+\e\eta'\big), \eta\big\ra +  \big\la \pa_a L \big(\xi + \e\eta, \cL_{\xi+\e\eta}, \xi'+\e\eta'\big), \eta'\big\ra \Big]\Big|_{\e=0}\\
&=& \tilde\dbE\Big[ \big\la \pa_{xx} L \big(\xi, \cL_{\xi}, \xi'\big) \eta, \eta\big\ra +\big\la  \pa_{x\mu} L \big(\xi, \cL_{\xi}, \tilde \xi, \xi'\big) \tilde \eta, \eta\big\ra\\
&&+ 2 \big\la \pa_{ax} L \big(\xi, \cL_{\xi}, \xi'\big)\eta, \eta'\big\ra +  \big\la \pa_{a\mu} L \big(\xi, \cL_{\xi}, \tilde \xi, \xi'\big)\tilde \eta, \eta'\big\ra+  \big\la \pa_{aa} L \big(\xi, \cL_{\xi}, \xi'\big)\eta', \eta'\big\ra\Big].
\eeaa
The expression $\L$ remains non-negative in particular when 
$$
\xi':= \psi(\xi),\ \ \text{and} \ \  
\eta' := -\big(\pa_{aa}L(\xi, \cL_\xi, \xi')\big)^{-1}\bigg(\frac{1}{2} \tilde \dbE_{\cF^1_T}[\pa_{a\mu}L(\xi, \cL_\xi, \tilde \xi, \xi')\tilde\eta]+\pa_{ax}L(\xi, \cL_\xi, \xi')\eta\bigg).
$$ 
Omitting the variables  $(\xi, \cL_\xi, \tilde \xi, \xi')$ inside the derivatives of $L$, we have
\beaa
0&\le& \L =\dbE\Big[ \big\la \pa_{xx} L \eta, \eta\big\ra +\big\la  \tilde \dbE_{\cF^1_T}[\pa_{x\mu} L  \tilde \eta], \eta\big\ra+ 2\big\la \pa_{ax} L \eta+{1\over 2} \tilde\dbE_{\cF^1_T}[\pa_{a\mu} L \tilde \eta], \eta'\big\ra\big]+  \big\la \pa_{aa} L \eta', \eta'\big\ra\Big]\\
&=&\dbE\Big[ \big\la \pa_{xx} L \eta, \eta\big\ra +\big\la  \tilde \dbE_{\cF^1_T}[\pa_{x\mu} L  \tilde \eta], \eta\big\ra + \big| [\pa_{aa}L]^{1\over 2} \eta' + [\pa_{aa}L]^{-{1\over 2}} \big[\pa_{ax} L \eta+{1\over 2} \tilde\dbE_{\cF^1_T}[\pa_{a\mu} L \tilde \eta]\big]\big|^2\\
&&- \big|  [\pa_{aa}L]^{-{1\over 2}} \big[\pa_{ax} L \eta+{1\over 2} \tilde\dbE_{\cF^1_T}[\pa_{a\mu} L \tilde \eta]\big]\big|^2\Big]\\
&=&\dbE\Big[ \big\la \pa_{xx} L \eta, \eta\big\ra +\big\la  \tilde \dbE_{\cF^1_T}[\pa_{x\mu} L  \tilde \eta], \eta\big\ra - \big|  [\pa_{aa}L]^{-{1\over 2}} \big[\pa_{ax} L \eta+{1\over 2} \tilde\dbE_{\cF^1_T}[\pa_{a\mu} L \tilde \eta]\big]\big|^2\Big].
\eeaa
This is exactly \eqref{Lconvex}.
 \qed

\begin{rem} Observe that \reff{Lconvex3} expresses a certain convexity property of $L$. To illustrate this, consider the separable case,  where $L(x,\mu, a) := L_0(x,\mu) + L_1(x, a).$ Then  $\L=\L_0 + \L_1$, where 
\begin{align*}
&\L_0 := {d^2\over d\e d\d}\dbE\Big[ L_0 \big(\xi + \e\eta+\d\eta, \cL_{\xi+\e\eta}\big)\Big]\Big|_{(\e,\d)=(0,0)},\\
&\dis\L_1:=  {d^2\over d\e d\d}\dbE\Big[ L_1 \big(\xi + \e\eta+\d\eta,  \xi'+\e\eta'+\d\eta'\big)\Big]\Big|_{(\e,\d)=(0,0)},
\end{align*}
and so, $\L_0 \ge 0, \L_1\ge 0$ implies \reff{Lconvex}. Note that $\L_1\ge 0$ exactly means $L_1$ is convex in $(x,a)$. Moreover, consider the potential game case for $L_0$: $\pa_x L_0(x, \mu) = \pa_\mu \wh L_0(\mu, x)$ for some function $\wh L_0(\mu)$. Then 
\begin{align*}
\L_0&= {d\over d\e} \dbE\Big[ \big\la\pa_x L_0 \big(\xi + \e\eta, \cL_{\xi+\e\eta}\big), \eta\big\ra \Big]\Big|_{\e=0} = 
 {d\over d\e} \dbE\Big[ \big\la\pa_\mu \wh L_0 \big(\cL_{\xi+\e\eta}, \xi + \e\eta\big), \eta\big\ra \Big]\Big|_{\e=0}\\ 
 &={d^2\over d\e^2} \wh L_0(\cL_{\xi+\e\eta})\Big|_{\e=0}.
\end{align*}
Thus $\L_0\ge 0$ exactly means $\wh L_0$ is displacement convex, namely the mapping $\xi \mapsto \wh L_0(\cL_\xi)$ is convex. 
These are the same displacement convexity assumptions on the data for potential deterministic mean field master equations, imposed in \cite{GM}. In particular, \eqref{Lconvex3} is reminiscent to the joint convexity assumption on the Lagrangian (Assumption (H7)) in \cite{GM}.
\end{rem}

%
%
%
%

\section{The displacement monotonicity of $V$}
\label{sect-Vconvex}
\setcounter{equation}{0}

In this section we show that under our standing assumptions the displacement monotonicity condition is propagated along any classical solution of the master equation. More precisely, let $H$ and $G$ satisfy our standing assumptions of the previous section and in particular suppose that they are displacement monotone in the sense of Assumption \ref{assum-convex}. 

\begin{thm}\label{thm-convex}
Let Assumptions \ref{assum-regG}  and  \ref{assum-regH}-(i)(iv) and \ref{assum-convex} hold, and $V$ be a classical solution of the master equation \eqref{master}. Assume further that 
\small
$$
V(t,\cdot,\cdot), \pa_{x}V(t,\cdot,\cdot), \pa_{xx}V(t,\cdot,\cdot)\in \cC^2(\mathbb R^d\times\mathcal{P}_2), \quad \pa_\mu V(t,\cdot,\cdot,\cdot), \pa_{x\mu}V(t,\cdot,\cdot,\cdot)\in \cC^2(\mathbb R^d\times\mathcal{P}_2\times\mathbb R^d),$$ \normalsize
and all their derivatives in the state and probability measure variables are also continuous in the time variable and are uniformly bounded. Then $V(t,\cdot,\cdot)$ satisfies \reff{displacement2} for all $t\in[0,T]$. 
\end{thm}
\proof Without loss of generality, we shall prove the thesis of the theorem only for $t_0=0$, i.e. that $V(0,\cdot,\cdot)$ satisfies \reff{displacement2}. 

Fix $\xi, \eta\in \dbL^{2}(\cF_0)$. Let us consider the following decoupled McKean-Vlasov SDEs.
\small\begin{align}\label{XY}
\nonumber X_t &= \xi -\int_0^t \pa_pH(X_s, \mu_s, \pa_x V(s, X_s, \mu_s)) ds +  B_t+\beta B_t^0,\q \mu_t := \cL_{X_t|\mathcal{F}_t^{0}};\\
\delta X_t &= \eta -\int_0^t \Big[H_{px} (X_s) \delta  X_s + {1\over 2}\tilde \dbE_{\mathcal{F}_s}[ H_{p\mu}(X_s,\tilde X_s) \delta  \tilde X_s] + H_{pp}(X_s)N_s\Big]ds,\q\mbox{where}\\
\nonumber N_t &:=\tilde \dbE_{\mathcal{F}_t}[\pa_{x\mu}  V(X_t, \tilde X_t) \delta \tilde X_t] + \pa_{xx} V(X_t)  \delta  X_t + {1\over 2}H_{pp}(X_t)^{-1}\tilde \dbE_{\mathcal{F}_t}[ H_{p\mu}(X_t,\tilde X_t) \delta  \tilde X_t].
\end{align}\normalsize
Here and in the sequel,  for simplicity of notation, we omit the variables $(t, \mu_t)$, as well as the dependence on $\partial_x V$ and denote 
\begin{equation}\label{eq:dec24.2020}
\left.\ba{c}
H_p(X_t) := \pa_p H(X_t, \mu_t, \pa_x V(t, X_t, \mu_t)), \\ H_{p\mu}(X_t, \tilde X_t) := \pa_{p\mu} H(X_t, \mu_t, \tilde X_t, \pa_x V(t, X_t, \mu_t)),
\ea\right.
\end{equation}
and similarly for $H_{xp}, H_{pp}$,  $H_{x\mu}$, $\partial_{xx}V$, $\partial_{x\mu}V$. Since $V$ is assumed to be regular enough with $\pa_xV,\pa_{xx}V, \pa_{x\mu}V$ uniformly bounded and $H$ satisfies Assumption \ref{assum-regH}-(i), the driving vector field is globally Lipschitz continuous. Therefore, classical results imply the existence of unique solutions $X_t$ and $\delta X_t$.  We also observe that $\delta  X_t$  can be interpreted as $\dis\lim_{\e\to 0}{1\over \e}[X^{\xi+\e\eta}_t -X^\xi_t]$ (cf. \cite{BLPR}). 

Below, we shall use the notation, for $\th\in \dbR^d$ 
\bea
\label{notation}
\th^{\top}\pa_{xx\mu}V(x,\tilde x):=\sum_{i=1}^d \th_i\pa_{x_ix\mu}V(x,\tilde x), \q \tr(\pa_{\mu\mu})\pa_{x\mu}V:=\sum_{i=1}^d \partial_{\mu_i \mu_i} \partial_{x \mu} V,
\eea
and similarly for other higher order derivatives of $V$. Introduce: 
\[
I(t):=\tilde\dbE\Big[\langle \pa_{x\mu} V(t, X_t, \mu_t, \tilde X_t)\delta  \tilde X_t , \delta  X_t\rangle \Big], \quad 
\bar I(t):=\dbE\Big[\langle\pa_{xx} V(t, X_t, \mu_t) \delta  X_t,\delta  X_t\rangle\Big].
\]
We remark that, since  $(\tilde X_t, \d\tilde X_t)$ is a conditionally independent copy of $(X_t, \d X_t)$ and $\mu_t$ is $\cF^0_t$-measurable, for the notations in Section \ref{sect-space} we have, 
\bea
\label{conddisplacement}
I(t) + \bar I(t) =\dbE \Big[ (d_x d)_{X_t(\o^0,\cd)}V(t, \cd) (\d X_t(\o^0, \cd), \d X_t(\o^0,\cd))\Big].
\eea
 Our plan is to show that 
\bea
\label{dotI}
\dot I(t)+\dot {\bar I}(t) \leq 0.
\eea
Then, recalling $V(T,\cd) = G$ and applying Assumption \ref{assum-convex} (i), 
\begin{align*}
(d_x d)_{\xi} V(0,\cd)(\eta, \eta)& = I(0)+\bar I(0) \ge I(T)+\bar I(T)\\ 
& = \dbE\Big[ (d_x d)_{X_T(\o^0,\cd)}G(\d X_T(\o^0, \cd), \d X_T(\o^0,\cd))\Big] \ge 0.
\end{align*}
That is, $V(0,\cd)$ satisfies \reff{displacement2}.

To show \reff{dotI}, we apply  It\^o's formula \reff{Ito} to obtain 
\begin{equation}\label{eq:ito_double}
 \dot I(t)= I_1 + I_2 + I_3,
\end{equation}
where, introducing another conditionally independent copy $\hat X$ of $X$  and defining $\hat{\bar{\tilde \dbE}}$ in the manner of \reff{product3},  
\begin{align*}
I_1&:=\hat{\bar{\tilde \dbE}}\bigg[\bigg\langle \Big\{\pa_{tx\mu} V(X_t, \tilde X_t)  + {\widehat \beta^2\over 2} ((\tr\pa_{xx})\pa_{x\mu} V)(X_t,  \tilde X_t)  -H_p (X_t)^\top\pa_{xx\mu} V(X_t,  \tilde X_t)\\
&+\beta^2(\tr(\pa_{x\mu})\pa_{x\mu}V)(X_t,\bar X_t,\tilde X_t)+\beta^2(\tr(\pa_{\tilde x\mu})\pa_{x\mu}V)(X_t,\bar X_t,\tilde X_t)\\
&+\beta^2(\tr(\pa_{\tilde x x})\pa_{x\mu}V)(X_t,\tilde X_t)+\frac{\beta^2}{2}(\tr(\pa_{\mu\mu})\pa_{x\mu}V)(X_t,\hat X_t,\bar X_t,\tilde X_t) \\
&+  {\widehat\beta^2\over 2}(\tr(\pa_{\bar x\mu}) \pa_{ x\mu} V)(X_t, \bar X_t, \tilde X_t)  - H_p(\bar X_t)^\top\pa_{\mu x\mu} V(X_t, \bar X_t, \tilde X_t)\\
& + {\widehat\beta^2\over 2} (\tr(\pa_{\tilde x\tilde x})\pa_{x\mu} V)(X_t, \tilde X_t) - H_p(\tilde X_t)^\top\pa_{\tilde xx\mu} V(X_t, \tilde X_t)   \Big\}\delta  \tilde X_t, \delta  X_t\bigg\rangle \bigg],
\end{align*}
and,  rewriting $(\tilde X, \d\tilde X, \tilde \dbE)$ in the expression of $N$  as $(\bar X, \d\bar X, \bar \dbE)$ (which does not change the value of $N$),
\begin{align*}
I_2&:=-\bar{\tilde \dbE}\bigg[\Big\langle \pa_{\mu x} V(X_t, \tilde X_t) \Big\{ \big[H_{px} (X_t) + H_{pp} (X_t)\pa_{xx} V(X_t) \big]\delta  X_t +{\bf II_2} 
 \Big\},\delta \tilde X_t\Big\rangle \bigg],\\
{\bf II_2}&:=\Big[ H_{p\mu}(X_t, \bar X_t) + H_{pp}(X_t)\pa_{x\mu} V(X_t, \bar X_t) \Big]\delta \bar X_t,\\
I_3&:=- \bar{\tilde \dbE}\bigg[\Big\langle \pa_{x\mu} V(X_t,  \tilde X_t) \Big\{ \big[H_{px}(\tilde X_t)+ H_{pp}(\tilde X_t)\pa_{xx} V(\tilde X_t)\big] \delta \tilde X_t + {\bf III_3} \Big\}, \delta  X_t\Big\rangle\bigg],\\
 {\bf III_3}&:= \Big[H_{p\mu}(\tilde X_t, \bar X_t)+H_{pp}(\tilde X_t) \pa_{x\mu} V(\tilde X_t,\bar X_t) \Big]\delta \bar X_t.
\end{align*}
 We apply $-\partial_{x\mu}$ to \reff{master} and rewrite $(\tilde \xi, \bar \xi, \bar{\tilde\dbE})$ in \reff{eq:M-operator} as $(\bar\xi, \hat\xi, \hat{\bar \dbE})$ to obtain 
\begin{equation}\label{paxmucLV}
 0 = -(\pa_{x\mu} \sL V)(t, x, \mu, \tilde x) = J_1+ J_2 +J_3.
\end{equation} 
\noindent Here, we have set, recalling the notation in \reff{notation}, 
\begin{align*}
J_1:= &\pa_{tx\mu } V (x, \tilde x)+ {\widehat \beta^2\over 2}(\tr({\pa_{xx}})\pa_{x\mu} V)(x, \tilde x)  - H_{x\mu}(x,\tilde x) - \pa_{xx} V(x) H_{p\mu}(x,\tilde x)\\
- &\Big(H_{xp}(x) +\pa_{xx} V(x) H_{pp}(x) \Big) \pa_{x\mu} V(x, \tilde x) - H_p(x)^\top \pa_{xx\mu} V(x, \tilde x), \\
J_2:= &
{{\widehat\beta^2\over 2} (\pa_{x \tilde x } \tr(\pa_{\tilde x\mu})V)(x, \tilde x)}- H_p(\tilde x)^\top\pa_{\tilde x x \mu} V(x,\tilde x)\\ 
& - \pa_{x\mu} V(x,\tilde x) \Big(H_{px} (\tilde x) + H_{pp}(\tilde x) \pa_{xx} V(\tilde x)\Big) \nonumber\\
&+ \beta^2(\pa_{x\tilde x}\tr(\pa_{x\mu})V)(x,\tilde x)+\beta^2\bar\dbE\big[(\pa_{x\tilde x}\tr(\pa_{\mu\mu})V)(x,\bar\xi,\tilde x)\big]\nonumber\\
J_3:=& \hat {\bar\dbE}\bigg[{\widehat \beta^2\over 2}(\tr(\pa_{\bar x\mu}) \pa_{x \mu  } V)(x, \tilde x, \bar \xi) - H_p(\bar \xi)^\top \pa_{\mu x \mu} V(x,  \tilde x, \bar \xi)\\ 
&- \pa_{x\mu} V(x, \bar \xi) \Big[ H_{p\mu}(\bar\xi, \tilde x) + H_{pp} (\bar \xi) \pa_{x\mu} V(\bar \xi, \tilde x)\Big] \nonumber\\
&+\beta^2(\tr(\pa_{x\mu})\pa_{x\mu}V)(x,\tilde x,\bar \xi)+\frac{\beta^2}{2}(\tr(\pa_{\mu\mu})\pa_{x\mu}V)(x,\tilde x,\hat\xi,\bar\xi)\bigg].
\end{align*}


Note that we can switch the order of the differentiation in $\partial_{x}\partial_\mu,$ $\partial_x\partial_{\tilde x} $, etc. We emphasize that special care is needed when considering $\pa_\mu \pa_{\tilde x}$ (since we cannot change here the order of differentiation). For such terms, we use their symmetric properties given in Lemma \ref{lem-order}. 
By evaluating \eqref{paxmucLV} along $(X_t,\mu_t,\tilde X_t)$ and plugging into  \reff{eq:ito_double}, one can cancel many terms and simplify the previous derivation as
\begin{align*}
\dot I(t)&=\bar{\tilde \dbE}\bigg[ - \Big\langle  \pa_{\mu x} V(X_t, \tilde X_t) \big[ H_{p\mu}(X_t, \bar X_t) + H_{pp}(X_t)\pa_{x\mu} V(X_t, \bar X_t)\big]\delta \bar X_t, \delta \tilde X_t\Big\rangle \nonumber \\
&+ \Big\langle \big[H_{x\mu}(X_t,\tilde X_t) + \pa_{xx} V(X_t)H_{p\mu}(X_t,\tilde X_t) \big]\delta \tilde X_t,\delta  X_t\Big\rangle \bigg].
\end{align*}
Thus, by using the tower property of conditional expectations and the conditional i.i.d. property of $(X, \d X), (\tilde X, \d \tilde X), (\bar X, \d \bar X)$, we have
\begin{align}
\label{ExmuV}
\dot I(t) &=\dbE\bigg[-\Big\langle H_{pp}(X_t) \tilde \dbE_{\mathcal{F}_t}\big[\pa_{x\mu}  V(X_t, \tilde X_t) \delta \tilde X_t\big],~\tilde \dbE_{\mathcal{F}_t}\big[\pa_{x\mu}  V(X_t, \tilde X_t) \delta \tilde X_t\big]\Big\rangle \\
&~ -\Big\langle \tilde \dbE_{\mathcal{F}_t}\big[ H_{p\mu}(X_t,\tilde X_t)\delta  \tilde X_t\big] , \tilde \dbE_{\mathcal{F}_t}\big[ \pa_{x\mu} V(X_t, \tilde X_t)  \delta  \tilde X_t\big] - \pa_{xx} V(X_t) \delta  X_t\Big\rangle  \nonumber\\
&+ \Big\langle \tilde{\mathbb E}_{\mathcal{F}_t}\big[H_{x\mu}(X_t,\tilde X_t) \delta  \tilde X_t\big], \delta  X_t\Big\rangle\bigg].\nonumber
\end{align}

Similarly as above, we apply It\^o formula {\reff{Ito}} to $\bar I(t)$ to obtain 
\[
\dot {\bar I}(t)= \bar I_1 + \overline{I}_2 + \overline{I}_3,
\] 
\no where, 
\beaa
\bar I_1&:=& \tilde\dbE\bigg[\Big\langle\Big\{\pa_{txx} V(X_t)  + {\widehat\beta^2\over 2} (\tr(\pa_{xx})\pa_{xx} V)(X_t)  -H_p(X_t)^\top \pa_{xxx} V(X_t)\\
&& +\beta^2(\tr(\pa_{x\mu})\pa_{xx}V)(X_t,\tilde X_t) \Big\} \delta X_t, \delta X_t \Big\rangle \bigg],\\
%
\overline{I}_2 &:=&\bar{\tilde\dbE}\bigg[\Big\langle\Big\{\frac{\beta^2}{2}(\tr(\pa_{\mu\mu})\pa_{xx}V)(X_t,\tilde X_t,\bar X_t)\\
&&+  {\widehat\beta^2\over 2} (\tr(\pa_{\tilde x\mu})\pa_{ xx} V)(X_t,\tilde X_t) -  H_p(\tilde X_t)^\top\pa_{\mu xx} V(X_t,  \tilde X_t)\Big\} \delta  X_t, \delta  X_t\Big\rangle\bigg],\\
\overline{I}_3&:=& \tilde\dbE\bigg[- 2\Big\langle\pa_{xx} V(X_t)  \Big\{ \big[H_{px} (X_t) + H_{pp} (X_t)\pa_{xx} V(X_t) \big]\delta  X_t\\
&&+  \big[ H_{p\mu}(X_t, \tilde X_t) + H_{pp}(X_t)\pa_{x\mu} V(X_t, \tilde X_t) \big] \delta \tilde X_t\Big\},\delta  X_t \Big\rangle \bigg].
\eeaa
On the other hand, applying $-\partial_{xx}$ to \reff{master} we obtain  
\begin{equation}\label{paxxVmaster}
0=  -(\pa_{xx} \sL V)(t, x, \mu)= \bar J_1 + \bar{J_2},
\end{equation}
\no where 
\begin{align*}
 {\bar J_1}&:= \pa_{txx} V + {\widehat\beta^2\over 2}(\tr(\pa_{xx})\pa_{xx} V) - H_{xx}(x) - 2H_{xp}(x)\pa_{xx} V(x)\\
&- \pa_{xx} V(x)H_{pp}(x) \pa_{xx} V(x)-H_p(x)^{\top} \pa_{xxx} V(x),\\
 {\bar J_2}&:=\bar{\tilde \dbE}\Big[{\widehat\beta^2\over 2} (\tr(\pa_{\tilde x\mu})\pa_{xx} V)(x,  \tilde  \xi) - H_p( \tilde \xi)^\top\pa_{\mu xx} V(x,  \tilde \xi)\\ 
&+\beta^2(\tr(\pa_{x\mu})\pa_{xx}V)(x,\tilde\xi)+\frac{\beta^2}{2}(\tr(\pa_{\mu\mu})\pa_{xx}V)(x,\bar\xi,\tilde\xi)\Big].
\end{align*}
\normalsize
We evaluate the previous expression along $(X_t,\mu_t)$ to obtain after a simplification 
\small
\begin{align}\label{ExxV}
\dot {\bar I}(t) & =\dbE\bigg[ - \Big\langle H_{pp}(X_t)\pa_{xx} V(X_t)  \delta  X_t, \pa_{xx} V(X_t)\delta  X_t\Big\rangle\\ 
&\nonumber- 2 \Big\langle H_{pp}(X_t) \pa_{xx} V(X_t)  \delta  X_t,  \tilde\dbE_{\mathcal{F}_t} \Big[\pa_{x\mu} V(X_t,\tilde X_t)\delta   \tilde X_t\Big]\Big\rangle\\
& - 2\Big\langle \pa_{xx} V(X_t)  \delta  X_t,\tilde \dbE_{\mathcal{F}_t}\Big[H_{p\mu}(X_t, \tilde X_t)\delta  \tilde X_t\Big]\Big\rangle +  \Big\langle H_{xx}(X_t) \delta  X_t, \delta  X_t\Big\rangle\bigg].\nonumber
\end{align}
\normalsize

We combine \reff{ExmuV} and \reff{ExxV} to deduce that,
\small
\begin{align}\label{Vconvex0}
\dot I(t) +\dot {\bar I}(t)&= \dbE\bigg[-\Big|H^{\frac{1}{2}}_{pp}(X_t) \Big\{\tilde \dbE_{\mathcal{F}_t}\big[\pa_{x\mu}  V(X_t, \tilde X_t) \delta  \tilde X_t\big] + \pa_{xx} V(X_t)  \delta  X_t\Big\}\Big|^2\\
& -\Big\langle \tilde \dbE_{\mathcal{F}_t}\big[ H_{p\mu}(X_t,\tilde X_t) \delta  \tilde X_t\big],  \tilde \dbE_{\mathcal{F}_t}\big[ \pa_{x\mu} V(X_t, \tilde X_t)  \delta  \tilde X_t\big] + \pa_{xx} V(X_t) \delta  X_t\Big\rangle \nonumber \\
& +\Big\langle \tilde \dbE_{\mathcal{F}_t}\big[H_{x\mu}(X_t,\tilde X_t) \delta  \tilde X_t\big] + H_{xx}(X_t) \delta  X_t, \delta  X_t\Big\rangle\bigg]\nonumber\\
&= \dbE\bigg[-\Big|H^{\frac{1}{2}}_{pp}(X_t) \Big\{\tilde \dbE_{\mathcal{F}_t}\big[\pa_{x\mu}  V(X_t, \tilde X_t)\delta  \tilde X_t\big] + \pa_{xx} V(X_t) \delta  X_t\Big\}\nonumber\\ 
&+\frac{1}{2}H^{-\frac{1}{2}}_{pp}(X_t) \tilde \dbE_{\mathcal{F}_t}\big[ H_{p\mu}(X_t,\tilde X_t)\delta  \tilde X_t\big] \Big|^2 \nonumber+ \Big\langle \tilde{\mathbb E}_{\mathcal{F}_t}\big[H_{x\mu}(X_t,\tilde X_t) \delta   \tilde X_t \big],\delta  X_t\Big\rangle\\ 
&+\Big\langle H_{xx}(X_t) \delta  X_t,\delta  X_t\Big\rangle +\frac{1}{4}\Big|H_{pp}^{-\frac{1}{2}}(X_t)\tilde \dbE_{\mathcal{F}_t}\big[ H_{p\mu}(X_t,\tilde X_t) \delta  \tilde X_t\big]\Big|^2 \bigg]\nonumber\\
&=- \dbE\Big[\Big|H^{\frac{1}{2}}_{pp}(X_t) N_t \Big|^2\Big] + \dbE^{\dbP_0}\Big[({ displ}^{\f}_{X_t(\o^0,\cd)} H)\big(\d\tilde X_t(\o^0,\cd), \d\tilde X_t(\o^0,\cd)\big)\Big],\nonumber
\end{align}
\normalsize
where we have set $~\f(x):= \pa_x V(t,x, \mu_t(\o^0))$ and the last line is in the spirit of \reff{conddisplacement}. 
Applying \reff{Hconvex} 
we obtain \reff{dotI} immediately.
%
\qed

\begin{rem}
\label{rem-Vconvex}
 (i) The main trick here is that we may complete the square in \reff{Vconvex0} for the terms involving $\pa_{xx} V$ and more importantly $\pa_{x\mu} V$, which is hard to estimate a priori. Since the identity is exact, \reff{Hconvex} seems essential for not loosing displacement monotonicity. 
 
Moreover, recalling \reff{conddisplacement} we see that
\bea\label{rate}
\left.\ba{c}
\dis{d\over dt} \dbE^{\dbP_0}\Big[ (d_x d)_{X_t(\o^0,\cd)}V(t, \cd, \cd) (\d X_t(\o^0, \cd), \d X_t(\o^0,\cd))\Big]\ss\\ 
\dis\le  \dbE^{\dbP_0}\Big[({\rm{displ}}^{\pa_xV(t,\cd, \mu_t(\o^0))}_{X_t(\o^0,\cd)} H)\big(\d\tilde X_t(\o^0,\cd), \d\tilde X_t(\o^0,\cd)\big)\Big].
\ea\right.
\eea
\no So, roughly speaking, ${\rm{displ}}\,H$ measures the rate of dissipation of the displacement monotonicity of $V$, through the bilinear form $(d_x d)V(t, \cdot, \cdot)$. 

(ii) In the separable case, i.e. $H(x,\mu, p) = H_0(x,p)-F(x,\mu) $ for some $H_0$ and $F$, \reff{ExmuV} becomes 
\begin{align*}
\dot I(t)&=-\dbE\Big[\big|[(H_0)_{pp}(X_t)]^{1\over 2} \tilde \dbE_{\mathcal{F}_t}[\pa_{x\mu}  V(X_t, \tilde X_t) \delta\tilde X_t]\big|^2+\tilde{\mathbb E}_{\mathcal{F}_t}\langle[\partial_{x\mu}F(X_t,\tilde X_t)\delta \tilde X_t],\delta X_t\rangle\Big].
\end{align*}
The term involving $\pa_{x\mu} V$ is again in a complete square, and it is no surprise that $V(t,\cdot,\cdot)$ would satisfy Lasry-Lions monotonicity condition \reff{mon2} provided that the data $G$ and $F$ also satisfy 
\reff{mon2}. So, our arguments provide an alternative  proof for the propagation of the  Lasry-Lions monotonicity along $V(t,\cdot,\cdot)$ (and for the global well-posedness of the master equation, as a consequence of it, just as in the rest of the paper) in the case of separable Hamiltonians and Lasry-Lions monotone data.

(iii) When $H$ is non-separable, however, it remains a challenge to find  sufficient conditions on $H$ that could ensure the right hand side of \reff{ExmuV}  being negative (for arbitrary times). This makes the propagation of the Lasry-Lions monotonicity condition along $V(t,\cdot,\cdot)$, hard to envision.  In  \cite{Lions}, a notion of monotonicity condition for non-separable Hamiltonians that depend \emph{locally} on the measure variable was proposed  (see also \cite{AP}).  This condition allows to obtain uniqueness of solutions for the corresponding MFG system.
\end{rem}

\section{The uniform Lipschitz continuity of $V$ under $W_2$}
\label{sect-Lipschitz}
\setcounter{equation}{0}
 The main result in this section is that, the displacement semi-monotone solutions to the master equation \eqref{master} are always uniformly $W_2$--Lipschitz continuous. We note that this observation that the Lipschitz continuity of $V(t,\cdot,\cdot)$ in $\mu$ (under $W_2$) is the consequence of the displacement semi-monotonicity of $V(t,\cdot,\cdot)$ only (other than the technical conditions), seems to be a new even in the separable case. We remark again that  the displacement semi-monotonicity is weaker than both the displacement monotonicity and the Lasry-Lions monotonicity (if $\partial_{xx}V$ is uniformly bounded), see Remark \ref{rem-semi}.
\begin{thm}
\label{thm-Lipschitz}
Let all the conditions in Theorem \ref{thm-convex} hold, except that we do not require Assumption \ref{assum-convex} (ii). Assume further that $V(t,\cdot,\cdot)$ satisfies the displacement semi-monotonicity \eqref{semidisplacement1} for each $t\in[0,T]$. Then $V$ and $\pa_x V$ are uniformly Lipschitz continuous in $\mu$ under $W_2$ with  Lipschitz constant $C^\mu_2$, where $C^\mu_2>0$ depends only on $d, T, \|\pa_{x}V\|_{L^\infty},\|\pa_{xx}V\|_{L^\infty}$, the  $L^G_2$ in Remark \ref{rem:W2}-(ii),  the $L^H(\|\pa_xV\|_{L^\infty})$ in Assumption  \ref{assum-regH}-(i), and the $\lambda$ in \eqref{semidisplacement1}.
%
\end{thm}
\proof  In this proof, $C>0$ denotes a generic constant depending only on quantities mentioned in the statement of the theorem. Without loss of generality, we show the thesis of the theorem only for $t_0=0$. We fix $\xi, \eta\in \dbL^{2}(\cF_0)$ and continue to use the notation as in the proof of Theorem \ref{thm-convex}. In particular, $\delta X$ is defined by  \reff{XY}. First we emphasize that the equality \reff{Vconvex0} does not rely on \reff{Hconvex}. Then,  integrating \reff{Vconvex0} over $[0, t]$ we obtain:
\begin{align*}
&\dis \int_0^t \dbE\big[\big|H_{pp}(X_s)^{\frac{1}{2}} N_s \big|^2\big]ds = [I(0) + \bar I(0)] -  [I(t) +\bar I(t)]\\ 
&+  \int_0^t \dbE\bigg[ \Big\langle \tilde{\mathbb E}_{\mathcal{F}_s}\big[H_{x\mu}(X_s,\tilde X_s) \delta   \tilde X_s \big],\delta  X_s\Big\rangle 
\dis +\Big\langle H_{xx}(X_s) \delta  X_s,\delta  X_s\Big\rangle\\ 
&+\frac{1}{4}\Big|H_{pp}^{-\frac{1}{2}}(X_s)\tilde \dbE_{\mathcal{F}_s}\big[ H_{p\mu}(X_s,\tilde X_s) \delta  \tilde X_s\big]\Big|^2 \bigg]ds\\
&\dis\le I(0)  -  [I(t) +\bar I(t)] + C\dbE[|\eta|^2]+C\int_0^t \dbE[|\d X_s|^2]ds,
\end{align*}
where we used the bound of  $\pa_{xx} V,  H_{x\mu}, H_{xx}, H_{p\mu}$. Since $V(t,\cd,\cd)$  satisfies \eqref{semidisplacement1}, by \reff{conddisplacement} we have $I(t) +\bar I(t)\ge -\l\mathbb E[|\delta X_t|^2]$. Then
\begin{equation}\label{Nest}
  \int_0^t \dbE\big[\big|H_{pp}(X_s)^{\frac{1}{2}} N_s \big|^2\big]ds \le I(0)+C\mathbb E[|\delta X_t|^2]  + C\dbE[|\eta|^2]+ C\int_0^t \dbE[|\d X_s|^2]ds.
\end{equation}

Next, using \reff{XY} and Young's inequality, we have for any $\epsilon>0$
\[
|\delta X_t|^2 \leq |\eta|^2 + C_\epsilon\int_0^t |\delta X_s|^2ds+\epsilon\int_0^t |H_{pp}(X_s)^{\frac{1}{2}}N_s|^2ds.
\]
Taking expectations on both sides and choosing $\epsilon>0$ small enough, by \reff{Nest} we obtain 
\[
\mathbb E\big[|\delta X_t|^2\big] \leq   C \int_0^t \dbE\big[|\delta X_s|^2\big]ds+C \mathbb E\big[|\eta|^2\big]+ C |I(0)|.
\]
Then it follows from Gr\"onwall's inequality  that 
\bea\label{EY1} 
\left.\ba{c}
\dis \sup_{t \in [0,T]}\mathbb E\big[|\delta X_t|^2\big] \leq C \mathbb E\big[|\eta|^2\big]+C |I(0)|\le C \mathbb E\Big[|\eta|^2+ |\eta| |\Upsilon_0|\Big],\\
\dis\mbox{where}\q \Upsilon_t := \tilde \dbE_{\cF_t}\Big[\pa_{x\mu} V(t, X_t, \mu_t, \tilde X_t) \d\tilde X_t\Big].
\ea\right.
\eea

We shall follow the arguments in Theorem \ref{thm-convex} to estimate $\Upsilon$. We first observe that,
\bea
\label{EtT}
\Upsilon_t = \tilde \dbE_{\cF_T}\Big[\pa_{x\mu} V(t, X_t, \mu_t, \tilde X_t) \d\tilde X_t\Big].
\eea
So, by applying It\^o formula \reff{Ito} on $\pa_{x\mu} V(t, X_t, \mu_t, \tilde X_t) \d\tilde X_t$,  taking conditional expectation $\tilde \dbE_{\cF_T}$, and then changing back to $\tilde \dbE_{\cF_t}$ as in \reff{EtT}, we obtain 
\begin{equation}\label{dUpsilon}
d\Upsilon_t =  (dB_t)^\top K_1(t) + \b  (dB^0_t)^\top K_2(t)  + [K_3(t) - K_4(t)]dt,
\end{equation}
where, recalling the notation in \eqref{notation} (in particular the stochastic integral terms above are column vectors),
\beaa
&&\dis K_1(t) := \tilde{\mathbb E}_{\mathcal{F}_t}\Big[\pa_{xx\mu}V(X_t,\tilde X_t)\delta \tilde X_t\Big],\\
&&\dis K_2(t):= K_1(t) + \bar{\tilde \dbE}_{\mathcal{F}_t}\Big[\Big\{(\pa_{\mu x\mu}V)(X_t,\bar X_t,\tilde X_t)  +\pa_{\tilde xx\mu}V(X_t,\tilde X_t)\Big\}\delta\tilde X_t\Big],\\
&&\dis K_3(t) := \hat{\bar{\tilde \dbE}}_{\cF_t} \bigg[ \bigg\{\pa_{tx\mu}V(X_t,\tilde X_t) -H_p( X_t)^\top\pa_{xx\mu}V(X_t,\tilde X_t) \\
&&\q -H_p(\tilde X_t)^{\top}\pa_{\tilde xx\mu}V(X_t,\tilde X_t) -H_p(\bar X_t)^\top\pa_{\mu x\mu}V(X_t,\bar X_t,\tilde X_t)\\
&&\q  + {\widehat\beta^2\over 2}\Big[ (\tr(\pa_{xx})\pa_{x\mu}V)(X_t, \tilde  X_t)+   (\tr(\pa_{\tilde x\tilde x})\pa_{x\mu}V)(X_t, \tilde  X_t)+(\tr(\pa_{\bar x\mu})\pa_{x\mu}V)(X_t,\bar X_t,\tilde X_t)\Big] \\
&&\q + \beta^2 \Big[(\tr(\pa_{x\mu})\pa_{x\mu}V)(X_t,\bar X_t,\tilde X_t)+ (\tr(\pa_{\tilde x x})\pa_{x\mu}V)(X_t,\tilde X_t)\\
 &&\q + (\tr(\pa_{\tilde x\mu})\pa_{x\mu}V)(X_t,\bar X_t,\tilde X_t) +\frac{1}{2} (\tr(\pa_{\mu\mu})\pa_{x\mu}V)(X_t,\hat X_t,\bar X_t,\tilde X_t)\Big]\bigg\}\delta\tilde X_t\bigg],\\
 &&\dis K_4(t) :=\bar{\tilde \dbE}_{\cF_t}\bigg[  \pa_{x\mu} V(X_t, \tilde X_t) \Big\{\big[H_{px}(\tilde X_t)+ H_{pp}(\tilde X_t)\pa_{xx} V(\tilde X_t)\big]   \delta\tilde X_t \\
 &&\q + \big[H_{p\mu}(\tilde X_t, \bar X_t)+H_{pp}(\tilde X_t) \pa_{x\mu} V(\tilde X_t,\bar X_t) \big]\delta\bar X_t\Big\}\bigg].
\eeaa\normalsize
In light of \reff{paxmucLV},  by straightforward calculation and simplification and setting 
\begin{align*}
&K_5(t) := H_{xp}(X_t) + \pa_{xx} V(X_t) H_{pp}(X_t),\\
&K_6(t) := \tilde\dbE_{\mathcal{F}_t}\Big[  \big[H_{x\mu}(X_t,\tilde X_t) +  \pa_{xx} V(X_t)H_{p\mu}(X_t,\tilde X_t) \big]\delta\tilde X_t\Big], 
\end{align*}
we derive  that 
\begin{equation}\label{dUpsilon2}
d\Upsilon_t = (dB_t)^\top K_1(t) + \b  (dB^0_t)^\top K_2(t) + \Big[ K_5(t) \Upsilon_t + K_6(t)\Big] dt.
\end{equation}
We have 
\beaa
\Upsilon_t = \Upsilon_T - \int_t^T (dB_s)^\top K_1(s) - \int_t^T \b  (dB^0_s)^\top K_2(s) - \int_t^T \Big[ K_5(s) \Upsilon_s + K_6(s)\Big] ds.
\eeaa
Take conditional expectation $\tilde \dbE_{\cF_t}$ and recall \reff{EY1}, we have
\bea
\label{Upsilon1}
\Upsilon_t = \tilde \dbE_{\cF_t}\Big[\pa_{x\mu} G(X_T, \mu_T, \tilde X_T) \d\tilde X_T\Big] - \int_t^T \tilde \dbE_{\cF_t}\Big[ K_5(s) \Upsilon_s + K_6(s)\Big] ds.
\eea
Then by \reff{dUpsilon2} and the required regularity of $G, H$ and $V$, in particular  \reff{L2G}, we have
\beaa
|\Upsilon_t|^2 \le C \tilde \dbE_{\cF_t}\big[|\d\tilde X_T|^2\big] + C\int_t^T \tilde \dbE_{\cF_t}\big[ |\Upsilon_s|^2 + |\d \tilde X_s|^2\big] ds.
\eeaa
Now take conditional expectation $\tilde \dbE_{\cF_0}$, we get
\beaa
\tilde \dbE_{\cF_0}\big[|\Upsilon_t|^2\big] \le C \tilde \dbE_{\cF_0}\big[|\d\tilde X_T|^2\big] + C\int_t^T \tilde \dbE_{\cF_0}\big[ |\Upsilon_s|^2 + |\d \tilde X_s|^2\big] ds.
\eeaa
Thus, by the Gr\"onwall inequality we have
\bea
\label{Upsilon0}
|\Upsilon_0|^2= \tilde \dbE_{\cF_0}\big[|\Upsilon_0|^2\big] \le C \tilde \dbE_{\cF_0}\big[|\d\tilde X_T|^2\big] + C\int_0^T \tilde \dbE_{\cF_0}\big[|\d \tilde X_s|^2\big] ds.
\eea
Plug this into \reff{EY1}, for any $\e>0$ we have
\beaa
\sup_{t \in [0,T]}\mathbb E\big[|\delta X_t|^2\big] \leq C_\e \mathbb E\big[|\eta|^2\big]+\e \dbE\big[ |\Upsilon_0|^2\big] \le C_\e \mathbb E\big[|\eta|^2\big]+C\e  \sup_{t \in [0,T]}\mathbb E\big[|\delta X_t|^2\big]. 
\eeaa
Set $\e = {1\over 2C}$ at above, we have
\bea\label{EY2} 
 \sup_{t \in [0,T]}\mathbb E\big[|\delta X_t|^2\big] \le C \dbE\big[|\eta|^2\big].
 \eea
Note that, recalling the setting in Section \ref{sect-space}, $\d \tilde X_t$ is measurable with respect to $\cF^0_t \vee \tilde \cF^1_t$, which is independent of $\cF_0$ under $\tilde \dbP$. Then the conditional expectation in the right side of \reff{Upsilon0} is actually an expectation. Plug \reff{EY2} into \reff{Upsilon0}, we have
 \bea
 \label{EUpsilon0}
\Big| \tilde \dbE_{\cF_0}\Big[\pa_{x\mu} V(0, \xi, \mu, \tilde \xi) \tilde \eta\Big] \Big|^2 = |\Upsilon_0|^2  \le  C\dbE\big[|\eta|^2\big].
\eea
This implies
\beaa
\Big|\tilde\dbE\big[\pa_{x\mu} V(0, x,  \mu, \tilde \xi)  \tilde \eta \big]\Big|\le  C(\mathbb E|\eta|^2)^{\frac{1}{2}},\q\mu-\mbox{a.e.} ~x.
\eeaa
Since $\pa_{x\mu} V$ is continuous, we have
\beaa
\Big|\dbE\big[\pa_{x\mu} V(0, x,  \mu, \xi)  \eta \Big]\Big|\le  C(\mathbb E|\eta|^2)^{\frac{1}{2}},\q\mbox{for all}~ x, \mu, \xi, \eta.
\eeaa
In particular, this implies that there exists a constant $C^\mu_2>0$ such that
\beaa
\Big|\pa_xV(0, x, \cL_{\xi+\eta}) - \pa_x V(0, x, \cL_\xi)\Big| = \Big|\int_0^1 \dbE\big[\pa_{x\mu} V(0, x, \cL_{\xi+\th \eta}, \xi + \th \eta) \eta\big] d\th\Big|\le C^\mu_2(\mathbb E|\eta|^2)^{\frac{1}{2}}.
\eeaa
Now, taking random variables $\xi,\eta$ such that $W^2_2(\cL_{\xi+\eta},\cL_{\xi})=\mathbb E|\eta|^2$, the above inequality
exactly means that $\pa_x V(0,x,\cdot)$ is uniformly Lipschitz continuous in $\mu$ under $W_2$ with uniform Lipschitz constant $C^\mu_2$.

Finally, denote
\beaa
\bar \Upsilon_t := \tilde\dbE_{\mathcal{F}_t}\Big[\pa_{\mu} V(t, X_t, \mu_t, \tilde X_t)  \delta\tilde X_t \Big]. 
\eeaa
Following similar arguments as in \reff{dUpsilon} we have
\bea
\label{dbarUpsilon}
d\bar \Upsilon_t =  (dB_t)^\top \bar K_1(t) + \b  (dB^0_t)^\top\bar  K_2(t)  + [\bar K_3(t) - \bar K_4(t)]dt,
\eea
where, 
\small
\beaa
\bar K_1(t) &:=& \tilde{\mathbb E}_{\mathcal{F}_t}\Big[\pa_{x\mu}V(X_t,\tilde X_t)\delta \tilde X_t\Big],\\
\bar K_2(t)&:=& K_1(t) + \bar{\tilde \dbE}_{\mathcal{F}_t}\Big[\Big\{(\pa_{\mu \mu}V)(X_t,\bar X_t,\tilde X_t)  +\pa_{\tilde x\mu}V(X_t,\tilde X_t)\Big\}\delta\tilde X_t\Big];\\
\bar K_3(t) &:=& \hat{\bar{\tilde \dbE}}_{\cF_t} \bigg[ \bigg\{\pa_{t\mu}V(X_t,\tilde X_t) -H_p( X_t)^\top\pa_{x\mu}V(X_t,\tilde X_t) \\
&&\q -H_p(\tilde X_t)^{\top}\pa_{\tilde x\mu}V(X_t,\tilde X_t) -H_p(\bar X_t)^\top\pa_{\mu \mu}V(X_t,\bar X_t,\tilde X_t)\\
&&\q  + {\widehat\beta^2\over 2}\Big[ (\tr(\pa_{xx})\pa_{\mu}V)(X_t, \tilde  X_t)+   (\tr(\pa_{\tilde x\tilde x})\pa_{\mu}V)(X_t, \tilde  X_t)+(\tr(\pa_{\bar x\mu})\pa_{\mu}V)(X_t,\bar X_t,\tilde X_t)\Big] \\
&&\q + \beta^2 \Big[(\tr(\pa_{x\mu})\pa_{\mu}V)(X_t,\bar X_t,\tilde X_t)+ (\tr(\pa_{\tilde x x})\pa_{\mu}V)(X_t,\tilde X_t)\\
 &&\qq + (\tr(\pa_{\tilde x\mu})\pa_{\mu}V)(X_t,\bar X_t,\tilde X_t) +\frac{1}{2} (\tr(\pa_{\mu\mu})\pa_{\mu}V)(X_t,\hat X_t,\bar X_t,\tilde X_t)\Big]\bigg\}\delta\tilde X_t\bigg],\\
\bar  K_4(t) &:=&\bar{\tilde \dbE}_{\cF_t}\bigg[  \pa_{\mu} V(X_t, \tilde X_t) \Big\{\big[H_{px}(\tilde X_t)+ H_{pp}(\tilde X_t)\pa_{xx} V(\tilde X_t)\big]   \delta\tilde X_t \\
 &&\qq + \big[H_{p\mu}(\tilde X_t, \bar X_t)+H_{pp}(\tilde X_t) \pa_{x\mu} V(\tilde X_t,\bar X_t) \big]\delta\bar X_t\Big\}\bigg].
\eeaa\normalsize
On the other hand, by  taking $-\partial_\mu$ of \reff{master} and omitting the variables $(t,\mu)$ we have,
\small
\begin{align*}
0 &= -\pa_{\mu} (\sL V)(t, x, \mu, \tilde x)\nonumber\\
& = \pa_{t\mu } V (x, \tilde x)+ {\widehat \beta^2\over 2}(\tr(\pa_{xx})\pa_{\mu} V)(x, \tilde x)  - H_\mu(x, \tilde x)- H_p(x)^\top \pa_{x\mu} V(x, \tilde x) \nonumber\\
&+{\widehat\beta^2\over 2} (\pa_{\tilde x}\tr(\pa_{ \tilde x \mu}) V)(x, \tilde x)+\beta^2(\pa_{\tilde x}\tr(\pa_{x\mu})V)(x,\tilde x)+\beta^2\bar{\mathbb E}\big[(\pa_{\tilde x}\tr(\pa_{\mu\mu})V)(x,\bar \xi,\tilde x)\big]\nonumber\\
& -H_p(\tilde x)^{\top}\pa_{\tilde x  \mu} V(x,\tilde x) - \pa_{\mu} V(x,\tilde x) \left(H_{px} (\tilde x) + H_{pp}(\tilde x) \pa_{xx} V(\tilde x)\right)\nonumber\\
& +  \hat{\bar \dbE}\bigg[{\widehat\beta^2\over 2}( \tr(\pa_{\bar x\mu})\pa_{\mu } V)(x, \tilde x, \bar\xi) +\beta^2(\tr(\pa_{x\mu})\pa_\mu V)(x,\tilde x,\bar\xi)+\frac{\beta^2}{2}(\tr(\pa_{\mu\mu})\pa_\mu V)(x,\tilde x,\hat\xi,\bar \xi)\\
& -H_p(\bar \xi)^\top \pa_{\mu  \mu} V(x,  \tilde x, \bar \xi) \nonumber
  - \pa_{\mu} V(x, \bar\xi) \left( H_{p\mu}(\bar\xi, \tilde x) + H_{pp} (\bar \xi) \pa_{x\mu} V(\bar \xi, \tilde x)\right) \bigg].\nonumber
\end{align*}
\normalsize 
Plug this into \reff{dbarUpsilon}, we have
\beaa
d\bar \Upsilon_t =  (dB_t)^\top \bar K_1(t) + \b  (dB^0_t)^\top\bar  K_2(t)  + \tilde \dbE_{\mathcal{F}_t}\Big[ H_{\mu} (X_t, \tilde X_t)\delta \tilde X_t\Big]dt.
\eeaa  
Then\small
\beaa
\bar \Upsilon_0  =  \tilde\dbE_{\cF_0} \Big[ \bar \Upsilon_T - \int_0^T  H_{\mu} (X_t, \tilde X_t)\delta \tilde X_tdt\Big]= \tilde\dbE_{\cF_0} \Big[ \pa_{\mu} G(X_T,  \tilde X_T) \delta \tilde X_T - \int_0^T  H_{\mu} (X_t, \tilde X_t)\delta \tilde X_tdt\Big],
\eeaa\normalsize
and thus, by \reff{L2G} again,
\beaa
\Big|\tilde\dbE_{\cF_0}\big[\pa_{\mu} V(0, \xi, \mu, \tilde \xi)  \tilde \eta \big]\Big|^2 = \big|\bar \Upsilon_0\big|^2  \le  C\tilde\dbE_{\cF_0} \Big[ |\delta \tilde X_T|^2 + \int_0^T  |\delta \tilde X_t|^2dt\Big].
\eeaa
Now by \reff{EY2}, follow the arguments for \reff{EUpsilon0} and the analysis afterwards, we see that $V(0,x,\cdot)$ is also uniformly Lipschitz continuous in $\mu$ under $W_2$ with uniform Lipschitz constant $C^\mu_2$.
\qed

%


We note that the a priori $W_2$-Lipschitz continuity of $V$ in $\mu$ is not sufficient for the global well-posedness of the master equation. We shall prove in the next section that together with other properties, it actually implies the uniform $W_1$-Lipschitz continuity. 
The following proposition, which can be viewed as  an analogue of Theorem \ref{thm-Lipschitz} for the version of Lasry-Lions monotonicity, obtains the $W_1$-Lipschtiz estimate directly. Since the proof should be standard for experts and is very similar to that of Theorem \ref{thm-Lipschitz}, we only sketch it and   focus on the main differences.

\begin{prop}
\label{prop-Lipschitz}
Let all the conditions in Theorem \ref{thm-convex} hold, except that we do not require Assumption \ref{assum-convex}. 
Assume further that there exists $\l>0$ such that, for any $\xi,\eta\in \mathbb L^2(\mathcal{F}_T^1)$,
\begin{equation}\label{semimon1}
\mathbb E\Big[\langle\pa_{x\mu}V(t,\xi,\mathcal{L}_{\xi},\tilde\xi)\tilde\eta,\eta\rangle\Big]\ge -\l \big(\mathbb E[|\eta|]\big)^2,\quad\text{for each $t\in[0,T]$},
\end{equation}
where $(\tilde\xi,\tilde\eta)$ is an independent copy of $(\xi,\eta)$.
Then $V$ and $\pa_x V$ are uniformly Lipschitz continuous in $\mu$ under $W_1$ with  Lipschitz constant $C^\mu_1$, where $C^\mu_1>0$ depends only on $d, T, \|\pa_{x}V\|_{L^\infty},\|\pa_{xx}V\|_{L^\infty}$, the  $L^G_1$ in Assumption \ref{assum-regG} (i),  the $L^H(\|\pa_xV\|_{L^\infty})$ in Assumption  \ref{assum-regH}-(i),  and the $\lambda$ in \reff{semimon1}.
\end{prop}
\proof  Denote
\beaa
N_t' := \tilde \dbE_{\mathcal{F}_t}[\pa_{x\mu}  V(X_t, \tilde X_t) \delta \tilde X_t] + {1\over 2}H_{pp}(X_t)^{-1}\tilde \dbE_{\mathcal{F}_t}[ H_{p\mu}(X_t,\tilde X_t) \delta  \tilde X_t].
\eeaa
First, by using \reff{ExmuV} and  \reff{semimon1}, similar to \reff{Nest} we can show
\small\begin{align}\label{NestLL}
 \int_0^t \dbE[|H_{pp}(X_s)^{\frac{1}{2}} N'_s|^2]ds & \le   I(0) +C\dbE\Big[\big(\dbE_{\cF^0_t}[|\d X_t|]\big)^2\Big] +  C\int_0^t \dbE\Big[|\d X_s| ~\tilde \dbE_{\cF_s}[|\d \tilde X_s|]\Big]ds\\
 & = I(0) + C\dbE\Big[\big(\dbE_{\cF^0_t}[|\d X_t|]\big)^2\Big]+  C\int_0^t \dbE\Big[\big(\dbE_{\cF^0_s}[|\d X_s|]\big)^2\Big]ds .\nonumber
\end{align}\normalsize
Next, by \reff{XY}, Young's inequality and noting that $\cF_0$ is independent of $\cF^0_t$, we have
\beaa
\Big(\dbE_{\cF^0_t}[|\delta X_t|]\Big)^2 \leq \Big(\dbE[|\eta|]\Big)^2 + C_{\epsilon}\int_0^t \Big(\dbE_{\cF^0_s}\big[ |\d X_s|\big]\Big)^2ds +\epsilon\int_0^t\Big(\mathbb E_{\mathcal{F}^0_s}\big[ H_{pp}^{\frac{1}{2}}(X_s)|N_s'| \big]\Big)^2ds.
\eeaa
Taking expectation on both sides, choosing $\epsilon>0$ small enough, together with \reff{NestLL} and for the same $\Upsilon$ in \reff{EY1}, it then follows from Gr\"onwall's inequality that
\begin{equation}\label{EY1LL} 
\sup_{t \in [0,T]}\dbE\Big[ \Big(\mathbb E_{\cF^0_t}\big[|\delta X_t|\big]\Big)^2\Big] \leq C  \Big(\mathbb E\big[|\eta|\big]\Big)^2+C |I(0)|\le C  \Big(\mathbb E\big[|\eta|\big]\Big)^2 + C \mathbb E\big[|\eta| |\Upsilon_0|\big].
\end{equation}
Now by \reff{dUpsilon2} and \reff{Upsilon1}, and noting that $|\pa_{x\mu} G|\le L^G_1$, we have
\begin{align*}
|\Upsilon_t| &\le C\tilde \dbE_{\cF_t}\big[|\d \tilde X_T|\big] + C\int_t^T \tilde \dbE_{\cF_t}\big[ |\Upsilon_s| + |\d \tilde X_s|\big] ds\\
&=C\dbE_{\cF^0_t}\big[|\d \tilde X_T|\big] + C\int_t^T \Big[\dbE_{\cF_t}\big[ |\Upsilon_s|\big] +\dbE_{\cF^0_t}\big[ |\d \tilde X_s|\big] \Big]ds.
\end{align*}
Then, since $\cF^0_0$ is degenerate,  namely, it reduced to $\{\emptyset,\Omega_0\}$
\beaa
|\Upsilon_0|\le C\dbE\big[|\d X_T|\big] + C\int_0^T \dbE\big[ |\d  X_s|\big]ds \le C\sup_{0\le t\le T} \dbE\big[|\d X_t|\big].
\eeaa
Combine this with \reff{EY1LL}, we have
\begin{equation}\label{EY2LL} 
\Big| \tilde \dbE_{\cF_0}\Big[\pa_{x\mu} V(0, \xi, \mu, \tilde \xi) \tilde \eta\Big] \Big|^2=|\Upsilon_0|^2 \le C \sup_{t \in [0,T]}\dbE\Big[ \Big(\mathbb E_{\cF^0_t}\big[|\delta X_t|\big]\Big)^2\Big] \leq C  \Big(\mathbb E\big[|\eta|\big]\Big)^2.
\end{equation}
This is the counterpart of \reff{EUpsilon0}. Then, following similar arguments as after \reff{EUpsilon0}, we conclude as follows. First, one obtains
\beaa
\Big|\dbE\big[\pa_{x\mu} V(0, x,  \mu, \xi)  \eta \Big]\Big|\le  C\mathbb E\left[|\eta|\right],\qquad \forall x, \mu, \xi, \eta.
\eeaa
In particular, this implies that there exists a constant $C^\mu_1>0$ such that
\beaa
\Big|\pa_xV(0, x, \cL_{\xi+\eta}) - \pa_x V(0, x, \cL_\xi)\Big| = \Big|\int_0^1 \dbE\big[\pa_{x\mu} V(0, x, \cL_{\xi+\th \eta}, \xi + \th \eta) \eta\big] d\th\Big|\le C^\mu_1 \mathbb E[|\eta|].
\eeaa
Now, taking random variables $\xi,\eta$ such that $W_1(\cL_{\xi+\eta},\cL_{\xi})=\mathbb E[|\eta|]$, the above inequality
implies that $\pa_x V(0,x,\cdot)$ is uniformly Lipschitz continuous in $\mu$ under $W_1$ with uniform Lipschitz constant $C^\mu_1$. 

By analyzing $\bar \Upsilon$ similarly, we show that $V$ is also uniformly Lipschitz continuous in $\mu$ under $W_1$.
\qed

\begin{rem}
(i) Proposition \ref{prop-Lipschitz} indicates that the a priori $W_1$--Lipschitz continuity of $V$ and $\pa_xV$ is a consequence of \eqref{semimon1}, even if $H$ is non-separable. However, we emphasize that, although \eqref{semimon1} is weaker than the Lasry-Lions monotonicity \reff{mon2}, it is stronger than  the displacement semi-monotonicity  condition \reff{semidisplacement1}. Unfortunately, for non-separable $H$, we are not able to find sufficient conditions to ensure \reff{semimon1} a priori for $V$.

(ii) We note that in general for displacement monotone Hamiltonians considered in this manuscript, the displacement semi-monotonicity of the terminal datum is not propagated in time along the solution of the master equation (cf. \cite[Section B.4]{GM}).
\end{rem}

\section{The global well-posedness}
\label{sect-global}
\setcounter{equation}{0}

In this section we establish the global well-posedness of master equation \reff{master}. As illustrated in \cite{CD2,CCD,MZ2}, the key to extend a local classical solution to a global one is the a priori uniform Lipschitz continuity estimate of the solution. We first investigate the regularity of $V$ with respect to $x$. The following result is somewhat standard, while our technical conditions could be slightly different from those in the literature. For completeness we provide a proof in Appendix \ref{sect-appendix}.  We remark that the regularity of $G$ and $H$ in $\mu$ is actually not needed in this result. 
\begin{prop}\label{prop:xbdd}
Let Assumptions \ref{assum-regG}-(i) and \ref{assum-regH}-(i), (iii) hold and $\rho:[0,T]\times\Omega\to\mathcal{P}_2$ be $\mathbb F^0$-progressively measurable (not necessarily a solution to \eqref{FBSDE1}) with 
$\dis\sup_{t\in[0,T]}\mathbb E\big[M_2^2(\rho_t)\big]<+\infty.$ 

(i) For any $x\in\mathbb R^d$ and for the $X^x$ in \eqref{FBSDE2}, the following BSDE on $[t_0, T]$ has a unique solution with bounded $Z^x$:
\begin{equation}\label{FBSDE-x}
Y_t^x=G(X_T^x,\rho_T) - \int_t^TH(X_s^x,\rho_s,Z_s^x)ds-\int_t^TZ_s^x\cdot dB_s-\int_t^TZ_s^{0,x}\cdot dB_s^0.
\end{equation}

(ii) Denote $u(t_0,x):=Y_{t_0}^x$, then there exist $C^x_1, C^x_2>0$, depending only on $d$, $T$, the $C_0$ in  Assumption \ref{assum-regH}-(iii), the constant $L^G_0$ in Assumption \ref{assum-regG}, and the function $L^H$ in Assumption \ref{assum-regH},  such that 
\bea
\label{pauest1}
|\pa_x u(t_0,x)| \le C^x_1,\q  |\pa_{xx} u(t_0,x)| \le C^x_2.
\eea
%
\end{prop}
Here the notation $C^x_i$ denotes the bound of the $i$-th order derivative of $u$ with respect to $x$, in particular, it is {\it not} a function of $x$.

The above result, combined with Theorems \ref{thm-convex} and \ref{thm-Lipschitz}, implies immediately the uniform a priori Lipschitz continuity of $V$ with respect to $\mu$ under $W_2$, with the uniform Lipschitz estimate depending only on the parameters in the assumptions, but not on the additional regularities required in Theorem \ref{thm-convex}. However, the existence of local classical solutions to the master equation \reff{master} requires the Lipschitz continuity under $W_1$, cf.  \cite[Theorem 5.10]{CD2}. To show that eventually the $W_2$--Lipschitz continuity of $V$ together with our standing assumptions on the data imply its Lipschitz continuity under $W_1$, we rely on a pointwise representation formula for $\pa_\mu V$ developed in \cite{MZ2}, tailored to our setting.

 
For this purpose, we fix $t_0\in [0, T]$, $x\in \dbR^d$, $\xi\in \dbL^2(\cF_{t_0})$, and let $\rho$ be given in \eqref{FBSDE1}, provided its wellposedness. We then consider the following FBSDEs  on $[t_0, T]$, which can be interpreted as a formal differentiation of \eqref{FBSDE3} with respect to $x_k$: 
\small\begin{equation}
\label{tdYx}
\left\{
\ba{ll}
\dis \td_{k} X_t^{\xi,x}&=\dis e_k -\int_{t_0}^t \big[(\nabla_k X_s^{\xi,x})^\top\pa_{xp}H(X_s^{\xi,x},\rho_s,Z_s^{\xi,x})
+ (\nabla_k Z_s^{\xi,x})^\top\pa_{pp}H(X_s^{\xi,x},\rho_s,Z_s^{\xi,x})\big]ds; \\[5pt]
\dis \td_{k} Y_t^{\xi,x}&= \dis \pa_x G(X^{\xi,x}_T,\rho_T) \cdot \td_{k} X^{\xi,x}_T-\dis \int_t^T\td_{k} Z_s^{\xi,x}\cdot dB_s^{t_0}- \int_t^T \td_{k} Z_s^{0,\xi,x}\cdot dB_s^{0,t_0}\\[5pt]
\dis\ &+ \dis \int_t^T \pa_x\h L(X_s^{\xi,x},\rho_s,Z_s^{\xi,x})\cdot\td_k X_s^{\xi,x}+ \pa_{p}\h L(X_s^{\xi,x},\rho_s,Z_s^{\xi,x})\cdot\td_k Z_s^{\xi,x}ds,
\ea\right.
\end{equation}\normalsize
the following McKean-Vlasov FBSDE on $[t_0, T]$:  \small
\begin{equation}
\label{tdYx-}
\left\{\ba{ll}
\dis \td_{k} \mathcal X_t^{\xi,x}&=\dis-\int_{t_0}^t\Big\{(\td_k \mathcal X_s^{\xi,x})^{\top}\pa_{xp}H(X_s^{\xi},\rho_s,Z_s^\xi)+(\td_k \mathcal Z_s^{\xi,x})^{\top}\pa_{pp}H(X_s^{\xi},\rho_s,Z_s^\xi)\\[7pt]
\dis & +\tilde{\mathbb E}_{\mathcal{F}_s}\left[(\nabla_{k} \tilde X_s^{\xi,x})^{\top}(\pa_{\mu p} H)(X^{\xi}_s,\rho_s,\tilde X_s^{\xi,x}, Z^{\xi}_s)
\dis +(\nabla_{k} \tilde {\mathcal X}_s^{\xi,x})^{\top}\pa_{\mu p} H(X^{\xi}_s,\rho_s,\tilde X_s^{\xi},Z^{\xi}_s)\right] \Big\}ds; \\[5pt]
\dis \td_{k} \mathcal Y_t^{\xi,x}&= \dis \pa_x G(X^{\xi}_T,\rho_T) \cdot \td_{k} \mathcal X^{\xi,x}_T\\[5pt]
&+\dis\tilde{\mathbb E}_{\mathcal{F}_T}\big[\pa_{\mu}G(X_T^\xi,\rho_T,\tilde X_T^{\xi,x})\cdot \nabla_{k} \tilde X_T^{\xi,x}+\pa_{\mu}G(X_T^\xi,\rho_T, \tilde X_T^\xi)\cdot\nabla_{k}\tilde {\mathcal X}_T^{\xi,x}\big]\\[5pt]
 &\dis + \int_t^T\Big\{ \pa_x \h L\big(X_s^{\xi},\rho_s,Z_s^{\xi})\cdot \td_{k} \mathcal X^{\xi,x}_s+  \pa_p \h L\big(X_s^{\xi},\rho_s,Z_s^{\xi})\cdot \td_{k} \mathcal Z^{\xi,x}_s\\[5pt]
 &\dis +\tilde{\mathbb E}_{\mathcal{F}_s}\big[\pa_{\mu}\h L\big(X^{\xi}_s,\rho_s,\tilde X_s^{\xi,x}, Z^{\xi}_s )\cdot\nabla_{k} \tilde X_s^{\xi,x}+\pa_{\mu}\h L\big(X^{\xi}_s, \rho_s,\tilde X_s^{\xi},Z^{\xi}_s)\cdot\nabla_{k} \tilde {\mathcal X}_s^{\xi,x}\big]\Big\}ds\\[5pt]
  &\dis -\int_t^T\td_{k} \mathcal Z_s^{\xi,x}\cdot dB_s^{t_0}- \int_t^T \td_{k} \mathcal Z_s^{0,\xi,x}\cdot dB_s^{0,t_0},
\ea\right.
\end{equation}\normalsize
and the following McKean-Vlasov BSDE on $[t_0, T]$:
\small
\begin{align}\label{tdYmu}
\left.\ba{lll}
\dis \td_{\mu_k} Y_t^{x,\xi,\tilde x}=\tilde{\mathbb E}_{\mathcal{F}_T}\big[\pa_\mu G (X_T^{x},\rho_T,\tilde X_T^{\xi,\tilde x})\cdot \nabla_{k} \tilde X_T^{\xi,\tilde x}+\pa_\mu G (X_T^{x},\rho_T,\tilde X_T^\xi)\cdot \nabla_{k}\tilde {\mathcal X}_T^{\xi,\tilde x}\big]\\
\dis\qq - \int_t^T \Big\{\pa_pH(X_s^x,\rho_s,Z_s^{x,\xi})\cdot \nabla_{\mu_k} Z_s^{x,\xi,\tilde x}\\
\dis\qq +\tilde{\mathbb E}_{\mathcal{F}_s}\big[\pa_\mu H (X_s^{x},\rho_s,\tilde X_s^{\xi,\tilde x},Z_s^{x,\xi})\cdot\nabla_{k} \tilde X_s^{\xi,\tilde x}+\pa_\mu H (X_s^{x},\rho_s,\tilde X_s^\xi,Z_s^{x,\xi})\cdot \nabla_{k}\tilde {\mathcal X}_s^{\xi,\tilde x}\big]\Big\}ds\\
\dis\qq -\int_t^T\td_{\mu_k} Z_s^{x,\xi,\tilde x}\cdot dB_s- \int_t^T \td_{\mu_k} Z_s^{0,x,\xi,\tilde x}\cdot dB_s^{0}.
\ea\right.
\end{align}\normalsize
The following result 
provides the crucial $W_1$-Lipschitz continuity of $V$. In particular, this extends \cite[Theorem 9.2]{MZ2} to our setting. 

\begin{prop}\label{prop:pamuV}  Let Assumptions \ref{assum-regG}-(i) and \ref{assum-regH}-(i), (iii) hold. Recall the constants $C^x_1$ in \eqref{pauest1}, $L^G_0, L^G_1$ in Assumption \ref{assum-regG},  $L_2^G$  in Remark \ref{rem:W2}, and the function $L^H$ in Assumption \ref{assum-regH}.  Then there exists a constant $\d>0$, depending only $d$, $L^G_0$, $L^G_2$, $L^H(C_1^x)$, such that whenever $T-t_0\le \d$, 
the following hold.

(i) The McKean-Vlasov FBSDEs \eqref{FBSDE1}, \eqref{FBSDE2}, \eqref{FBSDE3}, \eqref{tdYx}, \eqref{tdYx-}, and \reff{tdYmu} are well-posed on $[t_0,T]$, for any $\mu \in \mathcal P_2$ and $\xi \in \mathbb L^2(\mathcal F_{t_0}, \mu)$. 

(ii)
 Define $V(t_0,x, \mu):= Y_{t_0}^{x,\xi}$. We have the pointwise representation:
\bea\label{pamuV}
\pa_{\mu_k} V(t_0,x,\mu,\tilde x)=\nabla_{\mu_k}Y_{t_0}^{x,\xi,\tilde x}.
\eea
Moreover, there exists a constant $C_1^\mu >0$, depending only on $d, L^G_0,  L^G_1, L^H(C_1^x)$ such that
\begin{equation}\label{W1Lip}
|\pa_{\mu}V(t_0,x,\mu,\tilde x)|\leq C_1^\mu,\q |\pa_{x\mu}V(t_0,x,\mu,\tilde x)|\leq C_1^\mu.
\end{equation}

(iii) Assume further that Assumptions \ref{assum-regG}-(ii) and \ref{assum-regH}-(ii) hold true.  Then the master equation \reff{master} has a unique classical solution $V$ on $[t_0, T]$ and \reff{YXV} holds. Moreover, \small
$$
V(t,\cdot,\cdot),\pa_{x}V(t,\cdot,\cdot),\;\; \pa_{xx}V(t,\cdot,\cdot)\in \cC^2(\mathbb R^d\times\mathcal{P}_2),\;\; \pa_\mu V(t,\cdot,\cdot,\cdot),\;\; \pa_{x\mu}V(t,\cdot,\cdot,\cdot)\in \cC^2(\mathbb R^d\times\mathcal{P}_2\times\mathbb R^d),
$$ \normalsize
and all their derivatives in the state and probability measure variables are continuous in the time variable and are uniformly bounded.
\end{prop}
\proof Since this proof is essentially the same as that in \cite[Section 9]{MZ2}, we postpone it to Appendix \ref{sect-appendix}.
\qed

 We emphasize that the $\d$ in the above result depends on $L^G_2$, but not on $L^G_1$, while the $C^\mu_1$ in \reff{W1Lip} depends on $L^G_1$. This observation is crucial. 
We now establish the main result of the paper.

\begin{thm}\label{thm:wellposedness}
Let Assumptions \ref{assum-regG}, \ref{assum-regH}, and  \ref{assum-convex} hold. Then the master equation \eqref{master} on $[0, T]$ admits a unique classical solution $V$ with bounded $\pa_x V$, $\pa_{xx} V$, $\pa_\mu V$, and $\pa_{x\mu} V$. 

Moreover, the McKean-Vlasov FBSDEs \eqref{FBSDE1}, \eqref{FBSDE2}, \eqref{FBSDE3}, \eqref{tdYx}, \eqref{tdYx-}, and \reff{tdYmu} are also well-posed on $[0,T]$  and the representation formula \reff{pamuV} remains true on $[0, T]$.
\end{thm}
\proof  Let $C^x_1, C^x_2$ be as in \reff{pauest1}, and $C^\mu_2$ be the a priori (global) uniform Lipschitz estimate of $V$ with respect to $\mu$ under $W_2$, as established by Theorems \ref{thm-convex} and \ref{thm-Lipschitz}. Let $\d>0$ be the constant in Proposition \ref{prop:pamuV}, but with $L^G_0$ replaced with $C^x_1 \vee C^x_2$ and $L^G_2$ replaced with $C^\mu_2$. Let $0=T_0<\cds<T_n=T$ be a partition such that $T_{i+1}-T_i \le {\d\over 2}$, $i=0,\cds, n-1$. We proceed in three steps.

{\it Step 1. Existence.} First, since $T_n - T_{n-2}\le \d$, by Proposition \ref{prop:pamuV}  the master equation \reff{master}  on $[T_{n-2}, T_n]$ with terminal condition $G$ has a unique classical solution $V$. For each $t\in [T_{n-2}, T_n]$, applying Proposition \ref{prop:xbdd} we have $|\pa_xV(T_{n-1},\cdot,\cdot)|\le C^x_1$, $|\pa_{xx}V(T_{n-1},\cdot,\cdot)|\le C^x_2$. Note that by Proposition \ref{prop:pamuV} (iii) $V(t,\cdot,\cdot)$ has further regularities, this enables us to apply Theorems \ref{thm-convex} and \ref{thm-Lipschitz} and obtain that $V(t,\cdot,\cdot)$ is uniform Lipschitz continuous in $\mu$ under $W_2$ with Lipschitz constant $C^\mu_2$. Moreover, by Proposition \ref{prop:pamuV} (ii) $V(T_{n-1},\cdot,\cdot)$ is also uniformly Lipschitz continuous in $\mu$ under $W_1$.

We next consider  the master equation \reff{master}  on $[T_{n-3}, T_{n-1}]$ with terminal condition $V(T_{n-1},\cdot,\cdot)$. We emphasize that  $V(T_{n-1},\cdot,\cdot)$ has the above uniform regularity with the same constants $C^x_1, C^x_2, C^\mu_2$, then we may apply Proposition \ref{prop:pamuV}  with the same $\d$ and obtain a classical solution $V$ on $[T_{n-3}, T_{n-1}]$ with the additional regularities specified in Proposition \ref{prop:pamuV} (iii). Clearly this extends the classical solution of the master equation to $[T_{n-3}, T_n]$. We emphasize again that, while the bound of $\pa_\mu V(t,\cd)$, $\pa_{x\mu} V(t,\cd)$ may become larger for $t\in [T_{n-3}, T_{n-2}]$ because the $C^\mu_1$ in \reff{W1Lip} now depends on $\|\pa_\mu V(T_{n-1}, \cd)\|_{L^\infty}$ instead of $\|\pa_\mu V(T_{n}, \cd)\|_{L^\infty}$, by the global a priori estimates in Theorems \ref{thm-convex} and \ref{thm-Lipschitz} we see that $V(t,\cd)$ corresponds to the same $C^x_1, C^x_2$ and $C^\mu_2$ for all $t\in [T_{n-3}, T_{n}]$. This enables us to consider the master equation \reff{master}  on $[T_{n-4}, T_{n-2}]$ with terminal condition $V(T_{n-2},\cdot,\cdot)$,  and then we obtain a classical solution on $[T_{n-4}, T_n]$ with the desired uniform estimates and additional regularities. 
 
Now repeat the arguments backwardly in time, we may construct a classical solution $V$ for the original master equation \reff{master} on $[0, T]$ with terminal condition $G$. Moreover, since this procedure is repeated only $n$ times, by applying \reff{W1Lip} repeatedly we see that \reff{W1Lip} indeed holds true on $[0, T]$.

{\it Step 2. Uniqueness.} This follows directly from the local uniqueness in Proposition \ref{prop:pamuV}. Indeed, assume $V'$ is another classical solution with bounded  $\pa_x V'$, $\pa_{xx} V'$, $\pa_\mu V'$, and $\pa_{x\mu} V'$. By otherwise choosing larger $C^x_1, C^x_2, C^\mu_2$, we assume $|\pa_x V'|\le C^x_1$, $|\pa_{xx} V'|\le C^x_2$, and $C^\mu_2$ also serves as a Lipschitz constant for the $W_2$-Lipschitz continuity of $V'$ in $\mu$. Then, applying Proposition \ref{prop:pamuV} on the master equation on $[T_{n-1}, T_n]$ with terminal condition $ G$, by the uniqueness in Proposition \ref{prop:pamuV} (iii) (or in (i)) we see that $V'(t,\cd) = V(t, \cd)$ for $t\in [T_{n-1}, T_n]$. Next consider the master equation on $[T_{n-2}, T_{n-1}]$ with terminal condition $V'(T_{n-1}, \cd) = V(T_{n-1},\cd)$, by the uniqueness in Proposition \ref{prop:pamuV} (iii) again we see that $V'(t,\cd) = V(t, \cd)$ for $t\in [T_{n-2}, T_{n-1}]$. Repeat the arguments backwardly in time we prove the uniqueness on $[0, T]$.

{\it Step 3.} Let $V$ be the unique classical solution to the master equation \eqref{master} on $[0, T]$ with bounded $\pa_x V$ and $\pa_{xx} V$. Then, for $t_0 \in [0,T]$ and $\xi \in \mathbb L^2(\mathcal F_{t_0})$, the McKean-Vlasov SDE \reff{FBSDE4} on $[t_0, T]$ has a unique solution $X^\xi$ and $\rho$. Set  \small
\begin{equation*}
Y^\xi_t: = V(t, X^\xi_t, \rho_t),~ Z^\xi_t := \pa_x V(t, X^\xi_t, \rho_t),~ Z^{0, \xi}_t := \beta \Big( \pa_x V(t, X^\xi_t, \rho_t) +  \tilde {\mathbb E}_{\mathcal F_t} \big[\partial_\mu V\big(t, X^\xi_t, \rho_t, \tilde X^{\xi}_t  \big) \big] \Big).
\end{equation*}\normalsize
By \reff{master} and It\^{o} formula \reff{Ito} one verifies that $(X^\xi, Y^\xi, Z^\xi, Z^{0, \xi})$ satisfies FBSDE \reff{FBSDE1}. The uniqueness follows from the same arguments as in Step 2. 

Similarly, by the above decoupling technique we can easily see that the other McKean-Vlasov FBSDEs  \eqref{FBSDE2}, \eqref{FBSDE3}, \eqref{tdYx}, \eqref{tdYx-}, and \reff{tdYmu}  are also well-posed. In particular, besides \reff{YXV} we have the following:
\begin{align*}
&\dis \td_{k} Y_t^{\xi,x} = \pa_{x_k} V(t,  X_t^{\xi,x}, \rho_t) \td_{k} X_t^{\xi,x},\\
&\dis \td_{k} \mathcal Y_t^{\xi,x}=  \pa_x V(t, X^{\xi}_t,\rho_t) \cdot \td_{k} \mathcal X^{\xi,x}_t\\
&\dis \qq +\tilde{\mathbb E}_{\mathcal{F}_t}\big[\pa_{\mu}V(t,X_t^\xi,\rho_t,\tilde X_t^{\xi,x})\cdot \nabla_{k} \tilde X_t^{\xi,x}+\pa_{\mu}V(X_t^\xi,\rho_t, \tilde X_t^\xi)\cdot\nabla_{k}\tilde {\mathcal X}_t^{\xi,x}\big],\nonumber\\
&\nonumber\dis \td_{\mu_k} Y_t^{x,\xi,\tilde x}=\tilde{\mathbb E}_{\mathcal{F}_t}\big[\pa_\mu V (t, X_t^{x},\rho_t,\tilde X_t^{\xi,\tilde x})\cdot \nabla_{k} \tilde X_t^{\xi,\tilde x}+\pa_\mu V (t, X_t^{x},\rho_t,\tilde X_t^\xi)\cdot \nabla_{k}\tilde {\mathcal X}_t^{\xi,\tilde x}\big].
\end{align*}
\no Set $t=t_0$ and note that $\nabla_{k} \tilde X_{t_0}^{\xi,\tilde x} = e_k$ and $\nabla_{k}\tilde {\mathcal X}_{t_0}^{\xi,\tilde x}=0$, then the last equation above implies 
\beaa
\td_{\mu_k} Y_{t_0}^{x,\xi,\tilde x}=\tilde{\mathbb E}_{\mathcal{F}_{t_0}}\big[\pa_\mu V (t_0, x, \rho_{t_0},\tilde x)\cdot e_k\big] = \pa_{\mu_k} V (t_0, x, \cL_\xi,\tilde x),
\eeaa
which is exactly \reff{pamuV}. 
\qed

\begin{appendix}

\renewcommand {\theequation}{\thesection.\arabic{equation}}
\section{Proofs of some technical results}\label{sect-appendix}
\setcounter{equation}{0}

\proof[Proof of Remark \ref{rem-monequiv}] We show the equivalence of \reff{mon1} and \reff{mon2} for any $U\in \cC^2(\dbR^d\times \cP_2)$.  In fact, by Remark \ref{rem:implicationDisp}(i), we only need to prove that \eqref{mon2} implies \eqref{mon1}.  We now assume \reff{mon2} holds, and we want to show the following which is equivalent to \reff{mon1}: for any $\mu_0, \mu_1 \in \cP_2$,
\bea
\label{mon4}
J_0:= \int_{\dbR^{d}} \big[U(x, \mu_1) - U(x, \mu_0)\big]\big[\mu_1(dx) - \mu_0(dx)\big] \ge 0. 
\eea
Recall \reff{BR}. Since $U$ is continuous, by the standard density argument, it suffices to show \eqref{mon4} for $\mu_i$, $i\in\{0, 1\}$, which have densities $\rho_i \in C^\infty(B_R)$ such that $\min_{B_R} \rho_i >0$.  Consider one of the $W_1$--geodesic interpolations such as in \cite{EvansG}: 
\vskip-0.15in 
$$
 \mu_t:=\rho_t \mathcal L^d, \quad \rho_t= (1-t) \rho_0  + t  \rho_1 \qquad \forall t \in [0, 1].
$$ 
Since $\rho$ is bounded away from $0$ on $[0,1] \times B_R$, then for each $t$, there is a unique solution  $ \phi_t \in H_0^1(B_R^o) \cap C^\infty(B_R)$  to the elliptic equation 
$$\nabla \cdot( \rho_t \nabla \phi_t)=\rho_0-\rho_1  \quad \text{i.e.} \quad \partial_t \rho_t +\nabla \cdot( \rho_t \nabla \phi_t)=0.$$ 
Note that $\phi_t$ and $\nabla \phi_t$ are continuous in $t$ too. Setting $v_t := \nabla \phi_t$, we see that $v$ is a velocity for $t \to \mu_t$, and thus the following chain rule holds, cf. \cite[Lemma 9.8]{GangboS2015}:
\bea
\label{chain}
{d\over dt} U(x, \mu_t) = \int_{\dbR^d} \big\langle \pa_\mu U(x, \mu_t, \tilde x), v_t(\tilde x)\rangle \rho_t(\tilde x) d\tilde x.
\eea

We now compute $J_0$:
\beaa
J_0 &=&\int_{\dbR^{d}} \big[U(x, \mu_1) - U(x, \mu_0)\big]\big[\rho_1(x) - \rho_0(x)\big] dx =-\int_0^1 \int_{\dbR^{d}} {d\over dt} U(x, \mu_t)  \nabla \cdot( \rho_t v_t) dxdt \\
&=&- \int_0^1 \int_{\dbR^{2d}} \big\langle \pa_\mu U(x, \mu_t, \tilde x), v_t(\tilde x)\rangle \rho_t(\tilde x)\big[ \nabla \cdot( \rho_t(x) v_t(x))\big]d\tilde xdxdt. 
\eeaa
By  integration by parts formula, and recalling that $\pa_{x\mu} U := \pa_x [\pa_\mu U]^\top$, we have
\beaa
J_0 =\int_0^1 \int_{\dbR^{2d}}  \big\langle \pa_x\pa_\mu U(x, \mu_t, \tilde x) v_t(\tilde x), v_t(x)\rangle \rho_t(\tilde x) \rho_t(x) d\tilde xdxdt.  
\eeaa
Choose $\xi \in \dbL^2(\cF^1_t, \mu_t)$ and set $\eta := v_t(\xi)$, we see immediately that
\beaa
J_0 = \tilde\dbE\Big[ \big\langle \pa_x\pa_\mu U(\xi, \mu_t, \tilde \xi) \tilde \eta, \eta\big\rangle\Big] \ge 0,
\eeaa
where the inequality is due to \reff{mon2}. This proves \reff{mon4}.
\qed

\medskip

\proof[Proof of Proposition \ref{prop:xbdd}]
We proceed in three steps.

{\it Step 1.} We first show the well-posedness of the BSDE \eqref{FBSDE-x} in the case when $[T-t_0]C_0<1$ where $C_0$ is given in Assumption \ref{assum-regH} (iii). We note that $H$ is only locally Lipschitz continuous. For this purpose,  let $R>0$ be a constant which will be specified later. Let $I_{R}\in C^\infty(\dbR^d)$ be a truncation function such that $I_{R}(p)= p$ for $|p|\le R$,   $|\pa_p I_{R}(p)|=0$ for $|p|\ge R+1$, and $|\pa_p I_{R}(p)|\le 1$ for $p\in \dbR^d$. Denote $H_R(x,\mu,p)= H(x,\mu, I_R(p))$. Then clearly $|\pa_p H_R(x, \mu,p)|\le L^H(R+1)$ and $|\pa_x H_R(x,\mu,p)|\le \tilde L^H(R+1)$ for all $(x,\mu, p)\in \dbR^{d}\times\mathcal{P}_2\times\dbR^d$, where $\tilde L^H(R):= \sup_{(x, \mu,p)\in D_R}|\pa_x H_R(x,\mu,p)|$.  Consider the following BSDE  on $[t_0, T]$ (abusing the notation here):
\bea
\label{FBSDE-x-modified}
\!\!\! Y_t^{x}=G(X_T^{x},\rho_T)-\!\!\int_t^T\!\! H_{R}(X_s^x,\rho_s,Z_s^x)ds-\!\!\int_t^T\!\! Z_s^{x} \!\cd\! dB_s^{t_0}-\!\!\int_t^T\!\! Z_s^{0,x}\!\cd\! dB_s^{0,t_0},
\eea
and denote $u(t_0,x) := Y_{t_0}^{x}$, which is $\cF^{0}_{t_0}$-measurable. By standard BSDE arguments, clearly the above system is well-posed, and it holds $ Y^x_t = u(t,  X^x_t)$, $Z^x_t = \pa_x u(t,  X^x_t)$.  Moreover, we have  $\pa_x u(t_0, x) = \td Y^{x}_{t_0}$, where
\bea
\label{FBSDE-pax-modified}
  \left.\begin{array}{ll} 
\dis \nabla Y_t^{x}=\pa_x G( X_T^{x},\rho_T)-\int_t^T[\pa_x H_{R}( X_s^{x}, \rho_s,Z_s^{x})+\nabla Z_s^{x}~\pa_p H_{R}( X_s^{x}, \rho_s,Z_s^{x})]ds\\
\dis \qquad\qquad\qquad\qquad\quad\, -\int_t^T\nabla Z_s^{x}dB_s^{t_0}-\int_t^T\nabla Z_s^{0,x}dB_s^{0,t_0},\q t_0\le t\le T.
\end{array}
  \right.
\eea
Note that $|\pa_x G|\le L^G_0$ and $|\pa_x H_R|\le \tilde L^H(R+1)$,  one can easily see that
\beaa
|\pa_x u(t_0, x)| = |\td  Y^{x}_{t_0}| \le L^G_0 + T\tilde L^H(R+1).
\eeaa
Note that
\beaa
\limsup_{R\to\infty} {L^G_0 + T\tilde L^H(R+1)\over R} =\limsup_{R\to\infty} { T\tilde L^H(R+1)\over R+1}\leq  TC_0 < 1.
\eeaa
We may choose $R>0$ large enough  such that 
\beaa
|\pa_x u(t_0, x)|  \le L^G_0 + T\tilde L^H(R+1) \le R.
\eeaa
This proves $|\pa_xu(t,x)|\leq C^x_1$ by setting $C^x_1:=R$.   Moreover, since $|Z^x_t|= |\pa_x u(t, X^x_t)|\le R$, we see that $H_R(X^x_t, \rho_t,Z^x_t)= H(X^x_t,\rho_t, Z^x_t)$. Thus $(Y^x, Z^x,Z^{0,x})$ actually satisfies \reff{FBSDE-x}. 

On the other hand, for any solution $(Y^x, Z^x,Z^{0,x})$ with bounded $Z^x$, let $R>0$ be larger than the bound of $Z^x$. Then we see that $(Y^x, Z^x,Z^{0,x})$ satisfies \reff{FBSDE-x-modified}. Now the uniqueness follows from the uniqueness of the BSDE  \reff{FBSDE-x-modified} which has Lipschitz continuous data.

{\it Step 2.} We next estimate  $\pa_{xx} u$, again in the case $[T-t_0]C_0<1$. First, applying standard BSDE estimates on \reff{FBSDE-pax-modified} we see that
\bea
\label{tdZest}
\dbE\Big[ \Big(\int_{t_0}^T |\nabla Z^x_s|^2ds\Big)^2\Big] \le C,\q \mbox{a.s.}
\eea
Then we have $\pa_{xx} u(t_0, x) = \td^2  Y^x_{t_0}$, where, by differentiating \reff{FBSDE-pax-modified} formally in $x$:
\bea
\label{FBSDE-paxx-modified}
  \left.\begin{array}{ll} 
\dis \nabla^2 Y_t^{x}=\pa_{xx} G( X_T^{x},\rho_T)  - \int_t^T\sum_{i=1}^d \big[\nabla^2 Z_s^{i,x}dB_s^{i, t_0} +\nabla^2Z_s^{0,i, x}dB_s^{0,i, t_0}\big] \\
\dis\q -\int_t^T\Big[ \pa_{xx} H_{R}(\cd) + 2\nabla Z_s^{x} \pa_{xp} H_{R}(\cd)+  \nabla  Z_s^{x} \pa_{pp} H_{R}(\cd)  [\nabla Z_s^{x}]^\top \\
\dis \q +  \sum_{i=1}^d \nabla^2 Z_s^{i, x} \pa_{p_i} H_{R}(\cd) \Big](X_s^{x},\rho_s,Z_s^{x})ds.
\end{array}
  \right.
\eea
Denote 
$$M^x_T := \exp\left(-\int_{t_0}^T \pa_{p} H_{R}( X_s^{x},\rho_s,Z_s^{x}) \cd dB_s^{t_0} - {1\over 2} \int_{t_0}^T| \pa_{p} H_{R}(X_s^{x}, \rho_s,Z_s^{x})|^2 ds\right).$$ Then
\begin{align*}
\nabla^2Y_{t_0}^{x} &= \dbE_{\mathcal{F}_{t_0}} \Big[M^x_T \pa_{xx} G( X_T^{x},\rho_T)  - M^x_T\int_{t_0}^T\Big\{\pa_{xx} H_{R}(\cd) + 2\nabla  Z_s^{x} \pa_{xp} H_{R}(\cd)\\
& +  \nabla Z_s^{x} \pa_{pp} H_{R}(\cd)  [\nabla Z_s^{x}]^\top \Big\}( X_s^{x},\rho_s, Z_s^{x})ds\Big].
\end{align*}
Thus, by \reff{tdZest}, there exists some $C^x_2>0$ such that
\beaa
|\pa_{xx} u(t_0, x)| &=& |\nabla^2  Y_{t_0}^{x}| \le C\dbE_{\mathcal{F}_{t_0}} \Big[M^x_T   + M^x_T \int_{t_0}^T\big[1+ |\nabla  Z_s^{x}|^2\big]ds\Big]\\
&\le& C + C\Big( \dbE_{\mathcal{F}_{t_0}} [|M^x_T|^2 ]\Big)^{1\over 2} \Big( \dbE_{\mathcal{F}_{t_0}}\Big[ \big(\int_{t_0}^T |\nabla  Z^x_s|^2ds\big)^2\Big] \Big)^{1\over 2} \le C^x_2.
\eeaa

{\it Step 3.} We now consider the general case. Fix a partition $t_0<\cds<t_n=T$ such that $[t_{i+1} - t_i] C_0 < 1$ for all $i=0,\cds, n-1$.  We proceed backwardly in time by induction. Denote $u(t_n, \cd):= G$. Assume we have defined $u(t_{i+1}, \cd)$ with bounded first and second order derivatives in $x$. Consider the BSDE \reff{FBSDE-x} on $[t_i, t_{i+1}]$ with terminal condition $u(t_{i+1},\cd)$. Applying  the well-posedness result in Step 1 we obtain $u(t_i, x)$ satisfying \reff{pauest1}, for a possibly larger $C^x_1, C^x_2$ which depend on the same parameters. Since $n$ is finite, we obtain \reff{pauest1} for all $i$. Now it follows from standard arguments in FBSDE literature, cf. \cite[Theorem 8.3.4]{Zhang}, that the BSDE \reff{FBSDE-x} on $[t_0, T]$ is wellposed, and \reff{pauest1} holds for all $t\in [t_0, T]$.
\qed

\bs

\proof[Proof of  Proposition  \ref{prop:pamuV}]
 (i)  We first note that, given another initial value $\xi'\in \dbL^2(\cF_{t_0})$ in \reff{FBSDE1}, by \eqref{L2G} the following estimate depends on $L^G_0, L^G_2$, but not on $L^G_1$:
\bea
\label{L2G2}
\left.\ba{c}
\dis \dbE\Big[\big|G(X^\xi_T, \rho^\xi_T) - G(X^{\xi'}_T, \rho^{\xi'}_T)\big|^2\Big] \le 2\dbE\Big[ |L^G_0|^2 |X^\xi_T-X^{\xi'}_T|^2 + |L^G_2|^2 W_2^2(\rho^\xi_T, \rho^{\xi'}_T)\Big]\ss \\
\dis \le 2 \big[ |L^G_0|^2+|L^G_2|^2\big]\dbE\Big[|X^\xi_T-X^{\xi'}_T|^2 \Big]. 
\ea\right.
\eea
By first replacing $H$ with $H_R$ as in the proof of Proposition \ref{prop:xbdd}, 
it follows from the standard contraction mapping argument in FBSDE literature, cf. \cite[Theorem 8.2.1]{Zhang}, there exists $\d=\d_R>0$ such that the McKean-Vlasov FBSDE \eqref{FBSDE1} with $H_R$ is well-posed whenever $T-t_0\le \d$. Noticing that in the contraction mapping argument $G$ is used exactly in the form of \reff{L2G2}, here $\d_R$ depends on $d$, $L^G_0$,  $L^G_2$,  the function $L^H$, and $R$, but not on $L^G_1$. By Proposition \ref{prop:xbdd}, we can see that $|Z^\xi|\le C^x_1$, for the $C^x_1$ in \reff{pauest1} which does not depend on $R$. Now set $R = C^x_1$ and hence $\d$ depends only on $d, L^G_0, L^G_2, L^H_2(C^x_1)$, we see that $H_R(X^\xi_s, \rho_s, Z^\xi_s) = H(X^\xi_s, \rho_s, Z^\xi_s)$ and thus \reff{FBSDE1} is well-posed, which includes existence, uniqueness, and in particular stability. Similarly in \eqref{FBSDE2}, \eqref{FBSDE3}, \eqref{tdYx}, \eqref{tdYx-}, and \reff{tdYmu}, the difference of the terminal condition also appears like \reff{L2G2}, and thus they are also well-posed when $T-t_0\le \d$. In particular, we point out that in  the terminal condition of \reff{tdYx-}: 
\beaa
\tilde{\mathbb E}_{\mathcal{F}_T}\big[\pa_{\mu}G(X_T^\xi,\rho_T,\tilde X_T^{\xi,x})\cdot \nabla_{k} \tilde X_T^{\xi,x}+\pa_{\mu}G(X_T^\xi,\rho_T,X_T^\xi)\cdot\nabla_{k}\tilde {\mathcal X}_T^{\xi,x}\big],
\eeaa
 only $\nabla_{k}\tilde {\mathcal X}_T^{\xi,x}$ is part of the solution while all other involved random variables are already obtained from the other FBSDEs. Then in the contraction mapping argument,  the random coefficient $\pa_{\mu}G(X_T^\xi,\rho_T,X_T^\xi)$ of the solution term $\nabla_{k}\tilde {\mathcal X}_T^{\xi,x}$ again appears in $\dbL^2$-sense: 
\beaa
\tilde \dbE\Big[\Big|\pa_{\mu}G(X_T^\xi,\rho_T, \tilde X_T^\xi)\cdot \tilde \zeta_1 - \pa_{\mu}G(X_T^\xi,\rho_T,\tilde X_T^\xi)\cdot \tilde \zeta_2\Big|^2\Big] \le |L^G_2|^2 \dbE\big[|\zeta_1-\zeta_2|^2\big].
\eeaa


(ii) The proof of this point is rather lengthy, but follows almost the same arguments as in \cite[Theorem 9.2]{MZ2}. Given \reff{pamuV}, the first estimate of \reff{W1Lip} follows directly from the estimate for BSDE \reff{tdYmu}. Indeed, by Proposition \ref{prop:xbdd}, we obtained $|Z^\xi|,|Z^{x,\xi}|,|Z^{\xi,x}|\leq C_1^x$. Then, using Assumptions \ref{assum-regG}-(i) and \ref{assum-regH}-(i), the coefficients in the linear FBSDEs \eqref{tdYx} and \eqref{tdYx-}, and the linear BSDE \eqref{tdYmu} are bounded by $L_0^G$, $L_1^G$, $L^H(C_1^x)$ correspondingly. Therefore, there exists $C_1^{\mu}>0$ depending on these constants and $d$ such that 
\[
|\pa_{\mu_k}V(t_0,x,\mu,\tilde x)|=|\nabla_{\mu_k}Y_{t_0}^{x,\xi,\tilde x}|\leq C_1^{\mu}.
\]

Moreover, by differentiating \reff{tdYmu} with respect to $x$, we can derive the representation formula for $\pa_{x\mu}V$ from \reff{pamuV} and then the second estimate of \reff{W1Lip} also follows directly from the estimate for the differentiated BSDE.

We now prove \reff{pamuV} in four steps. Without loss of generality we prove only the case that $k=1$ and $t_0=0$. Throughout the proof, it is sometimes convenient to use the notation $X^{x,\xi}:=X^x$. 

{\it Step 1.} For any  $\xi\in \dbL^2(\cF_0, \mu)$ and any scalar random variable $\eta\in \dbL^2(\cF_0, \dbR)$, following standard arguments and by  the stability property of the involved systems we have
\bea
\label{tdXconv}
 \lim_{\e\to 0} \dbE\Big[\sup_{0\le t\le T} \Big| {1\over \e}\big[X^{\xi + \e \eta e_1}_t - X^\xi_t\big] - \delta X^{\xi, \eta e_1}_t\Big|^2\Big] = 0,
\eea
where $\Big(\delta  X^{\xi, \eta e_1}, \delta  Y^{\xi, \eta e_1}, \delta  Z^{\xi, \eta e_1},\delta  Z^{0,\xi, \eta e_1}\Big)$ satisfies the linear McKean-Vlasov FBSDE:
\small
\begin{equation}\label{tdmuFBSDE}
\left\{\ba{ll}
&\dis  \delta  X^{\xi, \eta e_1}_t=\eta e_1 - \int_0^t ( \delta  X_s^{\xi,\eta e_1})^{\top}\pa_{xp} H\big(X^{\xi}_s,\rho_s,Z^{\xi}_s)+( \delta  Z_s^{\xi,\eta e_1})^{\top}\pa_{pp} H\big(X^{\xi}_s,\rho_s,Z^{\xi}_s)\\
&\dis +\tilde \dbE_{\cF_s}\big[\pa_{p\mu}H(X_s^\xi,\rho_s, \tilde X_s^\xi, Z_s^\xi)\cdot \delta \tilde X_s^{\xi,\eta e_1}\big]ds; \\
&\dis \delta  Y_t^{\xi, \eta e_1}= \pa_x G(X_T^{\xi},\rho_T) \cdot \delta  X^{\xi, \eta e_1}_T + \tilde \dbE_{\cF_T}\big[\pa_\mu G(X_T^{\xi},\rho_T, \tilde X^\xi_T)\cdot \delta  \tilde X^{\xi, \eta e_1}_T\big]\\
 & \dis +\int_t^T  \pa_x \h L\big(X^{\xi}_s, \rho_s,Z^{\xi}_s)\cdot \delta  X^{\xi,\eta e_1}_s+   \pa_p \h L\big(X^{\xi}_s,\rho_s, Z^{\xi}_s)\cdot \delta  Z^{\xi,\eta e_1}_s\\
&\dis +\tilde \dbE_{\cF_s}\big[\pa_\mu \h L(X_s^{\xi},\rho_s, \tilde X^\xi_s,Z_s^{\xi})\cdot\delta  \tilde X^{\xi, \eta e_1}_s\big]ds - \int_t^T\delta  Z_s^{\xi, \eta e_1}\cdot dB_s- \int_t^T\delta  Z_s^{0,\xi,\eta e_1}\cdot d B_s^0.
\ea\right.
\end{equation}\normalsize
Similarly, by \reff{tdXconv} and \reff{FBSDE2}, one can show that
\bea
\label{tdYconv}
 \lim_{\e\to 0} \dbE\Big[\sup_{0\le t\le T} \Big| {1\over \e}\big[Y^{x, \xi + \e \eta e_1}_t - Y^{x,\xi}_t\big] - \delta  Y^{x, \xi, \eta e_1}_t\Big|^2\Big] = 0,
\eea
where $\Big(\delta  Y^{x,\xi, \eta e_1}, \delta  Z^{x, \xi, \eta e_1}, \delta  Z^{0,x, \xi, \eta e_1}\Big)$ satisfies the linear (standard) BSDE:\small
\begin{align}\label{eq:nabla mu X}
&\dis \delta  Y_t^{x, \xi, \eta e_1}=  \tilde \dbE_{\cF^0_T}\big[\pa_\mu G(X_T^x,\rho_T, \tilde X^\xi_T) \cdot\delta  \tilde X^{\xi, \eta e_1}_T\big] \nonumber\\[7pt]
&\dis\q -\int_t^T \pa_p H(X_s^{x,\xi},\rho_s,Z_s^{x,\xi})\cdot \delta  Z^{x,\xi,\eta e_1}_s+\tilde \dbE_{\cF_s}\big[\pa_\mu H(X_s^{x,\xi},\rho_s,\tilde X_s^{\xi},Z_s^{x,\xi})\cdot \delta  \tilde X_s^{\xi,\eta e_1}\big]ds\\
&\dis \q -\int_t^T\delta  Z_s^{x,\xi, \eta e_1}\cdot dB_s-\int_t^T\delta  Z_s^{0,x,\xi, \eta e_1}\cdot dB_s^0. \nonumber
\end{align}\normalsize
In particular, \reff{tdYconv} implies,
\bea
\label{pamuVconv}
 \lim_{\e\to 0}  \Big| {1\over \e}\big[V(0,x, \cL_{\xi + \e \eta e_1}) - V(0,x,\cL_\xi)\big] - \delta  Y^{x, \xi, \eta e_1}_0\Big|^2\ = 0.
\eea
Thus, by the definition of $\pa_\mu V$,
\bea
\label{pamuV1}
\dbE\big[\pa_{\mu_1} V(0,x, \mu, \xi) \eta\big] =  \delta  Y^{x, \xi, \eta e_1}_0.
\eea

{\it Step 2.} In this step we assume $\xi$ (or say, $\mu$) is discrete: $p_i = \dbP(\xi = x_i)$, $i=1,\cds, n$. Fix $i$, consider the following system of McKean-Vlasov FBSDEs: for $j=1,\cds, n$,
\begin{equation}\label{tdmuYij}
\left\{\ba{ll}
\dis  \td_{\mu_1} X^{i,j}_t&\dis= \1_{\{i=j\}} e_1 - \int_0^t \sum_{k=1}^n p_k \tilde \dbE_{\cF_s}\Big[ (\td_{\mu_1} \tilde X^{i,k}_s)^{\top}\pa_{\mu p} H(X_s^{\xi,x_j},\rho_s, \tilde X^{\xi,x_k}_T,Z_s^{\xi,x_j}) \Big]\\[7pt]
\dis & \dis+(\td_{\mu_1} X_s^{i,j})^{\top}\pa_{xp} H\big(X^{\xi,x_j}_s,\rho_s,Z^{\xi,x_j}_s)+(\td_{\mu_1} Z_s^{i,j} )^{\top}\pa_{pp} H\big(X^{\xi,x_j}_s,\rho_s,Z^{\xi,x_j}_s)ds; \\[7pt]
\dis \td_{\mu_1} Y_t^{i,j}& = \dis\pa_x G(X_T^{\xi,x_j},\rho_T) \cdot\td_{\mu_1} X^{i,j}_T+\sum_{k=1}^n p_k \tilde \dbE_{\cF_T}\Big[\pa_\mu G(X_T^{\xi,x_j},\rho_T, \tilde X^{\xi,x_k}_T)\cdot \td_{\mu_1} \tilde X^{i,k}_T\Big] \\[7pt]
\dis &\dis + \int_t^T  \pa_x \h L\big(X^{\xi,x_j}_s,\rho_s,Z^{\xi,x_j}_s)\cdot\td_{\mu_1} X^{i,j}_s+\pa_p \h L\big(X^{\xi,x_j}_s,\rho_s,Z^{\xi,x_j}_s)\cdot\td_{\mu_1} Z^{i,j}_s\\[7pt]
\dis&\dis +\sum_{k=1}^n p_k \tilde \dbE_{\cF_s}\Big[\pa_\mu \h L(X_s^{\xi,x_j},\rho_s, \tilde X^{\xi,x_k}_s,Z_s^{\xi,x_j}) \cdot\td_{\mu_1} \tilde X^{i,k}_sds\\[7pt]
\dis& \dis - \int_t^T\td_{\mu_1} Z_s^{i,j}\cdot dB_s-\b \int_t^T\td_{\mu_1} Z_s^{0,i,j}\cdot d B_s^0.
\ea\right.
\end{equation}
For any $\Phi\in\{X, Y, Z, Z^0\}$, we define
\beaa
\td_{1} \Phi^{\xi, x_i} := \td_{\mu_1} \Phi^{i,i},\q  \td_{1} \Phi^{\xi, x_i, *}:= {1\over p_i} \sum_{j\neq i} \td_{\mu_1} \Phi^{i,j} \1_{\{\xi=x_j\}}.
\eeaa
 Note that $\Phi^\xi =  \sum_{j=1}^n \Phi^{\xi, x_j} \1_{\{\xi=x_j\}}$. Since \reff{tdmuYij} is linear, one can easily check that
\footnotesize
\bea\label{tdXxi-}
&&\dis  \td_1 X^{\xi, x_i}_t=e_1- \int_0^t\Big\{(\td_1 X_s^{\xi,x_i})^{\top} \pa_{xp} H\big(X^{\xi, x_i}_s,\rho_s,Z^{\xi, x_i}_s)+(\td_1 Z_s^{\xi,x_i})^{\top} \pa_{pp} H\big(X^{\xi, x_i}_s,\rho_s,Z^{\xi, x_i}_s)\nonumber\\
&&\dis\qq +p_i\tilde \dbE_{\cF_s}\Big[(\td_{1}\tilde X_s^{\xi,x_i})^{\top}\pa_{\mu p}H(X_s^{\xi,x_i},\rho_s,\tilde X_s^{\xi,x_i},Z_s^{\xi,x_i})\\
&&\qq +(\td_1\tilde X_s^{\xi,x_i, *})^{\top}\pa_{\mu p}H(X_s^{\xi,x_i},\rho_s,\tilde X_s^{\xi},Z_s^{\xi,x_i})\Big]\Big\}ds,\nonumber\\
 &&\dis  \td_1 X^{\xi, x_i, *}_t=- \int_0^t\Big\{ (\td_1 X_s^{\xi,x_i, *})^{\top} \pa_{xp} H\big(X^{\xi}_s,\rho_s,Z^{\xi}_s)+(\td_1 Z_s^{\xi,x_i, *})^{\top}\pa_{pp} H\big(X^{\xi}_s,\rho_s,Z^{\xi}_s)\nonumber\\
 &&\dis \qq+\tilde \dbE_{\cF_s}\Big[(\td_{1}\tilde X_s^{\xi,x_i})^{\top}\pa_{\mu p}H(X_s^{\xi},\rho_s,\tilde X_s^{\xi,x_i},Z_s^{\xi})\label{tdXxi-new} \\
&&\qq +(\td_1\tilde X_s^{\xi,x_i, *})^{\top}\pa_{\mu p}H(X_s^{\xi},\rho_s,\tilde X_s^{\xi},Z_s^{\xi})\Big]\1_{\{\xi\neq x_i\}}\Big\}ds,\nonumber\\
\label{tdXxi-2}
&&\dis \td_1 Y_t^{\xi, x_i}= \pa_x G(X_T^{\xi,x_i},\rho_T) \cdot \td_1 X^{\xi, x_i}_T - \int_t^T\td_1 Z_s^{\xi, x_i}dB_s- \int_t^T\td_1 Z_s^{0,\xi,x_i}d B_s^0\nonumber\\
&&\dis\qq + p_i\tilde \dbE_{\cF_T}\big[\pa_\mu G(X_T^{\xi,x_i},\rho_T, \tilde X^{\xi,x_i}_T)\cdot \td_1\tilde X_s^{\xi,x_i}+\pa_\mu G(X_T^{\xi,x_i},\rho_T, \tilde X^\xi_T)\cdot \td_1 \tilde X^{\xi, x_i, *}_T\big]\\
&&\dis\qq + \int_t^T \Big\{ \pa_x \h L\big(X^{\xi,x_i}_s,\rho_s,Z^{\xi,x_i}_s)\cdot \td_1 X^{\xi,x_i}_s+ \pa_p \h L\big(X^{\xi,x_i}_s,\rho_s,Z^{\xi,x_i}_s)\cdot \td_1 Z^{\xi,x_i}_s + p_i \times\nonumber\\
&&\dis\qq\q \tilde \dbE_{\cF_s}\big[\pa_\mu \h L(X_s^{\xi,x_i},\rho_s, \tilde X^{\xi,x_i}_s,Z_s^{\xi,x_i})\cdot \td_1\tilde X_s^{\xi,x_i} 
+ \pa_\mu \h L(X_s^{\xi,x_i},\rho_s,  \tilde X^\xi_s,Z_s^{\xi,x_i})\cdot \td_1 \tilde X^{\xi, x_i, *}_T)\big]\Big\}ds, \nonumber\\
\label{tdXxi-2new}
&&\dis \td_1 Y_t^{\xi, x_i, *}= \pa_x G(X_T^{\xi},\rho_T)\cdot \td_1 X^{\xi, x_i, *}_T - \int_t^T\td_1 Z_s^{\xi, x_i, *}\cdot dB_s- \int_t^T\td_1 Z_s^{0,\xi,x_i, *}\cdot d B_s^0
\nonumber\\
&&\dis\qq +\tilde \dbE_{\cF_T}\big[\pa_\mu G(X_T^{\xi},\rho_T, \tilde X^{\xi,x_i}_T)\cdot \td_1\tilde X_s^{\xi,x_i} +\pa_\mu G(X_T^{\xi},\rho_T, \tilde X^\xi_T)\cdot \td_1 \tilde X^{\xi, x_i, *}_T\big]\1_{\{\xi\neq x_i\}}
\\
&&\dis \qq+ \int_t^T\bigg\{  \pa_x \h L\big(X^{\xi}_s,\rho_s,Z^{\xi}_s)\cdot \td_1 X^{\xi,x_i, *}_s+ \pa_p \h L\big(X^{\xi}_s,\rho_s,Z^{\xi}_s)\cdot \td_1 Z^{\xi,x_i, *}_s]\nonumber\\
&&\dis\qq +\tilde \dbE_{\cF_s}\big[\pa_\mu \h L(X_s^{\xi},\rho_s, \tilde X^{\xi,x_i}_s,Z_s^{\xi})\cdot \td_1\tilde X_s^{\xi,x_i}
\ +\pa_\mu \h L(X_s^{\xi},\rho_s,  \tilde X^\xi_s,Z_s^{\xi})\cdot \td_1 \tilde X^{\xi, x_i-}_T)\big]\1_{\{\xi\neq x_i\}}\bigg\}ds.\nonumber
\eea
\normalsize
Since \reff{tdmuFBSDE} is also linear, one can easily check that, for $\Phi\in\{X, Y, Z, Z^0\}$,
\bea
\label{tdXdecom}
\delta  \Phi^{\xi, \1_{\{\xi=x_i\}}e_1} = \td_1 \Phi^{\xi, x_i} \1_{\{\xi=x_i\}}+ p_i \td_1   \Phi^{\xi, x_i, *}.
\eea
Moreover, note that
\begin{align*}
\dis & \tilde \dbE_{\cF_T}\big[\pa_\mu G(X_T^{x,\xi},\rho_T, \tilde X^\xi_T)\cdot\delta \tilde X^{\xi, \1_{\{\xi=x_i\}}e_1}_T\big]\ms\\
\dis &=   \tilde \dbE_{\cF_T}\Big[\pa_\mu G(X_T^{x,\xi},\rho_T, \tilde X^\xi_T)\cdot\big[ \td_1 \tilde X_T^{\xi, x_i} \1_{\{\tilde \xi=x_i\}}+ p_i \td_1   \tilde X_T^{\xi, x_i, *}\big]\Big] \ms\\
\dis&= p_i \tilde \dbE_{\cF_T}\Big[\pa_\mu G(X_T^{x,\xi},\rho_T, \tilde X^{\xi, x_i}_T)\cdot \td_1 \tilde X^{\xi, x_i}_T+\pa_\mu G(X_T^{x,\xi},\rho_T, \tilde X^\xi_T)\cdot\td_1 \tilde X^{\xi, x_i, *}_T\Big]
\end{align*}
and similarly
\begin{align*}
\dis&  \tilde \dbE_{\cF_s}\big[\pa_\mu H(X_s^{x,\xi},\rho_s, \tilde X^\xi_s,Z_s^{x,\xi})\cdot\delta \tilde X^{\xi, \1_{\{\xi=x_i\}}e_1}_s\big]\ms\\
\dis&= p_i \tilde \dbE_{\cF_s}\Big[\pa_\mu H(X_s^{x,\xi},\rho_s, \tilde X^{\xi, x_i}_s,Z_s^{x,\xi})\cdot\td_1 \tilde X^{\xi, x_i}_s+\pa_\mu H(X_s^{x,\xi},\rho_s,\tilde X^\xi_s,Z_s^{x,\xi})\cdot \td_1 \tilde X^{\xi, x_i, *}_s)\Big].
\end{align*}
Plug this into \reff{eq:nabla mu X}, we obtain 
\bea
\label{tdYxipi}
\delta \Phi_t^{x, \xi, \1_{\{\xi=x_i\}}e_1}=  p_i \td_{\mu_1}\Phi^{x,\xi, x_i}_t,
\eea
 where
 \small
\bea
\label{eq:nabla mu X2}
&&\dis  \td_{\mu_1} Y^{x,\xi, x_i}_t= \tilde \dbE_{\cF_T}\big[\pa_\mu G(X_T^x,\rho_T, \tilde X^{\xi,x_i}_T)\cdot  \td_1 \tilde X^{\xi, x_i}_T+\pa_\mu G(X_T^x,\rho_T, \tilde X^\xi_T)\cdot  \td_1 \tilde X^{\xi, x_i, *}_T)\big] \nonumber \\
&&\dis\qq\q -\int_t^T \bigg\{\pa_p H(X_s^{x,\xi},\rho_s,Z_s^{x,\xi})\cdot \td_{\mu_1} Z^{x,\xi,x_i}_s\\
&&\dis\qq\q+ \tilde \dbE_{\cF_s}\Big[\pa_\mu H(X_s^{x,\xi},\rho_s, \tilde X^{\xi, x_i}_s,Z_s^{x,\xi})\cdot\td_1 \tilde X^{\xi, x_i}_s+\pa_\mu H(X_s^{x,\xi},\rho_s,\tilde X^\xi_s,Z_s^{x,\xi})\cdot \td_1 \tilde X^{\xi, x_i, *}_s\Big]\bigg\}ds\nonumber\\
&&\dis\qq\q-\int_t^T\td_{\mu_1} Z_s^{x,\xi, x_i}\cdot dB_s-\int_t^T\td_{\mu_1} Z_s^{0,x,\xi, x_i}\cdot dB_s^0. \nonumber
\eea\normalsize
In particular, by setting $\eta = \1_{\{\xi=x_i\}}$ in \reff{pamuV1} we obtain:
\bea
\label{pamuVdiscrete}
\pa_{\mu_1} V(0, x, \mu, x_i) = \td_{\mu_1} Y^{x,\xi, x_i}_0.
\eea
We shall note that \reff{tdXxi-}-\reff{tdXxi-new}, \eqref{tdXxi-2}-\eqref{tdXxi-2new} is different from \reff{tdYx} and \reff{tdYx-}, so \reff{pamuVdiscrete} provides an alternative discrete representation. 

{\it Step 3.} We now prove \reff{pamuV} in the case that $\mu$ is absolutely continuous. For each $n\ge 3$, set 
\[
x^n_{\vec i} := \frac{\vec i}{n}, \quad
\D_{\vec i}^n:=\left[\frac{i_1}{n},\frac{i_1+1}{n}\right)\times\cdots\times\left[\frac{i_d}{n},\frac{i_d+1}{n}\right),\q \vec i=(i_1,\cdots,i_d)^{\top}\in\mathbb Z^d,
\]
For any $x\in\mathbb R^d$, there exists $\vec i(x):=(i_1(x),\cdots,i_d(x))\in \mathbb Z^d$ such that $x\in \D_{\vec i(x)}^n$. Let 
$$
\vec i^n(x):=(i^n_1(x),\cdots,i^n_d(x))\in \mathbb Z^d,\q\mbox{where}\q 
i^n_l(x):=\min\{\max\{i_l,-n^2\},n^2\}, ~ l=1,\cdots,d.
$$
Denote $Q_n:= \{x\in \dbR^d: |x_i|\le n, i=1,\cds, d\}$, $\mathbb Z_n^d:=\{\vec i\in\mathbb Z^d\,:\,\D^n_{\vec{i}}\cap Q_{n}\not=\emptyset\}$, and 
\bea
\label{xin}
\xi_n :=\sum_{\vec i\in \mathbb Z_n^d} x_{\vec i}^n\1_{\D_{\vec i}^n}(\xi)+\frac{\vec i^n(\xi)}{n}\1_{Q_n^c}(\xi). 
\eea
 It is clear that $\lim_{n\to+\infty}\mathbb E\big[|\xi_n - \xi|^2\big]=0$ and thus $\lim_{n\to\infty} W_2(\cL_{\xi_n}, \cL_\xi)=0$. Then for any  scalar random variable $\eta$, by stability of FBSDE \reff{tdmuFBSDE} and  BSDE \reff{eq:nabla mu X}, we derive from \reff{pamuV1} that
 \bea
 \label{tdmuYconv}
 \dbE\Big[\pa_{\mu_1} V(0,x,\mu, \xi)\eta \Big] = \delta Y^{x,\xi, \eta e_1}_0 = \lim_{n\to\infty} \delta Y^{x,\xi_n, \eta e_1}_0.
  \eea

For each $\tilde x\in \dbR^d$, let $\vec i(\tilde x)$ be the $i$ such that $\tilde x \in \D_{\vec i}^n$, which holds when $n> |\tilde x|$. Then $\big(\cL_{\xi_n}, \frac{\vec i(\tilde x)}{n}\big) \to (\mu, \tilde x)$ as $n\to \infty$. By the stability of FBSDEs \reff{FBSDE1}-\reff{FBSDE2}, we have 
$
\Big(X^{\xi_n, \frac{\vec{i}(\tilde x)}{n}}, Z^{\xi_n, \frac{\vec{i}(\tilde x)}{n}}\Big) \to (X^{\xi, \tilde x}, Z^{\xi, \tilde x})
$ 
under appropriate norms. Moreover, since $\mu$ is absolutely continuous, 
$$
\dbP\Big(\xi_n = \frac{\vec{i}(\tilde x)}{n}\Big) = \dbP\Big(\xi \in \D_{\vec i}^n\Big) \to 0, \quad \text{as} \quad n\to \infty.
$$ 
Then by the stability of \reff{tdXxi-}-\reff{tdXxi-new}, \eqref{tdXxi-2}-\eqref{tdXxi-2new} and \reff{eq:nabla mu X2} we can check that
\small \begin{equation}
 \label{xinconv}
 \lim_{n\to\infty}\bigg(\td_1 \Phi^{\xi_n,\frac{\vec{i}(\tilde x)}{n}},\; \td_1 \Phi^{\xi_n,\frac{\vec{i}(\tilde x)}{n}, *},\; \td_{\mu_1}\Phi^{x,\xi_n,\frac{\vec{i}(\tilde x)}{n}}\bigg) =\bigg(\td_1 \Phi^{\xi,\tilde x},\; \td_1\Phi^{\xi,\tilde x, *},\; \td_{\mu_1}\Phi^{x,\xi,\tilde x}\bigg).
  \end{equation}\normalsize
 Now for any bounded function $\f\in C(\dbR^d)$, set $\eta=\f(\xi)$ in \reff{tdmuYconv}, by  \reff{tdYxipi} we have
\beaa
 \dbE\Big[\pa_{\mu_1} V(0,x,\mu, \xi)\f(\xi)\Big] 
 =  \lim_{n\to\infty} \delta Y^{x,\xi_n, \f(\xi_n)e_1}_0 
 =  \lim_{n\to\infty} \sum_{\vec i\in \mathbb Z_n^d} \f\Big(x^n_{\vec i}\Big)\delta Y^{x,\xi_n, \1_{\{\xi_n=x^n_{\vec i}\}}e_1}_0
 \eeaa
and so, 
 \begin{align*}
\dbE\Big[\pa_{\mu_1} V(0,x,\mu, \xi)\f(\xi)\Big]&= \lim_{n\to\infty} \sum_{\vec i\in\mathbb Z^d_n} \f(x^n_{\vec i})\td_{\mu_1} Y^{x,\xi_n, x^n_{\vec i}}_0 \dbP(\xi\in \D_{\vec i})\\
&  = \int_{\mathbb R^d} \f(\tilde x)\td_{\mu_1} Y^{x,\xi, \tilde x}_0 \mu(d\tilde x).
 \end{align*}
 This implies \reff{pamuV} immediately.

 {\it Step 4.} We finally prove the general case. Denote $\psi(x,\mu, \tilde x):= \td_{\mu_1} Y^{x,\xi, \tilde x}_0$. By the stability of FBSDEs, $\psi$ is continuous in all the variables. Fix an arbitrary $(\mu, \xi)$. One can easily construct $\xi_n$ such that  $\cL_{\xi_n} $ is absolutely continuous and $\lim_{n\to\infty}\dbE[|\xi_n-\xi|^2]=0$. Then, for any $\eta= \f(\xi)$ as in Step 3, by \reff{pamuV1} and Step 3 we have 
 \begin{align*}
 \dis \dbE\big[\pa_{\mu_1} V(0,x, \mu, \xi) \f(\xi)\big]& =  \lim_{n\to \infty}  \delta Y^{x, \xi_n, \f(\xi_n)e_1}_0=
  \lim_{n\to \infty}\dbE\big[ \psi(x, \cL_{\xi_n}, \xi_n)\f(\xi_n)\big]\\
  & = \dbE\big[ \psi(x, \mu, \xi)\f(\xi)\big], 
 \end{align*}
which implies \reff{pamuV} in the general case. 

\ms
(iii) This result follows immediately from well-known facts. Indeed, given the uniform estimate of $\pa_\mu V$ in \reff{W1Lip}, the well-posedness of the master equation \reff{master} on $[t_0, T]$ follows from the arguments in \cite[Theorem 5.10]{CD2}. This, together with It\^{o} formula \reff{Ito}, will easily lead to \reff{YXV}.  Under the additional Assumptions \ref{assum-regG}-(ii) and \ref{assum-regH}-(ii), the representation formulas and the boundedness of higher order derivatives in state and probability variables can be proved by further differentiating the McKean-Vlasov FBSDEs \eqref{FBSDE1}-\eqref{FBSDE3}, \eqref{tdYx}-\eqref{tdYx-} and BSDE \eqref{tdYmu}, with respect to the state and probability variables. The calculation is lengthy but very similar to that in   \cite[Section 9.2]{MZ2}. We omit the details.
\qed

\end{appendix}

{\bf Acknowledgments.}
The research of WG was supported by NSF grant DMS--1700202 and Air Force grant FA9550-18-1-0502. ARM was partially supported by the King Abdullah University of Science and Technology Research Funding (KRF) under Award No. ORA-2021-CRG10-4674.2. CM gratefully acknowledges the support by CityU Start-up Grant 7200684 and Hong Kong RGC Grant ECS 9048215. The research of JZ was supported in part by NSF grant DMS-1908665. We thank the anonymous referees for their thoughtful  comments which helped to improve our manuscript greatly.

\end{document}